\documentclass[11pt]{amsart}

\usepackage{amsmath,amsthm, amscd, amssymb, amsfonts}
\usepackage[all]{xy}
\usepackage{verbatim}
\usepackage{color}
\usepackage{array}
\usepackage{multirow}
\usepackage{multicol}

\usepackage{hhline}
\usepackage{enumitem}
\usepackage{float}
\usepackage{hyperref}
\allowdisplaybreaks
\usepackage{mathtools} 

\usepackage{listings}

\newtheorem{step}{Step}
\newtheorem{stepcyc}{Step}
\newtheorem{theorem}{Theorem}[section]
\newtheorem{lemma}[theorem]{Lemma}

\newtheorem{conjecture}[theorem]{Conjecture}
\newtheorem{proposition}[theorem]{Proposition}

\newtheorem{criterium}{Criterium}

\theoremstyle{definition}

\newtheorem*{notation}{Notation}

\theoremstyle{remark}
\newtheorem{remark}[theorem]{Remark}

\def\pf{\begin{proof}}
\def\epf{\end{proof}}

\newcommand{\nc}{\newcommand}

\newcommand{\I}{{\mathbb I}}

\newcommand{\qti}[1]{\widetilde{q_{#1}}}
\newcommand{\rti}[1]{\widetilde{r_{#1}}}

\nc{\ub}{\mathfrak{u}}
\nc{\g}{\mathfrak{g}}

\newcommand{\G}{{\mathbb G}}
\newcommand{\Gf}{{\mathbb G}_{\mathrm{f}}}

\newcommand\GL{\operatorname{GL}}
\newcommand\id{\operatorname{id}}

\newcommand\gdim{\operatorname{GKdim}}

\newcommand\ad{\operatorname{ad}}

\def\ot{\otimes}

\def\C{\Bbbk}

\def\N{\mathbb{N}}
\def\Z{\mathbb{Z}}

\def\mX{\mathcal{X}}

\newcommand\hgt{\operatorname{ht}}

\newcommand{\bq}{{\bf q}}

\newcommand{\Gt}[1]{\mathtt{G}(#1)}
\newcommand{\Gtp}[1]{\mathtt{G}'(#1)}

\def\B{\mathfrak{B}}
\def\K{\mathcal{K}}
\def\F{\mathfrak{F}}

\newcommand{\hlist}{\mathsf{T}}

\def\gap{\texttt{GAP}}

\newcommand{\linethree}[5]{\xymatrix@C-4pt{
\overset{#1}
{\underset{\ }{\circ}}\ar@{-}[r]^{#2}
& \overset{#3}{\underset{\ }{\circ}}\ar@{-}[r]^{#4}
& \overset{#5}{\underset{\ }{\circ}}
}
}

\newcommand{\tri}[6]{
\xymatrix@C-4pt@R-18pt{ & \overset{#5}{\underset{\ }{\circ}} \ar @{-}[ld]_{#6} \ar @{-}[rd]^{#4} & \\
\overset{#1}{\underset{\ }{\circ}} \ar @{-}[rr]_{#2} & & \overset{#3}{\underset{\ }{\circ}}} }

\newcommand{\linefour}[7]{\xymatrix@C-4pt{ 
\overset{#1}{\underset{\ }{\circ}} \ar@{-}[r]^{#2} & \overset{#3}{\underset{\ }{\circ}}\ar@{-}[r]^{#4} &
\overset{#5}{\underset{\ }{\circ}} \ar@{-}[r]^{#6} & \overset{#7}{\underset{\ }{\circ}}     }}

\newcommand{\tripod}[8]{
\xymatrix@C-4pt@R-10pt{ & \overset{#3}{\circ} \ar @{-}[d]^{#4} & \\
\overset{#1}{\underset{\ }{\circ}} \ar @{-}[r]_{#6} & \overset{#5}{\underset{\ }{\circ}} \ar @{-}[r]_{#8} & \overset{#7}{\underset{\ }{\circ}}
}}

\newcommand{\tadpole}[8]{
\xymatrix@C-4pt@R-18pt{ & \overset{#3}{\underset{\ }{\circ}} \ar @{-}[ld]_{#2} \ar @{-}[rd]^{#4} & & & \\
\overset{#1}{\underset{\ }{\circ}} \ar @{-}[rr]_{#6} & & \overset{#5}{\underset{\ }{\circ}} \ar @{-}[rr]^{#8} & & \overset{#7}{\underset{\ }{\circ}}
}}

\begin{document}

\lstset{language=GAP,
showstringspaces=false,
xleftmargin=0.0cm,
xrightmargin=0.0cm,
basicstyle=\small\ttfamily,
frame=single,
framerule=0pt,
}



\title[Finite GKdim Nichols algebras and finite root systems]{Finite GK-dimensional Nichols algebras of diagonal type and finite root systems}

\author[Angiono, Garc\'ia Iglesias]{Iv\'an Angiono, Agust\'in Garc\'ia Iglesias}

\address{FaMAF-CIEM (CONICET), Universidad Nacional de C\'ordoba,
Medina A\-llen\-de s/n, Ciudad Universitaria, 5000 C\' ordoba, Rep\'ublica Argentina.}

\email{(ivan.angiono|agustingarcia)@unc.edu.ar}

\thanks{\noindent 2010 \emph{Mathematics Subject Classification.}
16T05. \newline The work was partially supported by CONICET,
FONCyT-ANPCyT, Secyt (UNC)}

\begin{abstract}
Let $(V,c)$ be a finite-dimensional braided vector space of diagonal type. 
We show that the Gelfand Kirillov dimension of the Nichols algebra $\mathfrak{B}(V)$ is finite if and only if the corresponding root system is finite, 
that is $\mathfrak{B}(V)$ admits a PBW basis with a finite number of generators.

This had been conjectured in \cite{AAH-conj} and proved for $\dim V=2,3$ in \cite{AAH-rank2}, \cite{AnG-rank3} respectively.
\end{abstract}

\maketitle

\tableofcontents

\section{Introduction}\label{sec:intro}
We fix an algebraically closed field $\Bbbk$ of characteristic zero. 
Short after the celebrated conference of Drinfeld \cite{Drinfeld} (presented by Cartier) in the ICM in 1986 Hopf algebras and quantum groups became an active and fruitful topic of research within many areas of mathematics. The study of Hopf algebras in general is a broad field and thus it was necessary to design a roadmap with the introduction of some restrictions, sufficiently lax to include all interesting examples and tight enough to achieve meaningful results. In this direction Rosso proposed to look at Hopf algebras with finite Gelfand-Kirillov dimension ($\gdim$ for short) in \cite{R quantum groups}.

In order to classify these objects and to obtain new examples along the way, there are at least three possible approaches. The first one is to classify Hopf algebras of low $\gdim$, see e.g. the survey articles \cite{BZ,Goodearl}. Another alternative is to assume some algebraic constrains such as notherian, AS regular, between others \cite{BG,BZ}. 

Finally one can focus on the coalgebra structure, namely to study connected Hopf algebras \cite{BGZ,LSZ}, or more generally pointed Hopf algebras.
In this last setting one can follow the lines of the Lifting Method for the finite dimensional case \cite{ASc}. Moreover, as in loc.cit.~one starts off by looking at Hopf algebras with abelian coradical. This leads to the consideration of Nichols algebras with finite $\gdim$. Here the action of the underlying abelian group can be semisimple or not. In the first case the corresponding braided vector space is of diagonal type and one considers the associated root system introduced in \cite{H-inv}. Indeed it was suggested in \cite{AAH-conj} that this notion is intimately related with finite $\gdim$:

\begin{conjecture}
\label{conjecture}
Let $(V,c)$ be a braided vector space of diagonal type such that $\gdim \B(V)<\infty$. Then its generalized root system is finite.
\end{conjecture}

The goal of this article is to establish the validity of this statement.
We recall that it is true when  $\dim V=2$ by \cite[Theorem 1.2]{AAH-rank2}, while $\dim V=3$ is shown in \cite[Theorem 6.1]{AnG-rank3}.

\medbreak

The collection of all braided vector spaces of diagonal type with finite root system is presented in \cite{H-full} in terms of Dynkin diagrams, listed in four tables according with the rank (the dimension of underlying vector space).

\begin{notation}
	We write $\hlist$ to denote the full set of diagrams appearing in any of these tables.
\end{notation}

The diagrams in $\hlist$ are decorated with elements in $\Bbbk$ as labels on each one of the vertices and the edges. We distinguish two kind of diagrams: \emph{finite} and \emph{parametric}. The parametric diagrams are those for which at least one label can be arbitrarily chosen from (an infinite subset of) $\Bbbk$. In turn, we say that a diagram is finite when all the labels are fixed.

\medbreak

Turning back to the action of the abelian group, the non semisimple case is considered in \cite{AAH-conj,AAM}. These classification results rely on the validity of the conjecture. Another context in which the conjecture is relevant is the second step of the Lifting Method, namely the classification of post (and pre) Nichols algebras with finite $\gdim$ \cite{ASa,ACS1,ACS2}.

\medspace

Next, we state our main result.

\begin{theorem}
Let $(V,c)$ be a finite-dimensional braided vector space of diagonal type.
Then  $\gdim \B(V)<\infty$ if and only if its generalized root system is finite.
\end{theorem}

As it was already stated in  \cite{AAH-conj}, we only need to prove the direct implication: the converse is already known. 
We nevertheless sketch a proof here for convenience of the reader: If the root system is finite, then  $\B(V)$ has a PBW basis with a finite set of generators, by definition.
This basis induces an algebra filtration $\F$ as in \cite[Proposition 10.1]{DP-quantum}. The associated graded object is a (truncated) $q$-polynomial algebra on the PBW generators. 
Hence $\gdim \B(V)<\infty$ as it coincides with that of this $q$-polynomial algebra.

We also remark that it is enough to consider connected diagrams (i.e.~~indecomposable braided vector spaces). Otherwise, the Nichols algebra factorizes as the twisted tensor product of the Nichols algebras of each connected component by \cite[Theorem 2.2]{Grana}: We remark that although that theorem is actually stated for finite-dimensional Nichols algebras, the proof extends naturally to the case of {\it any} finite collection of Nichols algebras associated to a decomposable braided vector space.    

Let us fix a connected braided vector space $(V,c)$, with $\dim V=\theta$ and such that $\gdim\B(V)<\infty$. If $\theta\leq 3$, then the root system is finite, as stated above. We extend this, recursively, to any dimension in Theorems \ref{thm:rank4-conj-true} ($\theta=4$), \ref{thm:rank5-conj-true} ($\theta=5$), \ref{thm:rank-6} ($\theta=6$),  \ref{thm:rank-7} ($\theta=7$) and \ref{thm:rank-ge-8} ($\theta\geq 8$). 

When $4\leq \theta\leq 7$ we proceed as follows. First we discard those (unlabeled) diagram shapes which necessarily imply infinite Gelfand-Kirillov dimension. For example, we show that there are no cycles of length $\geq 4$ in Proposition \ref{prop:n-cycles}. For each of the remaining {\it admissible} shapes of rank $\theta$, we prove the following trichotomy: either $\gdim \B(V)=\infty$, the diagram belongs to $\hlist$ or else the diagram belongs to a finite set of diagrams constructed by gluing diagrams of rank $\theta-1$ which are either finite or parametric, but evaluated in a finite set $\G_f$ of roots of 1. This now reduces the problem to this finite set of diagrams. We use \gap\, to write algorithms that construct this set and that discard the diagrams which do not belong to $\hlist$, by implementing some criteria based on ideas in \cite{AnG-rank3}.  

Finally, case $\theta\geq 8$ does not require the use of \gap: this is consistent with the fact that in Heckenberger's list higher ranks correspond to infinite families of Cartan, super or standard type. That is, there are no exceptional diagrams for a given rank in this case.

\section{Preliminaries}\label{sec:prels}

In this section we fix some notation and recall some of the previous results that will be needed along the work.

Let $m\leq n\in\N$. We set $\I_{m,n}=\{t\in\N : m\leq t\leq n\}$, $\I_n\coloneqq \I_{1,n}$. 
We write $\G_n$ for the group of $n$th roots of 1 in $\Bbbk$ and we let $\G_n'$ be the subset of primitive roots of order $n$. 

If $A$ is an associative algebra, we let $\gdim A$ denote its Gelfand-Kirillov dimension. 
We recall two basic properties of $\gdim$:
Let $A,B$ be two associative algebras: if either $A\subseteq B$ is a subalgebra or there exists a surjective algebra morphism $B\twoheadrightarrow A$, then $\gdim A\leq \gdim B$. 
See \cite{KL} for details.

\subsection{Nichols algebras of diagonal type}

A braided vector space is a pair $(V,c)$ such that $V$ is a vector space and $c\in\GL(V\ot V)$ is a solution of the braid equation:
\[
(c\ot\id)(\id\ot\, c)(c\ot\id)=(\id\ot\, c)(c\ot\id)(\id\ot\, c).
\]
Recall that $(V,c)$ is of diagonal type when there is a basis $\{x_1,\dots,x_\theta\}$ of $V$, together with a matrix $\bq\in(\Bbbk^{\times})^{\theta\times\theta}$
such that 
the braiding $c$ given by the formula:
\[
c(x_i\ot x_j)=q_{ij}x_j\ot x_i, \qquad i,j\in\I_\theta.
\]  
The matrix $\bq$ is called the braiding matrix of $(V,c)$; we write $\qti{ij}\coloneqq q_{ij}q_{ji}$.
As well, for each $\alpha=(a_i), \beta=(b_i)\in\Z^\theta$, we set
\begin{align}\label{eqn:q-alphabeta}
\bq_{\alpha,\beta}=\prod_{i,j\in\I_{\theta}}q_{ij}^{a_ibj}. 
\end{align} 
The Dynkin diagram of $\bq$ is the labeled graph with vertices $i\in\I_\theta$, labeled by $q_{ii}$, and an edge between $i$ and $j$ when $\qti{ij}\neq 1$, labeled with this scalar. 

\medbreak

In the particular case when there is $c_{ij}\in \Z$, $i,j\in\I_\theta$, such that $q_{ii}^{c_{ij}}=\qti{ij}$ for all $i\ne j\in\I_\theta$, we say that $\bq$ is of \emph{Cartan} type.

\medbreak

Given $(V,c)$, the Nichols algebra $\B(V)$ is a $\N_0$-graded connected braided Hopf algebra defined as the quotient of the tensor algebra $T(V)$ 
for which the subspace of primitive elements coincides with $V$. When $(V,c)$ is of diagonal type, with matrix $\bq$, we write $\B_{\bq}=\B(V)$.
We refer to \cite{A-leyva,AA-diag-survey} for more details on Nichols algebras. 

\subsection{Root systems}

Next, we recall the concept of root system $\varDelta^{\bq}$ of a Nichols algebra of diagonal type $\B_{\bq}$. To do this, we fix a PBW basis
\begin{equation*}
\{\ell_1^{m_1}\dots \ell_k^{m_k} | k\in\N_0, \ell_1>\dots >\ell_k\in L, 0<m_i<\hgt(\ell_i), i\in\I_k \},
\end{equation*}
where $L\subset \B_{\bq}$ is a subset of $\Z^{\theta}$-homogeneous elements and $\hgt\colon L\to \N\cup\{\infty\}$ is a height function, see \cite{K1}. 
Now, the set of roots of $\bq$ is
\begin{align*}
\varDelta^{\bq}&=\varDelta^{\bq}_+\cup -\varDelta^{\bq}_+, & \text{for } \quad \varDelta^{\bq}_+&=\{\deg \ell : \ell\in L  \}.
\end{align*}
By \cite[Lemma 4.7]{HS}, $\varDelta_+\coloneqq \varDelta_+^{\bq}$ does not depend on the choice of $L$.

\medbreak

Let $i\in \I\coloneqq \I_\theta$ be such that for all $j\neq i$ there exists $n\in\N_0$ such that $(n+1)_{q_{ii}}(1-q_{ii}^n q_{ij}q_{ji})=0$. 
In this case we say that  \emph{we can reflect $\bq$ at} $i$.

Henceforth we assume that we can reflect at every vertex $i\in \I$ and consider as in \cite{H-inv}:
\begin{itemize}[leftmargin=*]
\item The generalized Cartan matrix  $C^{\bq}=(c_{ij}^{\bq})$, given by 
\begin{align*}
c_{ii}^{\bq}&=2, & c_{ij}^{\bq}&:=-\min\{n\in\N_0:(n+1)_{q_{ii}}(1-q_{ii}^nq_{ij}q_{ji})=0\}, \quad j\neq i.
\end{align*}
Set $m_{ij}^{\bq}=-c_{ij}^{\bq}$, $i,j\in \I$. For simplicity, we shall often omit the superscript when $\bq$ is clear by the context.
\item The reflections $s_i^{\bq} \in \mathrm{GL}(\Z^{\theta})$, $i\in\I$, given by 
\begin{align*}
s_i^{\bq} (\alpha_j) &= \alpha_j - c_{ij}^{\bq}\alpha_i, &  j & \in \I.
\end{align*}
\item The braiding matrices $\rho_i \bq$, $i\in\I$, given by
\begin{align*}
(\rho_i \bq)_{jk}&:= q_{s_i^{\bq}(\alpha_j), s_i^{\bq}(\alpha_k)} 
=  q_{jk}q_{ik}^{-c_{ij}^{\bq}}q_{ji}^{-c_{ik}^{\bq}}q_{ii}^{c_{ij}^{\bq} c_{ik}^{\bq}}, & j,k&\in\I.
\end{align*}
\item The braided vector spaces of diagonal type $\rho_i V$ with braiding matrix $\rho_i \bq$. 
\end{itemize}
Both $\rho_i \bq$ and $\rho_i V$ are called the reflection at vertex $i$ of $\bq$.

Let $\mX\coloneqq \mX_{\bq}$ be the collection of all braiding matrices $\bq'=\rho_{i_k} \dots \rho_{i_1} \bq$ obtained from $\bq$ by a finite number of successive reflections. 
Each $\bq'\in\mX$ is said to be \emph{Weyl equivalent} to $\bq$. 

We say that $\bq$ \emph{admits all reflections} if we can reflect $\bq'$ at every $i\in \I$ for all $\bq'\in\mX$. In this case, the \emph{Weyl groupoid} of $\bq$ is defined as the subgroupoid of $\mX\times\GL(\Z^{\theta})\times \mX$ generated by
\begin{align*}
\sigma_i^{\bq'} & := (\bq', s_i^{\bq'},\rho_i\bq'), & & \bq'\in\mX, i\in\I.
\end{align*}
The collection $(\varDelta^{\bq'}_+)_{\bq' \in \mX}$ is the \emph{generalized root system} of $\bq$. We remark that by \cite{H-inv, AA}
$\varDelta^{\rho_i \bq }_+ = s_i^{\bq}(\varDelta^{\bq}_+-\{\alpha_i\})\cup \{\alpha_i\}$ for all $i\in\I$ and 
\begin{align}\label{eq:GKdim-reflection}
\gdim \B_{\rho_i \bq } &= \gdim \B_{\bq}.
\end{align}

We remark that when $\bq$ and ${\bf r}$ give rise to the same Dynkin diagram (namely, they are {\it twist equivalent}), then  $\B_{\bq}\simeq \B_{{\bf r}}$ as $\N_0$-graded vector spaces and thus $\gdim \B_{\bq} = \gdim\B_{{\bf r}}$. As well, the corresponding sets of roots $\varDelta_+^{\bq}$ and $\varDelta_+^{\bf r}$  coincide. Therefore, for our purposes it is enough to look at diagrams rather than matrices. 

The list of all diagrams with finite set of roots was described by Heckenberger in \cite{H-full}.

\subsection{Tools and criteria}\label{sec:criteria}

Equation \eqref{eq:GKdim-reflection} above states that the Gelfand-Kirillov dimension is preserved along a Weyl equivalence class. This will be a useful asset throughout the paper. 
Next we recall some other useful results, together with some tools developed in \cite{AnG-rank3}, to deal with  the proof of the conjecture.
We start with the following.
\begin{theorem}\cite[Theorem 5.1]{AnG-rank3}\label{thm:Cartan}
Let $\bq$ be a matrix of Cartan type. Then $\gdim\B_{\bq}<\infty$ if and only if its root system is finite.
\end{theorem}

Fix $\theta\in\N$ and a braiding matrix $\bq$ of rank $\theta$. 
In \cite[Proposition 3.3]{AnG-rank3}, we consider for each $\omega\in\Z^\theta$ a subquotient $\K_{\omega}$ of $\B_{\bq}$. 
As explained in \cite[Snapshot (9)]{AnG-rank3}, it is possible to obtain a Nichols algebra $\B_{\bar{\bq}}$ as a subquotient of $\K_{\omega}$, where the braiding matrix $\bar{\bq}$ has rank $\theta-1$. 
For this, we pick suitable primitive elements $x_{\beta_1}, \dots, x_{\beta_{\theta-1}}$ in $\K_\omega$, of degrees $\beta_1,\dots,\beta_{\theta-1}\in\Z^\theta$ such that $(\beta_i|\omega)=0$ for each $i$. We set ${\bf B}=\{\beta_i:i\in\I_{\theta-1}\}$ and 
$\bar{\bq}=(\bq_{\beta\gamma})_{\beta,\gamma\in{\bf B}}$ is computed as in \eqref{eqn:q-alphabeta}.

In particular, if $\gdim\B_{\bq}<\infty$, then $\gdim\B_{\bar{\bq}}<\infty$ for each possible choice of $\beta_1,\dots,\beta_{\theta-1}\in\Z^\theta$. This is the key idea behind our recursive step of the proof of the conjecture.

Next, we write down the explicit choices of $\omega=\sum_{i=1}^\theta c_i\alpha_i\in\Z^\theta$ and ${\bf B}$ we make and the corresponding criteria determined by each choice. In general, if $c_i=0$ for some $i$, then we  pick $\alpha_i\in {\bf B}$, with $x_{\alpha_i}=x_i$. Thus, for each choice, we just explicit the primitive elements $x_\beta$, for $\beta\neq\alpha_i\in {\bf B}$. When $\beta=m\alpha_i+\alpha_j$, we choose, depending on $\omega$ as in \cite[\S 4.1]{AnG-rank3}:
\begin{align*}
x_\beta&=\ad(x_i)^m(x_j) &&\text{or} & x_\beta&=\ad'(x_i)^m(x_j)\coloneqq  [\dots[[x_j,x_i]_c,x_i]_c,\dots,x_i]_c.
\end{align*}

\begin{criterium}[for $i,j,n$]\label{crit-A}
Let $i\neq j$ be such that $m_{ij}\neq 0$ and let $n\leq m_{ij}$. We choose ${\bf B}=\{\alpha_k|k\neq i,j\}\cup \{\alpha_j+n\alpha_i\}$.
\end{criterium}
Here $\omega=\alpha_i-n\alpha_j$.
This generalizes Criterium 1 in \cite[\S 4.1]{AnG-rank3} from rank three to an arbitrary rank $\theta$. The proof readily extends to this case.

\begin{criterium}\label{crit-B}
Let $i,j,k$ be pairwise different with $m_{ik}, m_{jk}\neq 0$. 

\noindent
We choose ${\bf B}=\{\alpha_p|p\neq i,j,k\}\cup \{\alpha_i+\alpha_k,\alpha_j+\alpha_k\}$.
\end{criterium}
Here $\omega=\alpha_i+\alpha_j-\alpha_k$.
This generalizes Criterium 3 in \cite[\S 4.1]{AnG-rank3}. 

\begin{criterium}\label{crit-C}
Let $\theta\geq 6$. We assume that $m_{i\,i+1}\neq 0$, $i\in I_{2,\theta-1}$.

\noindent
We choose ${\bf B}=\{\alpha_1,\alpha_\theta\}\cup \{\alpha_i+\alpha_{i+1}|i\in\I_{2,\theta-2}\}$.
\end{criterium}
Here $\omega=\sum_{i\in \I_{2,\theta-1}}(-1)^i\alpha_i$ and $x_{\alpha_i+\alpha_{i+1}}=\begin{cases}
\ad(x_i)(x_{i+1}), &i\text{ even,}\\
\ad'(x_{i})(x_{i+1}), &i\text{ odd.}
\end{cases}$

\begin{criterium}\label{crit-D}
Let $\theta=4$. We assume that $m_{21}, m_{23}\neq 0$ and let $m\leq m_{21}$, $n\leq m_{23}$.
We choose ${\bf B}=\{\alpha_1+m\alpha_2,n\alpha_2+\alpha_3,\alpha_4\}$. 
\end{criterium}
Here $\omega=\alpha_2-m\alpha_1-n\alpha_3$, $x_{\alpha_1+m\alpha_2}=\ad(x_2)^m(x_1)$, $x_{n\alpha_2+\alpha_3}=\ad(x_2)^n(x_3)$.

\begin{criterium}\label{crit-E}
Let $\theta=4$. We assume that $m_{12}, m_{23},m_{34}\neq 0$ and let $m\leq m_{12}$, $n\leq m_{23}$, $p\leq m_{34}$.
We choose ${\bf B}=\{m\alpha_1+\alpha_2,n\alpha_2+\alpha_3,p\alpha_3+\alpha_4\}$. 
\end{criterium}
Here $\omega=\alpha_1-m\alpha_2+mn\alpha_3-mnp\alpha_4$ and
\begin{align*}
x_{m\alpha_1+\alpha_2}&=\ad(x_1)^m(x_2), & x_{n\alpha_2+\alpha_3}&=\ad'(x_2)^n(x_3), & x_{p\alpha_3+\alpha_4}&=\ad(x_3)^p(x_4).
\end{align*}

\smallbreak

We shall also use Criteria \ref{crit-D} and \ref{crit-E} for proper permutations of $\I_4$. We shall also make use of the following extension of the later (for proper permutations of $\I_\theta$ as well).
\begin{criterium}\label{crit-F}
Let $\theta\geq 5$. We assume that $m_{12}, m_{23},m_{34}\neq 0$.

\noindent
We choose ${\bf B}=\{\alpha_1+\alpha_2,\alpha_2+\alpha_3,\alpha_3+\alpha_4\}\cup \{\alpha_p|p\geq 5\}$.
\end{criterium}
Here $\omega=\alpha_1-\alpha_2+\alpha_3-\alpha_4$ and
$x_{\alpha_1+\alpha_2}=[x_1,x_2]_c$, $x_{\alpha_2+\alpha_3}=[x_3,x_2]_c$, $x_{\alpha_3+\alpha_4}=[x_3,x_4]_c$.

\begin{remark}\label{rem:criterium-reflections}
Note that, to show that $\gdim \B_{\bq}=\infty$, it is enough to show that one of the criteria fails for one of the reflections $\rho_i\bq$, $i\in\I$, by \eqref{eq:GKdim-reflection}. 
\end{remark}

As an example, we exhibit the diagram associated to the matrix $\bar{\bq}=(\bar{q}_{ij})$ obtained after the application of Criterium \ref{crit-A} for $(1,2;n)$. The set of vertices is $\I'=\{\star\}\cup\I_{3,\theta}\,(\simeq \I_{\theta-1})$. The labels become, using \eqref{eqn:q-alphabeta}:
\begin{align*}
\bar{q}_{ii}&=\begin{cases}
q_{11}^{n^2}\qti{12}^nq_{22}, & i=\star,\\
q_{ii}, & i\in \I_{3,\theta};
\end{cases}
&
\widetilde{\bar{q}_{ij}}&=\begin{cases}
\qti{1j}^n\qti{2j}, & i=\star, j\in \I_{3,\theta},\\
\qti{ij}, & i,j\in \I_{3,\theta}.
\end{cases}´
\end{align*}

\section{Cycles}\label{sec:cycles}

In this section we focus on diagrams containing a cycle. We will prove that if the diagram of $\bq$ contains a cycle and $\gdim\B_{\bq}<\infty$, then the cycle is a triangle. For the proof, we first discard diagrams containing a $4$-cycle, which are obtained by gluing two diagrams of rank three in $\hlist$. When one of these subdiagrams is parametric we show that $\gdim\B_{\bq}=\infty$ unless the parameters belong to the following set of roots of 1:
\begin{align}\label{eq:defn-G-f}
\Gf &:= \left(\G_{10}\cup \G_{12}\cup \G_{18}\right)\setminus \{1\}.
\end{align}
After this reduction, we are left with a finite (though large) collection of diagrams with labels in $\Gf$, that we attack with \gap\,.

Once the case of rank 4 is solved, we deal with cycles of bigger length using a recursive argument.

\medspace

We consider \emph{one-parametric} and \emph{two-parametric} lines of rank 3 in $\hlist$ as functions
$\mathtt{G}:\C^{\times}\to (\C^{\times})^5$, respectively $\mathtt{G}:(\C^{\times})^2\to (\C^{\times})^5$. Explicitly these lines are
\begin{align}\label{eq:parametric-line-3}
\begin{aligned}
& \linethree{q}{q^{-1}}{q}{q^{-1}}{q}, &
& \linethree{q^2}{q^{-2}}{q^2}{q^{-2}}{q}, &
& \linethree{q}{q^{-1}}{q}{q^{-2}}{q^2},
\\ 
& \linethree{-1}{q^{-1}}{q}{q^{-1}}{q}, &
& \linethree{-1}{q}{-1}{q^{-1}}{q}, &
& \linethree{-1}{q^{-2}}{q^2}{q^{-2}}{q},
\\ 
& \linethree{-1}{q^2}{-1}{q^{-2}}{q}, &
& \linethree{q^2}{q^{-2}}{-1}{q^{2}}{-q^{-1}}, &
& \linethree{-1}{q^{-1}}{q}{q^{-2}}{q^2},
\\ 
& \linethree{-1}{q}{-1}{q^{-2}}{q^2}, &
& \linethree{-1}{q^{-1}}{q}{q^{-3}}{q^3}, &
& \linethree{-1}{q}{-1}{q^{-3}}{q^3},
\\ 
& \linethree{q^3}{q^{-3}}{-1}{q^{2}}{-q^{-1}}, &
& \linethree{q}{q^{-1}}{-1}{q}{q^{-1}}, &
& \linethree{-1}{q}{-1}{q^{-1}}{-1},
\\ 
& \linethree{-1}{q^{-1}}{q}{q^{-1}}{-1}, &
& \linethree{q}{q^{-1}}{-1}{r^{-1}}{r}.
\end{aligned}
\end{align}
For instance the last diagram corresponds to the function 
$$ \mathtt{G}(q,r)=(q,q^{-1},-1,r^{-1},r).$$
The same applies for one and two-parametric triangles:
\begin{align}\label{eq:parametric-triangle}
& \tri{q}{q^{-1}}{-1}{q^{2}}{-1}{q^{-1}}, &
& \tri{q}{q^{-1}}{-1}{q^{3}}{-1}{q^{-2}}, &
& \tri{-1}{q}{-1}{r}{-1}{q^{-1}r^{-1}}.
\end{align}

\begin{remark}\label{rem:rk-3-finite-in-Gf}
The entries of each finite line of rank 3 in $\hlist$ belong to $\G_{6}\cup \G_{9}$, 
while those of each finite triangle of rank 3 in $\hlist$ belong to $\G_{6}$.
\end{remark}

\begin{lemma}\label{lem:4-cycles}
If the diagram of $\bq$ contains a $4$-cycle, then $\gdim \B_{\bq}=\infty$.
\end{lemma}

\pf
Up to relabeling, we may assume that $\qti{12},\qti{23},\qti{34},\qti{14}\neq 1$. It is enough to reduce to rank 4 by considering the submatrix $\bq_{|\I_4}$.
We fix the following numbering of the vertices:
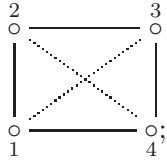
\begin{align*}
&\xymatrix@C40pt@R-5pt{ \overset{2}{\circ} \ar @{-}[d] \ar @{-}[r] \ar @{.}[rd]& 
\overset{3}{\circ} \ar @{-}[d] \\
\underset{1}{\circ} \ar@{-}[r] \ar@{.}[ru] & \underset{4}{\circ};}
\end{align*}
the dot line means that it may or may not appear, depending on $\qti{13},\qti{24}$. 
Even more, we may assume that all connected proper subdiagrams of rank 3 belong to $\hlist$ by \cite[Theorem 6.1]{AnG-rank3}. We proceed in several steps.

\begin{stepcyc}\label{stepcyc:diagonal+parametric}
If $\qti{13}\ne 1$ and $\bq_{|\I_{3}}$ is parametric, namely a function $\mathtt{G}(q)$ or $\mathtt{G}(q,r)$, then either $\gdim\B_{\bq}=\infty$ or else $q,r\in\Gf$.
\end{stepcyc}
\pf
If the diagram of $\bq_{|\I_{3}}$ is the first triangle in \eqref{eq:parametric-triangle}, then the diagram of $\bq$ has one of the following shapes:
\begin{align}\label{eq:cycles-case-i}
&\xymatrix@C40pt@R-5pt{ \overset{q}{\circ} \ar @{-}[d]_{q^{-1}} \ar @{-}[r]^{q^{-1}} \ar @{.}[rd]& 
\overset{-1}{\circ}\ar @{-}[d]^r \\
\overset{-1}{\circ} \ar@{-}[r]_t \ar@{-}[ru]_{q^2} & \overset{s}{\circ},}
&
&\xymatrix@C40pt@R-5pt{ \overset{-1}{\circ} \ar @{-}[d]_{q^2} \ar @{-}[r]^{q^{-1}} \ar @{.}[rd]& 
\overset{q}{\circ}\ar @{-}[d]^r \\
\overset{-1}{\circ} \ar@{-}[r]_t \ar@{-}[ru]_{q^{-1}} & \overset{s}{\circ},}   & & r,s,t \ne 1.
\end{align}

If $\bq_{|\{1,3,4\}}$ is finite, then $q\in\Gf$: indeed, by Remark \ref{rem:rk-3-finite-in-Gf} either $q^{-1}\in\G_6$ or $q^2\in\G_6$. Thus we assume that $\bq_{|\{1,3,4\}}$ is parametric. For the first diagram in \eqref{eq:cycles-case-i}, we have three cases:
\begin{itemize}[leftmargin=*]\renewcommand{\labelitemi}{$\circ$}
\item $st=sr=1$, $s^2=q^2$. Criterium \ref{crit-A} for $(2,3;1)$, that is $i=2,j=3$ and $n=1$, gives the diagram
$$ \tri{-1}{s^{-1}}{s}{s^{-1}\qti{24}}{-1}{q}.$$
If $\qti{24}\ne s, q$, then this diagram is a triangle not in \eqref{eq:parametric-triangle}, since the product of the labels of the edges is not 1. 
If $\qti{24}=s$ or $\qti{24}=q$, then $\bq_{|\{1,2,4\}}$ is a triangle with exactly one vertex labeled with -1, so it is not in \eqref{eq:parametric-triangle}. In any case, $\gdim \B_{\bq}=\infty$.
\item $s^2t=sr=1$, $s^3=q^2$. Criterium \ref{crit-A} for $(2,3;1)$ gives the diagram
$$ \tri{-1}{s^{-2}}{s}{s^{-1}\qti{24}}{-1}{q}.$$
If $\qti{24}\ne s, q$, then this diagram is a triangle not in \eqref{eq:parametric-triangle}. 
Otherwise $\bq_{|\{1,2,4\}}$ is a triangle not in \eqref{eq:parametric-triangle}. In any case, $\gdim \B_{\bq}=\infty$.
\item $trq^2=1$, $s=-1$. Criterium \ref{crit-A} for $(2,3;1)$ gives the diagram
$$ \tri{-1}{t}{-1}{r\qti{24}}{-1}{q}.$$
If $\qti{24}\ne q, r^{-1}$, then this diagram is a triangle not in \eqref{eq:parametric-triangle} since the products of the labels in the edges is not 1. 
If $\qti{24}=q$, then $\bq_{|\{1,2,4\}}$ is a triangle not in \eqref{eq:parametric-triangle} by the same reason. 
If $\qti{24}=r^{-1}$, then $\bq_{|\{2,3,4\}}$ is a triangle not in \eqref{eq:parametric-triangle}.
In any case, $\gdim \B_{\bq}=\infty$.
\end{itemize}

For the second diagram in \eqref{eq:cycles-case-i}, $s=-1$, $r=q^{-a}$, $t=q^{a+1}$ for $a\in\I_2$.
Criterium \ref{crit-A} for $(2,3;1)$ gives the diagram
$$ \tri{-1}{q^{a+1}}{-1}{q^{-a}\qti{24}}{-1}{q}.$$
If $\qti{24}\ne q^{-2}, q^a$, then this diagram is a triangle not in \eqref{eq:parametric-triangle} since the products of the labels in the edges is not 1. 
If $\qti{24}=q^{-2}$, respectively $\qti{24}=q^a$, then $\bq_{|\{1,2,4\}}$, respectively $\bq_{|\{2,3,4\}}$, is a triangle not in \eqref{eq:parametric-triangle}. In any case, $\gdim \B_{\bq}=\infty$.

\medspace

If the diagram of $\bq_{|\I_{3}}$ is the second triangle in \eqref{eq:parametric-triangle}, then the diagram of $\bq$ has one of the following shapes:
\begin{align*}
&\xymatrix@C40pt@R-5pt{ \overset{q}{\circ} \ar @{-}[d]_{q^{-2}} \ar @{-}[r]^{q^{-1}} \ar @{.}[rd]& 
\overset{-1}{\circ}\ar @{-}[d]^r \\
\overset{-1}{\circ} \ar@{-}[r]_t \ar@{-}[ru]_{q^3} & \overset{s}{\circ},}
&
&\xymatrix@C40pt@R-5pt{ \overset{-1}{\circ} \ar @{-}[d]_{q^3} \ar @{-}[r]^{q^{-b}} \ar @{.}[rd]& 
\overset{q}{\circ}\ar @{-}[d]^r \\
\overset{-1}{\circ} \ar@{-}[r]_t \ar@{-}[ru]_{q^{-a}} & \overset{s}{\circ},}   & & r,s,t \ne 1, \, \{a,b\}=\{1,2\}.
\end{align*}
If $\bq_{|\{1,3,4\}}$ is finite, then $q\in\Gf$ by Remark \ref{rem:rk-3-finite-in-Gf}. If $\bq_{|\{1,3,4\}}$ is parametric, then the proof is analogous to the one for the first triangle, case-by-case.

\medspace

Otherwise the diagram of $\bq$ is as follows:
\begin{align*}
&\xymatrix@C40pt@R-5pt{ \overset{-1}{\circ} \ar @{-}[d]_{q} \ar @{-}[r]^{q^{-1}r^{-1}} \ar @{.}[rd]& 
\overset{-1}{\circ}\ar @{-}[d]^s \\
\overset{-1}{\circ} \ar@{-}[r]_u \ar@{-}[ru]_{r} & \overset{t}{\circ},}  & & q,r,s,t,u, qr \ne 1.
\end{align*}
As we assume that  $\bq_{|\{1,3,4\}}$ is in $\hlist$, we have that $rsu=1$. Again, Criterium \ref{crit-A} for $(2,3;1)$ gives the diagram
\begin{align}\label{eq:4-cycle-criterium-last-case}
\begin{aligned}
\tri{-1}{u}{t}{s\qti{24}}{q^{-1}r^{-1}}{qr}.
\end{aligned}
\end{align}
We assume that this diagram is in $\hlist$ since otherwise $\gdim \B_{\bq}=\infty$. 
Suppose that $s\qti{24}\ne 1$, so \eqref{eq:4-cycle-criterium-last-case} is a triangle with $us\qti{24}qr=1$, which implies that $q\qti{24}=1$;
thus $\bq_{|\{1,2,4\}}$ is a triangle where the product of the three edges is $u\ne 1$, so $\gdim \B_{\bq}=\infty$. 
But if $s\qti{24}=1$, then $\bq_{|\{2,3,4\}}$ is a triangle where the product of the three edges is $q^{-1}r^{-1}\ne 1$, so $\gdim \B_{\bq}=\infty$. 
\epf

\begin{stepcyc}\label{stepcyc:diagonal}
If either $\qti{13}\ne 1$ or $\qti{24}\ne 1$, then $\gdim\B_{\bq}=\infty$.
\end{stepcyc}
\pf
Up to relabeling the vertices we may assume that $\qti{13}\ne 1$. If $\bq_{|\I_{3}}$ is parametric, with parameters not in $\Gf$, then $\gdim\B_{\bq}=\infty$ by Step \ref{stepcyc:diagonal+parametric}. Otherwise we may assume that both triangles $\bq_{|\I_{3}}$ and $\bq_{|\{1,3,4\}}$ are either finite or parametric evaluated in $\Gf$. The result follows using \gap: We discard all squares with one or two diagonals as explained in \S \ref{gap:criteria4}, p.~\pageref{page:squares}.
\epf

Thus we can restrict to the case $\qti{13}=\qti{24}=1$ (squares without diagonals).

\begin{stepcyc}\label{stepcyc:wo-diagonal+parametric}
If $\qti{13}=\qti{24}=1$ and $\bq_{|\I_{3}}$ is parametric, that is a function $\mathtt{G}(q)$ or $\mathtt{G}(q,r)$, then either $\gdim\B_{\bq}=\infty$ or else $q,r\in\Gf$.
\end{stepcyc}
\pf
Assume that $\gdim\B_{\bq}<\infty$, so we want to prove that the parameters are evaluated in $\Gf$. 

If $\bq_{|\I_{3}}$ is the tenth, twelfth, thirteenth or seventeenth diagram in \eqref{eq:parametric-line-3}, then $\rho_2\bq$ has a parametric triangle with vertices $1,2,3$. By Step \ref{stepcyc:diagonal}, the parameters are evaluated in $\Gf$.
If $\bq_{|\I_{3}}$ is the ninth, respectively eleventh, diagram in \eqref{eq:parametric-line-3}, then $\rho_1\bq_{|\I_{3}}$ is the tenth, respectively the twelfth, diagram in \eqref{eq:parametric-line-3}, so the parameter is evaluated in $\Gf$ by the argument above.

If $\bq_{|\I_{3}}$ is the fifth diagram in \eqref{eq:parametric-line-3} (so $q\ne-1$), then the diagram of $\bq$ is
\begin{align}\label{eq:cycles-wo-diagonal-case-i}
&\begin{aligned}
\xymatrix@C40pt@R-5pt{ \overset{-1}{\circ} \ar @{-}[d]_{q} \ar@{-}[r]^{q^{-1}} & \overset{q}{\circ}\ar @{-}[d]^r \\
\overset{-1}{\circ} \ar@{-}[r]_t & \overset{s}{\circ},}
\end{aligned}
& &r,s,t \ne 1.
\end{align}
Criterium \ref{crit-A} for $(1,2;1)$ gives the diagram
\[
\tri{q}{t}{s}{r}{q}{q^{-1}}
\] 
We assume that this triangle is in $\hlist$ since otherwise $\gdim \B_{\bq}=\infty$. As $q\ne -1$, the triangle does not appear in \eqref{eq:parametric-triangle}; thus it is finite, so $q\in\Gf$.

If $\bq_{|\I_{3}}$ is the fourth diagram in \eqref{eq:parametric-line-3}, then $\rho_1\bq_{|\I_{3}}$ is the fifth diagram in \eqref{eq:parametric-line-3}, 
so the parameter is evaluated in $\Gf$.

If $\bq_{|\I_{3}}$ is the seventh diagram in \eqref{eq:parametric-line-3} ($q^2\ne-1$), then the diagram of $\bq$ is
\begin{align}\label{eq:cycles-wo-diagonal-case-ii}
&\begin{aligned}
\xymatrix@C40pt@R-5pt{ \overset{-1}{\circ} \ar @{-}[d]_{q^2} \ar@{-}[r]^{q^{-2}} & \overset{q}{\circ}\ar @{-}[d]^r \\
\overset{-1}{\circ} \ar@{-}[r]_t & \overset{s}{\circ},}
\end{aligned}
& &r,s,t\ne 1.
\end{align}
Criterium \ref{crit-A} for $(1,2;1)$  gives the triangle
$$\tri{q^2}{t}{s}{r}{q}{q^{-2}}$$ and the proof is analogous to the fifth diagram.

If $\bq_{|\I_{3}}$ is the sixth, resp. eighth, diagram in \eqref{eq:parametric-line-3}, then $\rho_1\bq_{|\I_{3}}$, resp. $\rho_2\bq_{|\I_{3}}$, is the seventh diagram, so the parameter is evaluated in $\Gf$.

If $\bq_{|\I_{3}}$ is the fifteenth diagram in \eqref{eq:parametric-line-3} ($q\ne-1$), then the diagram of $\bq$ is
\begin{align}\label{eq:cycles-wo-diagonal-case-iii}
&\begin{aligned}
\xymatrix@C40pt@R-5pt{ \overset{-1}{\circ} \ar @{-}[d]_{q} \ar@{-}[r]^{q^{-1}} & \overset{-1}{\circ}\ar @{-}[d]^r \\
\overset{-1}{\circ} \ar@{-}[r]_t & \overset{s}{\circ},}
\end{aligned}
& &r,s,t\ne 1.
\end{align}
Criterium \ref{crit-A} gives for $(1,2;1)$, respectively for $(2,3;1)$:
\begin{align*}
&\tri{q}{t}{s}{r}{-1}{q^{-1}} & 
&\tri{-1}{t}{s}{r}{q^{-1}}{q}
\end{align*}
As $q^2\ne 1$, either $q^{-1}rt\ne 1$ or else $qrt\ne 1$. Thus either the first or the second triangle is not in $\hlist$ since the product of the edges is not 1. 

If $\bq_{|\I_{3}}$ is the fourteenth, resp. sixteenth, diagram in \eqref{eq:parametric-line-3}, then $\rho_2\bq_{|\I_{3}}$, resp. $\rho_1\bq_{|\I_{3}}$, is the fifteenth diagram, so $\gdim \B_{\bq}=\infty$. 

Finally, assume that $\bq_{|\I_{3}}$ is either the first, the second or the third diagram in \eqref{eq:parametric-line-3} (that is, of Cartan type). Up to exchange vertices $1$ and $3$, we may assume that $q_{33}=q$. By Theorem \ref{thm:Cartan} either 1, 3 or 4 is not Cartan. If 3 is not Cartan, then $q=q_{33}=-1$ since $\bq_{|\I_{2,4}}$ is in $\hlist$. Similarly, if 1 is not Cartan, then $q_{11}=-1$, so $q\in\G_4$. Finally, if 4 is not Cartan, then $\bq_{|\{1,3,4\}}$ is either finite (in which case $q=q_{33}\in\Gf$), or one of the diagrams considered before and again $q=q_{33}\in\Gf$.
\epf

\begin{stepcyc}\label{stepcyc:wo-diagonal}
If $\qti{13}=\qti{24}=1$, then $\gdim\B_{\bq}=\infty$.
\end{stepcyc}

If $\bq_{|\I_{3}}$ is parametric with parameters not in $\Gf$, then $\gdim\B_{\bq}=\infty$ by Step \ref{stepcyc:wo-diagonal+parametric}. Otherwise we may assume that the four lines are either finite or parametric evaluated in $\Gf$. We discard these cases using \gap\, in \S \ref{gap:criteria4}, p.~\pageref{page:clean-squares}, removing all squares.
\epf

The previous result shows that the diagram of a matrix $\bq$ with finite $\gdim$ cannot contain a $4$-cycle. 
Next we show that the diagram cannot contain $n$-cycles for any $n\ge 4$.

\begin{proposition}\label{prop:n-cycles}
If the diagram of $\bq$ contains a $n$-cycle for some $n\ge 4$, then $\gdim \B_{\bq}=\infty$.
\end{proposition}
\pf
The proof is by induction on $n$. The case $n=4$ is Lemma \ref{lem:4-cycles}. 

Assume that $n\ge 4$ and the statement holds for all $k\in\I_{4,n}$. Up to permuting the vertices we may assume that $\qti{i\, i+1}, \qti{1\,n+1}\ne 1$ for $i\in\I_{n}$. 
\begin{itemize}[leftmargin=*]\renewcommand{\labelitemi}{$\diamond$}
\item If the cycle has a diagonal, then up to permutation we may assume that one of the vertices of this diagonal is $1$: that is, $\qti{1k}\ne 1$ for $k\in\I_{3,n}$. If $4\le k\le n$, then $\bq_{|\I_{k}}$ is a $k$-cycle; if $k=3$, then $\bq_{|\{1,3,4,\dots,n+1\}}$ is a $n$-cycle.
\item If the cycle does not have diagonals, then Criterium \ref{crit-A} for $(1,2;1)$  gives an $n$-cycle.
\end{itemize}
In any case, $\gdim \B_{\bq}=\infty$ by inductive hypothesis.
\epf

\section{Rank 4}\label{sec:rank4}

In this section we prove Conjecture \ref{conjecture} for matrices $\bq$ of rank 4 with connected diagram. 
We can restrict to those $\bq$ whose diagram does not contain a 4-cycle by Lemma \ref{lem:4-cycles}. 
The strategy follows the lines of \cite{AnG-rank3}.
In the first step we show that a trichotomy holds: $\bq$ is in $\hlist$, $\gdim\B_{\bq}=\infty$ or the diagram of $\bq$ is obtained by gluing two diagrams of rank 3 in $\hlist$ which are either finite or parametric evaluated in the set $\Gf$ as in \eqref{eq:defn-G-f}.
Once this is done, we are left with a finite (though large) collection of diagrams, that we attack with \gap: we see that $\bq$ is in $\hlist$ or $\gdim\B_{\bq}=\infty$.

\subsection{On parametric subdiagrams}\label{subsec:eval-param}
We start with the case in which there is a parametric subdiagram of rank 3. As in Section \ref{sec:cycles}, the parametric subdiagram is considered as a function $\mathtt{G}$ on one or two parameters, see \eqref{eq:parametric-line-3} and \eqref{eq:parametric-triangle}. To establish the trichotomy, we show a series of lemmas dealing with all combinations of finite and parametric subdiagrams of rank three:
\begin{enumerate}
\item A tadpole whose triangle is parametric;
\item A tadpole with a parametric line of rank 3;
\item A tripod with a parametric line of rank 3;
\item A line with a parametric line of rank 3.
\end{enumerate}
We fix the following numerations for vertices of tadpoles, tripods and lines:
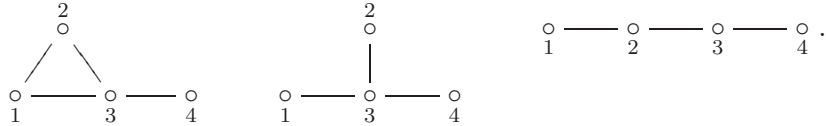
\begin{align*}
& \xymatrix@C-18pt@R-18pt{ & \overset{2}{\circ} \ar @{-}[ld] \ar @{-}[rd] & & & \\
\underset{1}{\circ} \ar @{-}[rr] & & \underset{3}{\circ} \ar @{-}[rr] & & \underset{4}{\circ}}
&
&\xymatrix@C-4pt@R-18pt{ & \overset{2}{\circ} \ar @{-}[d] & \\
\underset{1}{\circ} \ar @{-}[r] & \underset{3}{\circ} \ar @{-}[r] & \underset{4}{\circ}}
&
&\xymatrix@C-4pt@R-18pt{ \underset{1}{\circ} \ar@{-}[r] & \underset{2}{\circ} \ar @{-}[r] & \underset{3}{\circ} \ar@{-}[r] & \underset{4}{\circ}}.
\end{align*}

\begin{lemma}\label{lem:tadpole-with-param-triang}
Let $\bq$ be a tadpole such that $\bq_{|\I_3}$ is parametric of type $\mathtt{G}$. Then either $\bq$ belongs to $\hlist$, or $\bq_{|\I_3}=\Gt{q}$ for some $q\in\Gf$, or $\bq_{|\I_3}=\Gt{q,r}$ for some $q,r\in\Gf$, or else $\gdim \B_{\bq}=\infty$.
\end{lemma}
\pf
We assume that  $\gdim \B_{\bq}<\infty$. Thus each connected subdiagram of rank 3 belongs to $\hlist$ by \cite{AnG-rank3}.
We analyse each one of the three triangles in \eqref{eq:parametric-triangle}, including rotations.

\noindent $\triangleright$ Assume that $\bq_{|\I_3}$ is the first triangle in \eqref{eq:parametric-triangle}; i.e.~ the diagram of $\bq$ is
\begin{align*}
&\tadpole{q}{q^{-1}}{-1}{q^{2}}{-1}{q^{-1}}{t}{s} &&\text{for some }s,t\ne 1.
\end{align*}
If $\bq_{|\{1,3,4\}}$ is finite, then $q\in\G_{6}\cup\G_{9}\subseteq \Gf$, since $q$ is the label of a vertex. Similarly, if $\bq_{|\{2,3,4\}}$ is finite, then $q^2\in\G_{6}\cup\G_{9}\subseteq \Gf$, so $q\in\Gf$. Therefore we may assume that both $\bq_{|\{1,3,4\}}$ and $\bq_{|\{2,3,4\}}$ are parametric lines, so they are one of those in \eqref{eq:parametric-line-3}. Looking at $\bq_{|\{1,3,4\}}$, there are 6 possibilities:
\begin{description}[leftmargin=10pt]
\item[$s=q$, $t=-1$] The diagram of $\bq_{|\{2,3,4\}}$ is $\linethree{-1}{q^2}{-1}{q}{-1}$, which appears in \eqref{eq:parametric-line-3} iff either $q\in \G_4'$ (10th diagram) or $q\in \G_3'$ (15th diagram).

\item[$s^2=q$, $t=-1$] The diagram of $\bq_{|\{2,3,4\}}$ is $\linethree{-1}{s^4}{-1}{s}{-1}$, which appears in \eqref{eq:parametric-line-3} iff $s\in \G_5'$ (15th diagram).

\item[$s^3=q$, $t=-1$] The diagram of $\bq_{|\{2,3,4\}}$ is $\linethree{-1}{s^6}{-1}{s}{-1}$, which appears in \eqref{eq:parametric-line-3} iff either $s\in \G_4'$ (10th diagram) or $s\in \G_7'$ (15th diagram). If $s\in \G_7'$, then Criterium \ref{crit-A} for $(1,2;1)$ gives $\linethree{-1}{s^3}{-1}{s}{-1}$, which is not in $\hlist$.

\item[$s=r^2$, $t=-r^{-1}$, $q=r^3$] The diagram of $\bq_{|\{2,3,4\}}$ is $\linethree{-1}{r^6}{-1}{r^2}{-r^{-1}}$, which appears in \eqref{eq:parametric-line-3} iff either $r\in \G_{4}'$ (2nd diagram) or $r\in \G_8'$ (7th diagram). If $r\in \G_8'$, then  Criterium \ref{crit-A} for $(1,2;1)$ gives $\linethree{-1}{r^3}{-1}{r^2}{r^3}$, which is not in $\hlist$.

\item[$s=q$, $t=q^{-1}$] The diagram of $\bq_{|\{2,3,4\}}$ is $\linethree{-1}{q^2}{-1}{q}{q^{-1}}$, which appears in \eqref{eq:parametric-line-3} iff either $q\in \G_3'$ (5th diagram), $q\in \G_5'$ (10th diagram), $q\in \G_{7}'$ (12th diagram) or $q\in \G_4'$ (17th diagram). If $q\in \G_{7}'$, then  $\gdim \B_{\bq}=\infty$ applying  Criterium \ref{crit-A} for $(1,2;1)$.

\item[$st=1$] The diagram of $\bq_{|\{2,3,4\}}$ is $\linethree{-1}{q^2}{-1}{t^{-1}}{t}$, which appears in \eqref{eq:parametric-line-3} iff either 
$q\in \G_4'$, $t=-1$ (1st diagram), $t=q^2$ (5th diagram), $t=q^4$ (10th diagram), $t=q^6$ (12th diagram) or $q\in \G_4'$ (17th diagram). The second case appears in \cite[Table 3, row 9]{H-full}; for the third and the fourth cases,  Criterium \ref{crit-A} for $(1,2;1)$ gives, respectively, $\linethree{-1}{q}{-1}{q^{-4}}{q^4}$ and $\linethree{-1}{q}{-1}{q^{-6}}{q^6}$, so $q\in \G_3'\cup\G_4'\cup\G_5'  \subset \Gf$.
\end{description}

\medbreak

\noindent $\triangleright$ If $\bq_{|\I_3}$ is the first triangle in \eqref{eq:parametric-triangle}, rotated, then the diagram of $\bq$ is
\begin{align*}
&\tadpole{-1}{q^{2}}{-1}{q^{-1}}{q}{q^{-1}}{t}{s} &&\text{for some }s,t\ne 1.
\end{align*}
Again, we reduce to study the case in which both $\bq_{|\{1,3,4\}}$ and $\bq_{|\{2,3,4\}}$ are parametric lines. 
Looking at $\bq_{|\{1,3,4\}}$, there are 5 possibilities. Three of them are $s=q^{-2}=t^{-1}$, $s=q^{-3}=t^{-1}$, and
$s=t^{-2}=q^{-1}$:  Criterium \ref{crit-A} for $(1,2;1)$ gives $\linethree{q^2}{q^{-2}}{q}{s}{t}$, a braiding of Cartan type whose Cartan matrix is not finite, a contradiction with Theorem \ref{thm:Cartan}. The other two cases are $s=q^{-1}$, $t\in\{q,-1\}$, which belong to rows 8 and 12 of \cite[Table 3]{H-full}.

\medbreak

\noindent $\triangleright$ Assume that $\bq_{|\I_3}$ is the second triangle in \eqref{eq:parametric-triangle}; i.e.~ the diagram of $\bq$ is
\begin{align*}
&\tadpole{q}{q^{-1}}{-1}{q^{3}}{-1}{q^{-2}}{t}{s} &&\text{for some }s,t\ne 1.
\end{align*}
As above we assume that both $\bq_{|\{1,3,4\}}$ and $\bq_{|\{2,3,4\}}$ are parametric lines. 
Looking at $\bq_{|\{1,3,4\}}$, there are 4 possibilities:
\begin{description}[leftmargin=10pt]
\item[$s=q^{-2}$, $t=q^2=-1$] Thus, $q\in \G_4'\subset \Gf$.
\item[$s=q^{2}$, $t=-1$]  Criterium \ref{crit-A} for $(1,2;1)$ gives $\linethree{-1}{q}{-1}{q^2}{-1}$. 
This diagram is in $\hlist$ if and only if either it is finite (thus $q\in\Gf$), or 
it appears in \eqref{eq:parametric-line-3}, which implies that $q\in\G_3'\cup\G_4'$.
\item[$s=q^{2}$, $t=q^{-2}$]  Criterium \ref{crit-A} for $(1,2;1)$ gives $\linethree{-1}{q}{-1}{q^2}{q^{-2}}$. 
This diagram is in $\hlist$ if and only if either it is finite (in which case $q\in\Gf$), or 
it appears in \eqref{eq:parametric-line-3}: the three possibilities are $q^2=q^{-1}$, $q^2=q^{-2}$, $q^2=q^{-3}$, so $q\in\Gf$.
\item[$s=-q^{3}$, $t=-q^{-3}$]  Criterium \ref{crit-A} for $(1,2;1)$ gives $\linethree{-1}{q}{-1}{-q^3}{-q^{-3}}$. 
This diagram is in $\hlist$ if and only if either it is finite (in which case $q\in\Gf$), or 
it appears in \eqref{eq:parametric-line-3}: the possibilities are $-q^3=q^{-1}$, $-q^3=q^{-2}$, $-q^3=q^{-3}$, but $\bq_{|\{2,3,4\}}$ is not in $\hlist$ in any case, a contradiction with \cite{AnG-rank3}.
\end{description}

\medbreak

\noindent $\triangleright$ If $\bq_{|\I_3}$ is the second triangle in \eqref{eq:parametric-triangle}, rotated counter-clockwise, then the diagram of $\bq$ is
\begin{align*}
&\tadpole{-1}{q^{3}}{-1}{q^{-2}}{q}{q^{-1}}{t}{s} &&\text{for some }s,t\ne 1.
\end{align*}
Again, we reduce to study the case in which both $\bq_{|\{1,3,4\}}$ and $\bq_{|\{2,3,4\}}$ are parametric lines. 
Looking at $\bq_{|\{2,3,4\}}$, there are 2 possibilities: $s=q^{-1}$, $t\in\{q,-1\}$ and $q^2=-1$, so $q\in\G_4' \subset \Gf$.

\medbreak

\noindent $\triangleright$ If $\bq_{|\I_3}$ is the second triangle in \eqref{eq:parametric-triangle}, rotated clockwise, then the diagram is
\begin{align*}
&\tadpole{-1}{q^{-2}}{q}{q^{-1}}{-1}{q^{3}}{t}{s} &&\text{for some }s,t\ne 1.
\end{align*}
We reduce to study the case in which $\bq_{|\{1,3,4\}}$ and $\bq_{|\{2,3,4\}}$ are parametric lines. There are 6 possibilities for $\bq_{|\{2,3,4\}}$ as in the analysis of the first triangle:
\begin{description}[leftmargin=10pt]
\item[$s=q$, $t=-1$] $\bq_{|\{1,3,4\}}$ is in \eqref{eq:parametric-line-3} iff $q\in \G_6'$ (12th diagram) or $q\in \G_4'$ (15th).

\item[$s^2=q$, $t=-1$] $\bq_{|\{1,3,4\}}$ is in \eqref{eq:parametric-line-3} iff $s\in \G_7'$ (15th), but Criterium \ref{crit-A} for $(1,2;1)$ gives $\linethree{-s^{-2}}{s^4}{-1}{s}{-1}$, which is not in $\hlist$.

\item[$s^3=q$, $t=-1$] $\bq_{|\{1,3,4\}}$ is in \eqref{eq:parametric-line-3} iff $s\in \G_{10}'$ (15th), so $q\in\Gf$.

\item[$s=r^2$, $t=-r^{-1}$, $q=r^3$] $\bq_{|\{1,3,4\}}$ is in \eqref{eq:parametric-line-3} iff $r\in \G_{11}'$ (7th), but Criterium \ref{crit-A} for $(1,2;1)$ gives $\linethree{-r^{-3}}{r^6}{-1}{r^2}{-r^{-1}}$, which is not in $\hlist$.

\item[$s=q$, $t=q^{-1}$] $\bq_{|\{1,3,4\}}$ is in \eqref{eq:parametric-line-3} iff $q\in \G_4'$ (5th), $q\in \G_7'$ (10th diagram), $q\in \G_{10}'$ (12th diagram) or $q\in \G_6'$ (17th diagram). If $q\in \G_{7}'$, then Criterium \ref{crit-A} for $(1,2;1)$ gives a contradiction.

\item[$st=1$] $\bq_{|\{1,3,4\}}$ is in \eqref{eq:parametric-line-3} iff 
$q\in \G_6'$, $t=-1$ (1st), $t=q^3$ (5th), $t=q^6$ (10th), $t=q^9$ (12th) or $q\in \G_6'$ (17th). Criterium \ref{crit-A} for $(1,3;1)$ rules out cases 2, 3 and 4 as it gives a braiding of Cartan type which is not finite: $\linethree{q}{q^{-3}}{q^3}{q^{-3n}}{q^{3n}}$, $n\in\I_3$.
\end{description}

\medbreak

\noindent $\triangleright$ Finally, if $\bq_{|\I_3}$ is the third triangle in \eqref{eq:parametric-triangle}, then the diagram of $\bq$ is
\begin{align*}
&\tadpole{-1}{q^{-1}r^{-1}}{-1}{r}{-1}{q}{t}{s} &&\text{for some }s,t\ne 1.
\end{align*}
Up to exchange vertices $1$ and $2$ we may assume that $q\ne -1$. Assume that $\bq_{|\{1,3,4\}}$ is in $\hlist$: if it is finite, then $r\in\Gf$. Otherwise $\bq_{|\{1,3,4\}}$ is parametric and there are 5 possibilities:
\begin{description}[leftmargin=10pt]
\item[$s=q^{-1}$, $t=q$]  Criterium \ref{crit-A} for $(2,3;1)$ gives $\linethree{-1}{r^{-1}}{r}{q^{-1}}{q}$. If this line is finite, then $q,r\in\Gf$. If it is parametric, then either $r=-1$ and $\bq$ appears in row 14 of \cite[Table 3]{H-full}, or $q=r^a$ for some $a\in\I_3$: the cases $a=1,2$ belong to rows 13 and 9 of $\hlist$ respectively, and for $a=3$ we use Criterium \ref{crit-B} with $\alpha_1+\alpha_3$, $\alpha_2+\alpha_3$, $\alpha_4$ to obtain $\linethree{r}{r^{-3}}{r^3}{r^{-3}}{r^3}$, which is not in $\hlist$.
\item[$s=q^{-1}$, $q=t^2$]  Criterium \ref{crit-A} for $(3,4;1)$ gives $$\tri{-1}{t^{2}}{-t^{-1}}{r}{-1}{t^{-2}r^{-1}.}$$ 
This diagram is in $\hlist$ if and only if either it is finite (in which case $q,r\in\Gf$), or 
it appears in \eqref{eq:parametric-triangle}, in which case $r=-t$; it follows that $\bq_{|\{2,3,4\}}$ must be finite, so $q,r\in\Gf$.
\item[$s=q^{-a}$, $t=q^a$, $a\in\{2,3\}$] Criterium \ref{crit-B} for $\alpha_1+\alpha_3$, $\alpha_2+\alpha_3$, $\alpha_4$ gives $\linethree{q}{q^{-a}}{q^a}{q^{-a}}{r}$. If this line is finite, then $q,r\in\Gf$. Otherwise $a=2$, $r=q^2$, in which case $\bq$ is in row 9 of \cite[Table 3]{H-full}.
\item[$s=q^{-1}$, $t=-1$] Criterium \ref{crit-A} for $(3,4;1)$ gives 
$$\tri{-1}{q}{q^{-1}}{r}{-1}{q^{-1}r^{-1}.}$$ 
This diagram is in $\hlist$ if and only if either it is finite (hence $q,r\in\Gf$), or 
it appears in \eqref{eq:parametric-triangle}, in which case $r=q$ and $\bq$ appears in \cite[Table 3, row 9]{H-full}, or $r=q^2$: 
Criterium \ref{crit-A} for $(2,3;1)$ gives $\linethree{-1}{q^{-2}}{q^2}{q^{-1}}{-1}$, so $q\in\Gf$.
\end{description}

Therefore, either $\bq$ is in $\hlist$ or the parameters involved are in $\Gf$.
\epf

\begin{lemma}\label{lem:tadpole-with-param-line}
Let $\bq$ be a tadpole such that $\bq_{|\I_{2,4}}$ is parametric of type $\mathtt{G}$. Then either $\bq$ belongs to $\hlist$, or $\bq_{|\I_{2,4}}=\Gt{q}$ for some $q\in\Gf$, or $\bq_{|\I_{2,4}}=\Gt{q,r}$ for some $q,r\in\Gf$, or else $\gdim \B_{\bq}=\infty$.
\end{lemma}
\pf
Assume that $\gdim \B_{\bq}<\infty$. We consider the different possibilities for the triangle $\bq_{|\I_3}$ (necessarily in $\hlist$) and the parametric subdiagram $\bq_{|\I_{2,4}}$.

If $\bq_{|\I_3}$ is finite and $\bq_{|\I_{2,4}}=\Gt{q}$ for some $q\in\Bbbk^{\times}$, then $\qti{23}\in \G_6$ by Remark \ref{rem:rk-3-finite-in-Gf} and at the same time it coincides with $q$, $q^2$, or $q^3$; thus $q\in\Gf$. 

Next we assume that $\bq_{|\I_3}$ is parametric $\Gtp{r}$ and $\bq_{|\I_{2,4}}=\Gt{q}$ for some $q, r\in\Bbbk^{\times}$. 
By Lemma \ref{lem:tadpole-with-param-triang}, either $\bq$ is in $\hlist$ or else $r\in\Gf$. Hence we are done when
$\pm q^{\pm 1}\in \{q_{22},\qti{23},q_{33}\}$. Otherwise, one of the following holds:
\begin{multicols}{2}
\begin{itemize}[leftmargin=*]
\item $q_{22}=q_{33}=q^2$, $\qti{23}=q^{-2}$; 
\item $q_{22}=-1$, $q_{33}=q^2$, $\qti{23}=q^{-2}$;  
\item $q_{22}=q_{33}=-1$, $\qti{23}=q^2$; 
\item $q_{22}=q^2$, $q_{33}=-1$, $\qti{23}=q^{-2}$;  
\item $q_{22}=q^3$, $q_{33}=-1$, $\qti{23}=q^{-3}$.
\end{itemize}
\end{multicols}
Assume that $\bq$ is not in $\hlist$. Criterium \ref{crit-A} for $(3,4;1)$ gives either
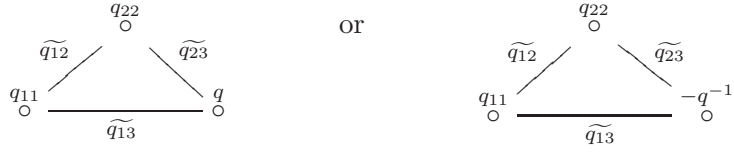
\begin{align*}
&\tri{q_{11}}{\qti{13}}{q}{\qti{23}}{q_{22}}{\qti{12}} &
&\text{or} &
&\tri{q_{11}}{\qti{13}}{-q^{-1}}{\qti{23}}{q_{22}}{\qti{12}}
\end{align*}
In any case, if the diagram is finite, then $q\in\G_6\subset \Gf$. Otherwise, it is one-parametric (since it has a vertex $\ne -1$), with parameter $\pm r^{\pm 1}$. Thus $q=\pm r^{\pm 1}\in\Gf$.

Hence it remains to check the case in which $\bq_{|\I_{2,4}}=\Gt{q,r}$, the last diagram in \eqref{eq:parametric-line-3}. Assume that $\bq$ is not in $\hlist$. By Lemma \ref{lem:tadpole-with-param-triang}, $q\in \Gf$ since it coincides with $q_{22}$. Set $\rho_3 \bq=(r_{ij})_{i,j\in\I_4}$: by direct computation,
\begin{align*}
\rti{12}&=1, &
\rti{14}&=\qti{13}r^{-1}, &
\rti{24}&=q^{-1}r^{-1}\ne 1, &
\rti{3i}&=\qti{3i}^{-1}\ne 1, \, i\ne 3.
\end{align*}
This forces $\qti{13}r=1$ since $\rho_3 \bq$ has no 4-cycles by Lemma \ref{lem:4-cycles}. Thus $\rho_3 \bq$ is a tadpole, where the triangle is $\rho_3\bq_{|\I_{2,4}}$, the last one in \eqref{eq:parametric-triangle}. Notice that $\rho_3\bq$ is not in $\hlist$ since we assume that $\bq$ is not in $\hlist$. By Lemma \ref{lem:tadpole-with-param-triang}, $r=\rti{34}\in\Gf$.
\epf

\begin{lemma}\label{lem:tripod-1param}
Let $\bq$ be a tripod such that $\bq_{|\I_3}$ is parametric of type $\mathtt{G}$. Then either $\bq$ belongs to $\hlist$, or $\bq_{|\I_3}=\Gt{q}$ for some $q\in\Gf$, or $\bq_{|\I_3}=\Gt{q,r}$ for some $q,r\in\Gf$, or else $\gdim \B_{\bq}=\infty$.
\end{lemma}
\pf
Assume that $\gdim \B_{\bq}<\infty$. If $\bq_{|\I_3}=\Gt{q,r}$, then $\rho_3 \bq=(r_{ij})_{i,j\in\I_4}$ is such that
\begin{align*}
\rti{12}&=q^{-1}r^{-1}\ne 1, & 
\rti{14}&=q^{-1}\qti{34}, &
\rti{24}&=r^{-1}\qti{34}, &
\rti{3i}&=\qti{3i}^{-1}\ne 1, \, i\ne 3.
\end{align*}
Then $q=r=\qti{34}$, since $\rho_3 \bq$ has no 4-cycles by Lemma \ref{lem:4-cycles}. Thus the diagrams of $\bq$ and $\rho_3\bq$ are, respectively,
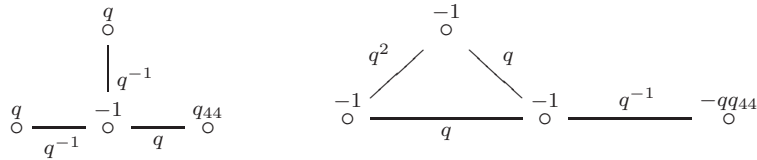
\begin{align*}
&\tripod{q}{1}{q}{q^{-1}}{-1}{q^{-1}}{q_{44}}{q} &
&\tadpole{-1}{q^2}{-1}{q}{-1}{q}{-qq_{44}}{q^{-1}}
\end{align*}
By Lemma \ref{lem:tadpole-with-param-triang}, either $\rho_3\bq$ is in $\hlist$, and so is $\bq$, or $q\in\Gf$.

\smallbreak
Next we assume that $\bq_{|\I_3}=\Gt{q}$ (one-parametric). If $\bq_{|\I_{2,4}}$ is finite, then either $q^{\pm 1}\in \{q_{22},q_{33},\qti{23}\}$, 
or $q^{\pm 2}\in \{q_{22},q_{33},\qti{23}\}$ (and for all these cases, $q\in\Gf$), or the diagram of $\bq$ is one of the following:
\begin{align*}
&\tripod{-1}{1}{q^3}{q^{-3}}{-1}{q}{t}{s} &
&\tripod{-q^{-1}}{1}{q^3}{q^{-3}}{-1}{q^2}{t}{s} &&q,s,t\in\Bbbk^{\times}.
\end{align*}
Set $\rho_3 \bq=(r_{ij})_{i,j\in\I_4}$. By Lemma \ref{lem:4-cycles}, $\qti{13}=\qti{23}=\qti{34}^{-1}$, which implies that $q\in\Gf$.
It remains to consider the case $\bq_{|\I_{2,4}}$ parametric. 
If $3$ is not Cartan, then $q_{33}=-1$. By Lemma \ref{lem:4-cycles}, $\rho_3\bq$ cannot contain a $4$-cycle. If $\rho_3\bq$ is a tripod, then $\qti{i3}=-1$ for all $i\in\{1,2,4\}$, so $\bq_{|\I_3}=\Gt{q}$ for some $q\in \Gf$. Otherwise $\rho_3\bq$ is a tadpole and thus there is a single pair $i<j\in\{1,2,4\}$ such that $\qti{i3}\qti{j3}\ne 1$. By Lemmas \ref{lem:tadpole-with-param-triang} and \ref{lem:tadpole-with-param-line}, either $\bq$ is in $\hlist$ or $q\in \Gf$.
\epf

\begin{lemma}\label{lem:line-1param}
Let $\bq$ be a line such that $\bq_{|\I_3}$ is parametric of type $\mathtt{G}$. Then either $\bq$ belongs to $\hlist$, or $\bq_{|\I_3}=\Gt{q}$ for some $q\in\Gf$, or $\bq_{|\I_3}=\Gt{q,r}$ for some $q,r\in\Gf$, or else $\gdim \B_{\bq}=\infty$.
\end{lemma}
\pf
Assume that $\gdim \B_{\bq}<\infty$. If $\bq_{|\I_3}=\Gt{q,r}$, then $\rho_2\bq$ has diagram
$$\tadpole{-1}{q}{-1}{r}{-1}{q^{-1}r^{-1}}{q_{44}}{\qti{34}}$$
and the result follows by Lemma \ref{lem:tadpole-with-param-triang}.

Now we deal with the case $\bq_{|\I_3}=\Gt{q}$ and $\bq_{|\I_{2,4}}$ finite. If one of the scalars $q_{22}$, $\qti{23}$, $q_{33}$ is either $q$ or $q^2$, then $q\in\Gf$. Otherwise, 
$q_{22}=-1$, $\qti{23}=q^{-3}$, $q_{33}=q^3$, $\rho_2\bq$ is a tadpole and Lemma \ref{lem:tadpole-with-param-triang} applies.

Finally, assume that $\bq_{|\I_3}=\Gt{q}$ and $\bq_{|\I_{2,4}}$ is parametric. By Theorem \ref{thm:Cartan} we may assume that at least one vertex $i\in\I_4$ is not Cartan, so $q_{ii}=-1$. By symmetry we may assume that $i\in\I_2$. If $q_{11}=-1$, then $\rho_1\bq=(r_{ij})$ is a line with $r_{22}=-1$. 
Hence we fix $i=2$. If $\qti{12}\qti{23}\ne 1$, then $\rho_2\bq$ is a tadpole and Lemma \ref{lem:tadpole-with-param-triang} applies, so we also assume $\qti{12}\qti{23}=1$. As $\bq_{|\I_3}$ is in \eqref{eq:parametric-line-3}, we have that either $q_{33}=-1$, $q_{33}\qti{23}=1$ or $q_{33}^2\qti{23}=1$. If $q_{33}=-1$, then we have three possibilities:
\begin{itemize}[leftmargin=*]
\item $\qti{23}\qti{34}\ne 1$. Then $\rho_3\bq$ is a tadpole, so Lemma \ref{lem:tadpole-with-param-line} applies.
\item $\qti{23}\qti{34}=1$, $q_{44}^2\qti{34}=1$. If $q_{11}\in\{-1,\qti{12}^{-1}\}$, then $\bq$ is in $\hlist$. Otherwise $q_{11}^2\qti{12}=1$, which is a contradiction as Criterium \ref{crit-A} for $(3,4;1)$ gives a diagram not in $\hlist$.
\item $\qti{23}\qti{34}=1$, $q_{44} \in\{-1,\qti{34}^{-1}\}$. Here, $\bq$ is in $\hlist$.
\end{itemize}
If $q_{33}\qti{23}=1$, then $\rho_2\bq=(r_{ij})$ is a line with $r_{22}=r_{33}=-1$ and this is the case just analysed. Finally, when $q_{33}^2\qti{23}=1$, necessarily $q_{33}=q$ and looking at the possibilities for $\bq_{|\I_{2,4}}$ in \eqref{eq:parametric-line-3} we have that $q^2=-1$.
\epf

\subsection{Proof of the conjecture in rank 4}

\begin{theorem}\label{thm:rank4-conj-true}
A rank 4 matrix $\bq$ is in $\hlist$ if and only if $\gdim \B_{\bq}<\infty$.
\end{theorem}
\pf
Let $\bq=(q_{ij})_{i,j\in\I_4}$ be such that $\gdim \B_{\bq}<\infty$. By Theorem \ref{prop:n-cycles} the diagram does not contain a 4-cycle, so 
there are three possible shapes:
\begin{description}[leftmargin=*]
\item[Tadpole] If one of the connected subdiagrams of rank three is a parametric $\Gt{q}$ or $\Gt{q,r}$ with $q,r\notin \Gf$, then $\bq$ is in $\hlist$ by Lemmas \ref{lem:tadpole-with-param-triang} and \ref{lem:tadpole-with-param-line}. Otherwise the three subdiagrams are finite or parametric evaluated in $\Gf$ and we apply \S \ref{gap:rank4t} to get that $\bq$ is in $\hlist$.
\item[Tripod] Similarly, if one of the connected subdiagrams of rank three is parametric not evaluated in $\Gf$, then $\bq$ is in $\hlist$ by Lemma \ref{lem:tripod-1param}; otherwise $\bq$ is in $\hlist$ by \S \ref{gap:rank4t}.
\item[Line] This follows by either Lemma \ref{lem:line-1param} or else \S \ref{gap:rank4t}.
\end{description}
In any case, $\bq$ is in $\hlist$.
\epf

\section{Rank 5}\label{sec:rank5}

Now we prove Conjecture \ref{conjecture} for matrices $\bq$ of rank 5, with connected diagram. 
To start we analyse the admissible shapes for the diagram of $\bq$. First it cannot contain cycles by Proposition \ref{prop:n-cycles}.
We also show that no vertex can have degree bigger than three. This leads to the four possible shapes depicted in \eqref{eq:rank5-possible-shapes}.
We show that the same trichotomy from \S \ref{sec:rank4} holds and apply \gap\, to deal with the finite set of diagrams left.

\subsection{Vertices of high degree}\label{subsec:highdegree}

We begin by stating the result that restricts the degree of the vertices of those $\bq$ such that $\gdim\B_{\bq}<\infty$.

\begin{remark}\label{rem:rank4-parametric-tadpoles}
There are seven parametric rank four tadpoles in $\hlist$, they belong to Rows 8, 9, 12, 13 and 14.
We observe that for each one of them, $m_{ij}\{0,1\}$ for all $i\ne j \in\I_4$.

Similarly, there are four parametric tripods in $\hlist$, in Rows 5, 12 and 13, and again $m_{ij}\{0,1\}$ for all $i\ne j \in\I_4$.
\end{remark}

\begin{proposition}\label{prop:rk5-degree-four}
If the diagram of $\bq$ has a vertex of degree bigger than 3, then $\gdim\B_{\bq}=\infty$.
\end{proposition}
\pf
By relabelling the vertices we may assume that vertex $2$ has degree $\ge 4$, and $\qti{2j}\neq 1$ for all $j\in\I_5$, $j\ne 2$. Now it is enough to prove the result for $\bq_{|\I_5}$. If the diagram contains a 4-cycle, then $\gdim\B_{\bq}=\infty$, by Lemma \ref{lem:4-cycles}. Thus we may assume that $\qti{13}=\qti{45}=1$. 
Hence the diagram can only be of one of the following three shapes:
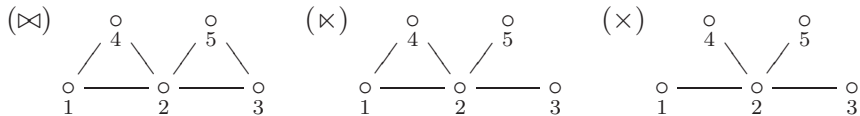
\begin{align*}
(\bowtie) & 
\xymatrix@C-18pt@R-18pt{ & \underset{4}{\circ} \ar @{-}[ld] \ar @{-}[rd] & & \underset{5}{\circ} \ar @{-}[ld] \ar @{-}[rd] & \\
\underset{1}{\circ} \ar @{-}[rr] & & \underset{2}{\circ} \ar @{-}[rr] & & \underset{3}{\circ}}
&  (\ltimes)& 
\xymatrix@C-18pt@R-18pt{ & \underset{4}{\circ} \ar @{-}[ld] \ar @{-}[rd] & & \underset{5}{\circ} \ar @{-}[ld]  & \\
\underset{1}{\circ} \ar @{-}[rr] & & \underset{2}{\circ} \ar @{-}[rr] & & \underset{3}{\circ}}
&  (\times)& 
\xymatrix@C-18pt@R-18pt{ & \underset{4}{\circ} \ar @{-}[rd] & & \underset{5}{\circ} \ar @{-}[ld]  & \\
\underset{1}{\circ} \ar @{-}[rr] & & \underset{2}{\circ} \ar @{-}[rr] & & \underset{3}{\circ}} & 
\end{align*}

We analyze the case $(\bowtie)$. We start by showing that if the tadpole $\bq_{|\I_4}$ is parametric $\Gt{q}$, then  $\gdim\B_{\bq}=\infty$ or $q\in\Gf$. The same holds if any of the submatrices $\bq_{|\I_5\setminus j}$, $j=1,3,4$, is a parametric tadpole. 

If $\bq_{|\I_4}$ is in Row 8, then the diagram of $\bq$ is:
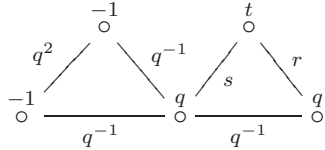
\begin{align*}
&\xymatrix@C-10pt@R-10pt{ & \overset{-1}{\circ} \ar @{-}[dl]_{q^2} \ar @{-}[dr]^{q^{-1}} & & \overset{t}{\circ} \ar @{-}[dr]^{r} \ar @{-}[dl]^{s} & \\
\overset{-1}{\circ} \ar @{-}[rr]_{q^{-1}} & & \overset{q}{\circ} \ar @{-}[rr]_{q^{-1}}  & & \overset{q}{\circ}}
\end{align*}
for some $r,s,t\ne 0,1$. As $\bq_{|\{2,3,5\}}$ is a triangle with two of the labels of the vertices $\ne -1$, $\bq_{|\{2,3,5\}}$ is finite since it does not appear in \eqref{eq:parametric-triangle}, thus $q_{22}=q\in\Gf$.

If $\bq_{|\I_4}$ is the first tadpole in Row 9, then the diagram of $\bq$ is:
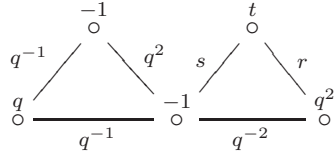
\begin{align*}
&\xymatrix@C-10pt@R-10pt{ & \overset{-1}{\circ} \ar @{-}[dl]_{q^{-1}} \ar @{-}[rd]^{q^2} & & \overset{t}{\circ} \ar @{-}[dr]^{r} \ar @{-}[dl]_{s} & \\
\overset{q}{\circ} \ar @{-}[rr]_{q^{-1}}  & & \overset{-1}{\circ} \ar @{-}[rr]_{q^{-2}} & & \overset{q^2}{\circ}}
\end{align*}
for some $r,s,t\ne 0,1$. Criterium \ref{crit-A} for $(2,4;1)$ gives
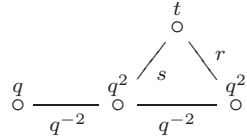
\begin{align*}
&\xymatrix@C-15pt@R-15pt{ & & & \overset{t}{\circ} \ar @{-}[dr]^{r} \ar @{-}[dl]^{s} & \\
	\overset{q}{\circ} \ar @{-}[rr]_{q^{-2}}  & & \overset{q^2}{\circ} \ar @{-}[rr]_{q^{-2}} & & \overset{q^2}{\circ}}
\end{align*}
which cannot be one of the seven parametric tadpoles by Remark \ref{rem:rank4-parametric-tadpoles}, so it should be finite and thus $q\in\Gf$.

If $\bq_{|\I_4}$ is the second tadpole in Row 9, then the diagram of $\bq$ is:
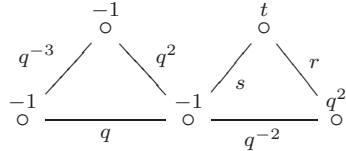
\begin{align*}
&\xymatrix@C-10pt@R-10pt{ & \overset{-1}{\circ} \ar @{-}[dl]_{q^{-3}} \ar @{-}[rd]^{q^2} & & \overset{t}{\circ} \ar @{-}[dr]^{r} \ar @{-}[dl]^{s} & \\
\overset{-1}{\circ} \ar @{-}[rr]_{q}  & & \overset{-1}{\circ} \ar @{-}[rr]_{q^{-2}} & & \overset{q^2}{\circ}}
\end{align*}
for some $r,s,t\ne 0,1$. Criterium \ref{crit-A} for $(1,2;1)$ gives a necessarily non-parametric tadpole as above, so it should be finite and thus $q\in\Gf$.

If $\bq_{|\I_4}$ is in Row 12, then the diagram of $\bq$ is one of the following:
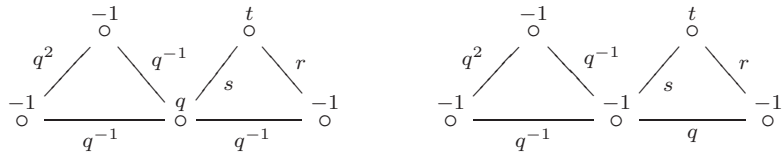
\begin{align*}
&\xymatrix@C-10pt@R-10pt{ & \overset{-1}{\circ} \ar @{-}[dl]_{q^2} \ar @{-}[dr]^{q^{-1}} & & \overset{t}{\circ} \ar @{-}[dr]^{r} \ar @{-}[dl]^{s} & \\
	\overset{-1}{\circ} \ar @{-}[rr]_{q^{-1}} & & \overset{q}{\circ} \ar @{-}[rr]_{q^{-1}}  & & \overset{-1}{\circ} }
&\xymatrix@C-10pt@R-10pt{ & \overset{-1}{\circ} \ar @{-}[dl]_{q^2} \ar @{-}[dr]^{q^{-1}} & & \overset{t}{\circ} \ar @{-}[dr]^{r} \ar @{-}[dl]^{s} & \\
	\overset{-1}{\circ} \ar @{-}[rr]_{q^{-1}} & & \overset{-1}{\circ} \ar @{-}[rr]_{q}  & & \overset{-1}{\circ} }
\end{align*}
for some $r,s,t\ne 0,1$. In the first case, Criterium \ref{crit-A} for $(1,4;1)$ gives again a non-parametric tadpole, so $q\in\Gf$. 
In the second case, Criterium \ref{crit-A} for $(2,4;1)$ gives a finite tadpole, which says that $q\in\Gf$, or the first tadpole in row 12, which implies $r=q^{-2}$, $s=q$, $t=-1$. In this case, Criterium \ref{crit-A} for $(1,4;1)$ necessarily gives non-parametric tadpole, so $q\in\Gf$.

If $\bq_{|\I_4}$ is in Row 13, then the diagram of $\bq$ is:
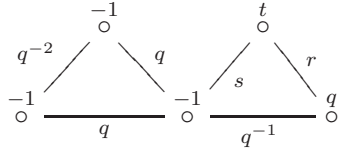
\begin{align*}
&\xymatrix@C-10pt@R-10pt{ & \overset{-1}{\circ} \ar @{-}[dl]_{q^{-2}} \ar @{-}[dr]^{q} & & \overset{t}{\circ} \ar @{-}[dr]^{r} \ar @{-}[dl]^{s} & \\
	\overset{-1}{\circ} \ar @{-}[rr]_{q} & & \overset{-1}{\circ} \ar @{-}[rr]_{q^{-1}}  & & \overset{q}{\circ}}
\end{align*}
for some $r,s,t\ne 0,1$. Criterium \ref{crit-A} for $(2,4;1)$ gives a tadpole with two vertices $\ne -1$ in the triangle, so it should be finite and thus $q\in\Gf$.

Assume that $\bq_{|\I_4}$ is in Row 14. Then the diagram of $\bq$ is as follows:
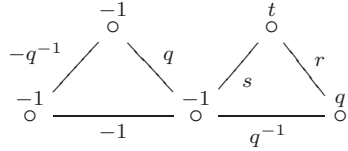
\begin{align*}
&\xymatrix@C-10pt@R-10pt{ & \overset{-1}{\circ} \ar @{-}[dr]^{q} \ar @{-}[dl]_{-q^{-1}} & & \overset{t}{\circ} \ar @{-}[dl]^{s} \ar @{-}[rd]^{r} & \\
	\overset{-1}{\circ} \ar @{-}[rr]_{-1} & & \overset{-1}{\circ} \ar @{-}[rr]_{q^{-1}}  & & \overset{q}{\circ}}
\end{align*}
for some $r,s,t\ne 0,1$. Criterium \ref{crit-A} for $(1,4;1)$ gives
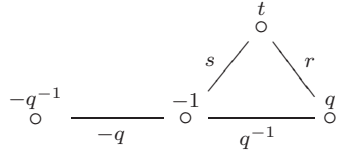
\begin{align*}
&\xymatrix@C-10pt@R-10pt{ &  & & \overset{t}{\circ} \ar @{-}[dl]_{s} \ar @{-}[rd]^{r} & \\
	\overset{-q^{-1}}{\circ} \ar @{-}[rr]_{-q} & & \overset{-1}{\circ} \ar @{-}[rr]_{q^{-1}}  & & \overset{q}{\circ}}
\end{align*}
which is the first tadpole in Row 9 with $q\in\G_6'$ or finite. In any case, $q\in\Gf$.

We are left with the setting in which all tadpoles are finite or parametric evaluated in $\Gf$. We use \gap\ to show that $\gdim\B_{\bq}=\infty$ in \ref{gap:valence4}. 

\medbreak

Now we turn to $(\ltimes)$. If $\bq_{|\I_4}$ (or analogously $\bq_{|\I_5\setminus 3}$) is a parametric tadpole $\Gt{q}$, then we prove that $\gdim\B_{\bq}=\infty$ or the parameter $q$ is in $\Gf$.

If $\bq_{|\I_4}$ is in Row 8, then the diagram of $\bq$ is:
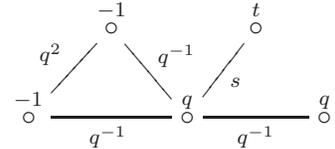
\begin{align*}
&\xymatrix@C-10pt@R-10pt{ & \overset{-1}{\circ} \ar @{-}[dl]_{q^2} \ar @{-}[dr]^{q^{-1}} & & \overset{t}{\circ}  \ar @{-}[dl]^{s} & \\
	\overset{-1}{\circ} \ar @{-}[rr]_{q^{-1}} & & \overset{q}{\circ} \ar @{-}[rr]_{q^{-1}}  & & \overset{q}{\circ}}
\end{align*}
for some $s,t\ne 0,1$. If  $\bq_{|\I_5\setminus 3}$ is finite, then $q\in\Gf$. Otherwise $\bq_{|\I_5\setminus 3}$ is parametric, so $s=q^{-1}$, $t\in\{q,-1\}$.
If $t=-1$, then $\rho_2\rho_5\bq=(r_{ij})$ is such that $\rti{13}=\rti{34}=q^{-2}\ne 1$, thus its diagram contains a 4-cycle and we apply Proposition \ref{prop:n-cycles}. If $t=q$, then Criterium \ref{crit-A} for $(1,4;1)$ gives a braiding of affine Cartan type and we apply Theorem \ref{thm:Cartan}.

If $\bq_{|\I_4}$ is the tadpole in Row 9, then the diagram of $\bq$ is one of the following:
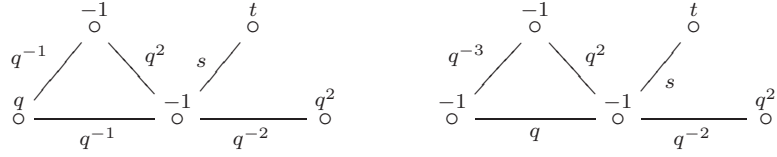
\begin{align*}
&\xymatrix@C-10pt@R-10pt{ & \overset{-1}{\circ} \ar @{-}[dl]_{q^{-1}} \ar @{-}[rd]^{q^2} & & \overset{t}{\circ}  \ar @{-}[dl]_{s} & \\
	\overset{q}{\circ} \ar @{-}[rr]_{q^{-1}}  & & \overset{-1}{\circ} \ar @{-}[rr]_{q^{-2}} & & \overset{q^2}{\circ}}
& &\xymatrix@C-10pt@R-10pt{ & \overset{-1}{\circ} \ar @{-}[dl]_{q^{-3}} \ar @{-}[rd]^{q^2} & & \overset{t}{\circ}  \ar @{-}[dl]^{s} & \\
	\overset{-1}{\circ} \ar @{-}[rr]_{q}  & & \overset{-1}{\circ} \ar @{-}[rr]_{q^{-2}} & & \overset{q^2}{\circ}}
\end{align*}
for some $s,t\ne 0,1$. Set $\rho_2\bq=(r_{ij})_{i,j\in\I_5}$. In the first case, if $s\ne q,q^2$, then $\rti{15},\rti{35}\ne 1$; if $s=q \ne q^{-2}$, then $\rti{35},\rti{45}\ne 1$; if $s=q^2 \ne q^{-2}$, then $\rti{15},\rti{45}\ne 1$. In all these cases $\rho_2\bq$ contains a 4-cycle, so $\gdim\B_{\bq}=\infty$ by Proposition \ref{prop:n-cycles}. Otherwise, $q\in\G_3'\cup \G_4'\subset \Gf$. The proof for the second diagram is analogous.

If $\bq_{|\I_4}$ is in Row 12, then the diagram of $\bq$ is one of the following:
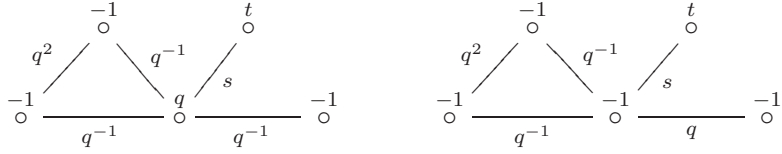
\begin{align*}
&\xymatrix@C-10pt@R-10pt{ & \overset{-1}{\circ} \ar @{-}[dl]_{q^2} \ar @{-}[dr]^{q^{-1}} & & \overset{t}{\circ}  \ar @{-}[dl]^{s} & \\
	\overset{-1}{\circ} \ar @{-}[rr]_{q^{-1}} & & \overset{q}{\circ} \ar @{-}[rr]_{q^{-1}}  & & \overset{-1}{\circ} }
&\xymatrix@C-10pt@R-10pt{ & \overset{-1}{\circ} \ar @{-}[dl]_{q^2} \ar @{-}[dr]^{q^{-1}} & & \overset{t}{\circ}  \ar @{-}[dl]^{s} & \\
	\overset{-1}{\circ} \ar @{-}[rr]_{q^{-1}} & & \overset{-1}{\circ} \ar @{-}[rr]_{q}  & & \overset{-1}{\circ} }
\end{align*}
for some $s,t\ne 0,1$. Observe that $\rho_3$ interchanges these diagrams, hence it is enough to prove the claim for the second. 
Set $\rho_2\bq=(r_{ij})_{i,j\in\I_5}$. If $s\ne q$, then $\rti{15},\rti{45}\ne 1$, so $\gdim\B_{\bq}=\infty$ by Proposition \ref{prop:n-cycles}. If $s=q$, then $r_{33}=q$, $\rti{35}=q^2$, so either $\bq_{|\{1,2,3,5\}}$ is a finite tadpole or else $q^3=1$ by Remark \ref{rem:rank4-parametric-tadpoles}. In any case $q\in\Gf$.

If $\bq_{|\I_4}$ is in Row 13, respectively 14, then the diagram of $\bq$ is:
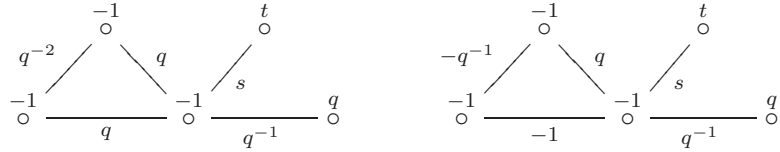
\begin{align*}
&\xymatrix@C-10pt@R-10pt{ & \overset{-1}{\circ} \ar @{-}[dl]_{q^{-2}} \ar @{-}[dr]^{q} & & \overset{t}{\circ}  \ar @{-}[dl]^{s} & \\
	\overset{-1}{\circ} \ar @{-}[rr]_{q} & & \overset{-1}{\circ} \ar @{-}[rr]_{q^{-1}}  & & \overset{q}{\circ}}
& &\xymatrix@C-10pt@R-10pt{ & \overset{-1}{\circ} \ar @{-}[dr]^{q} \ar @{-}[dl]_{-q^{-1}} & & \overset{t}{\circ} \ar @{-}[dl]^{s} & \\
	\overset{-1}{\circ} \ar @{-}[rr]_{-1} & & \overset{-1}{\circ} \ar @{-}[rr]_{q^{-1}}  & & \overset{q}{\circ}}
\end{align*}
for some $s,t\ne 0,1$. The proof for the first diagram is analogous to that for Row 12, while the second follows as for Row 9.

In the case in which the tripod $\bq_{|\I_5\setminus4}$ (or analogously $\bq_{|\I_{2,5}}$) is a parametric diagram $\Gt{q}$, either one of tadpoles is finite, in which case $q\in\Gf$ since $q^{\pm 1}$ appears in one of the edges shared between the tadpole and the tripod, or else both tadpoles are parametric: by the proof above, $\gdim\B_{\bq}=\infty$ or the tadpole is evaluated in $\Gf$ and thus $q\in\Gf$. 

Hence all tadpoles and tripods are finite or parametric evaluated in $\Gf$ and we use \gap\ to show that $\gdim\B_{\bq}=\infty$ in \ref{gap:valence4}.

\medbreak

Finally, in the case $(\times)$, we show that if $\bq_{|\I_4}$ is a parametric tripod $\Gt{q}$, then  $\gdim\B_{\bq}=\infty$ or $q\in \Gf$. The same holds if $\bq_{|\I_5\setminus j}$, $j=1,3,4$, is a parametric tripod. If one of the other tripods is finite, then $q\in\Gf$, so we assume that all tripods are parametric. If the four tripods are of Cartan $D_4$ type, then $\bq$ is of Cartan $D_4^{(1)}$ type, so $\gdim\B_{\bq}=\infty$ by Theorem \ref{thm:Cartan}. If $q_{22}=q\ne -1$, then the diagram of $\bq$ is, up to rotation,
\begin{align*}
&\xymatrix@C-10pt@R-10pt{ & \overset{q}{\circ} \ar @{-}[dr]^{q^{-1}} & & \overset{q}{\circ} \ar @{-}[dl]^{q^{-1}} & \\
	\overset{-1}{\circ} \ar @{-}[rr]_{q^{-1}} & & \overset{q}{\circ} \ar @{-}[rr]_{q^{-1}}  & & \overset{q}{\circ}}
\end{align*}
and the diagram of $\rho_2\rho_1\bq$ contains a 4-cycle. Otherwise $q_{22}=-1$, the four tadpoles are in Row 12 or 13 and the diagram of $\rho_2\bq$ contains a 4-cycle. In any case we apply Proposition \ref{prop:n-cycles} and $\gdim\B_{\bq}=\infty$.
We end up with two tripods which are either finite or parametric evaluated in $\Gf$. We discard this in \ref{gap:valence4}. 
\epf

\subsection{On parametric subdiagrams}\label{subsec:eval-param-rk5}

We begin by contemplating the possible shapes of the underlying graph. By Proposition \ref{prop:rk5-degree-four} all vertices have degree lesser or equal than 3. Observe that when all the vertices are of degree 1 and 2, then this corresponds to a line, as there can be no 5-cycles. 

Assume now that we have a vertex of degree 3. The rank 4 subdiagram formed by it and its neighbours it is either a tadpole or a tripod. 
This last case can only be extended to a graph of type $D_5$, since the missing vertex can only be joined to a single (degree 1) vertex of the tripod, to avoid a 4-cycle. 
Finally, if the subdiagram is a tadpole, then the missing vertex can only be joined to a single vertex, as above. Thus, we have two choices: to join it with the only vertex of degree 1 (obtaining a rank 5 tadpole) or a vertex of degree 2 (obtaining a diagram with a triangle in the middle).

We draw the four possible shapes we have obtained:
\begin{align}\label{eq:rank5-possible-shapes}
\begin{aligned}
& \xymatrix@C-18pt@R-18pt{  & &  &  &  & & \\
{\circ} \ar @{-}[rr] & & {\circ} \ar @{-}[rr] & & {\circ}\ar @{-}[rr] & & {\circ}  \ar @{-}[rr] & & {\circ} } &
& \xymatrix@C-18pt@R-18pt{ &  & & &{\circ}  \ar @{-}[d]   &  & \\
{\circ} \ar @{-}[rr] & & {\circ} \ar @{-}[rr] & & {\circ}\ar @{-}[rr] & & {\circ}   } \\
& \xymatrix@C-18pt@R-18pt{ & {\circ} \ar @{-}[ld] \ar @{-}[rd] & &  &  &  & \\
{\circ} \ar @{-}[rr] & & {\circ} \ar @{-}[rr] & & {\circ}\ar @{-}[rr] & & {\circ}   } &
& \xymatrix@C-18pt@R-18pt{& & & {\circ} \ar @{-}[ld] \ar @{-}[rd]   &  &  & \\
{\circ} \ar @{-}[rr] & &{\circ} \ar @{-}[rr] & & {\circ} \ar @{-}[rr] & & {\circ}  }
\end{aligned}
\end{align}

Next we prove again that the validity of the conjecture is reduced to checking a finite number of diagrams: those obtained by gluing either finite diagrams or parametric ones evaluated in $\Gf$, as in rank four.

\begin{lemma}\label{lemm:reduction-to-Gf-rank5}
Let $\bq=(q_{ij})_{i,j\in\I_5}$ be such that a connected subdiagram of rank four is parametric evaluated in $q\in\Bbbk^{\times}$. Then either $\bq$ is in $\hlist$, $\gdim \B_{\bq}=\infty$, or else $q\in\Gf$.
\end{lemma}
\pf
Assume that $\gdim \B_{\bq}<\infty$. Thus, by Theorem \ref{thm:rank4-conj-true}, all connected subdiagrams of rank 4 and those obtained by applying the criteria belong to $\hlist$. Also, the shape of the diagram is one of the four shapes in \eqref{eq:rank5-possible-shapes}. 
We also assume that all connected subdiagrams of rank 4 are parametric; indeed, if a subdiagram of rank four is finite, then it shares two edges with the parametric subdiagram of rank four, so $q\in\Gf$. 

Thus we are reduced to prove the following, for each possible shape: If all connected subdiagrams of rank 4 are parametric, then either $\bq$ is in $\hlist$ or else the parametric subdiagrams are evaluated in $\Gf$. For the last option, it is enough to prove that one of them is evaluated in an element in $\Gf$ since the connected subdiagrams share at least two edges.

\begin{enumerate}[leftmargin=*,label=\rm{(\Alph*)}]
\item\label{item:rank5-parametric-shapeA} The shape is 
$\xymatrix@C-18pt@R-18pt{& & & {\circ} \ar @{-}[ld] \ar @{-}[rd]   &  &  & \\
{\circ} \ar @{-}[rr] & &{\circ} \ar @{-}[rr] & & {\circ} \ar @{-}[rr] & & {\circ}  }$.
\end{enumerate}

If $\bq_{|\{1,2,3,5\}}$ is the first tadpole in Row 9, then the diagram of $\bq$ is:
\begin{align*}
&\xymatrix@C-15pt@R-15pt{ & & & \overset{-1}{\circ} \ar @{-}[dl]_{q^2} \ar @{-}[rd]^{q^{-1}} & \\
\overset{q^2}{\circ} \ar @{-}[rr]^{q^{-2}} & & \overset{-1}{\circ} \ar @{-}[rr]^{q^{-1}}  & & \overset{q}{\circ} \ar @{-}[rr]^{r} & & \overset{s}{\circ},}
&
&\xymatrix@C-15pt@R-15pt{ & & & \overset{q}{\circ} \ar @{-}[dl]_{q^{-1}} \ar @{-}[rd]^{q^{-1}} & \\
\overset{q^2}{\circ} \ar @{-}[rr]^{q^{-2}} & & \overset{-1}{\circ} \ar @{-}[rr]^{q^2}  & & \overset{-1}{\circ} \ar @{-}[rr]^{r} & & \overset{s}{\circ},}
\end{align*}
for some $r,s\ne 0,1$. Criterium \ref{crit-A} for $(2,3;1)$ gives, respectively,
\begin{align*}
&\xymatrix@C-10pt@R-10pt{ & & \overset{-1}{\circ} \ar @{-}[d]^{q} & & \\
\overset{q^2}{\circ} \ar @{-}[rr]^{q^{-2}} & & \overset{-1}{\circ} \ar @{-}[rr]^{r} & & \overset{s}{\circ},}
&
&\xymatrix@C-10pt@R-10pt{ & & \overset{q}{\circ} \ar @{-}[d]^{q^{-2}} & & \\
\overset{q^2}{\circ} \ar @{-}[rr]^{q^{-2}} & & \overset{q^2}{\circ} \ar @{-}[rr]^{r} & & \overset{s}{\circ}.}
\end{align*}
The second diagram is not parametric since the upper vertex is not -1 neither the inverse of the adjacent edge. As well, the first cannot be parametric as this would force to be in Row 13 with $q^2=q$, a contradiction. Thus, $q\in\Gf$.

If $\bq_{|\{1,2,3,5\}}$ is the second tadpole in Row 9, then the diagram of $\bq$ is:
\begin{align*}
&\xymatrix@C-15pt@R-15pt{ & & & \overset{-1}{\circ} \ar @{-}[dl]_{q^a} \ar @{-}[rd]^{q^{-3}} & \\
\overset{q^2}{\circ} \ar @{-}[rr]^{q^{-2}} & & \overset{-1}{\circ} \ar @{-}[rr]^{q^b}  & & \overset{-1}{\circ} \ar @{-}[rr]^{r} & & \overset{s}{\circ},}
& & \{a,b\}=\I_2, \, r,s\ne 0,1.
\end{align*}
Criterium \ref{crit-A} for $(2,3;1)$ gives 
$$ \xymatrix@C-10pt@R-10pt{ & & \overset{-1}{\circ} \ar @{-}[d]^{q^{-b}} & & 
\\ \overset{q^2}{\circ} \ar @{-}[rr]^{q^{-2}} & & \overset{q^b}{\circ} \ar @{-}[rr]^{r} & & \overset{s}{\circ}.}$$
If this diagram is finite, then $q\in\Gf$. If it is parametric, then either $q^2=-1$ (Row 5), $q^{b}=-1$ (Row 12 and the first diagram in Row 13), or else 
$b=2$, $r=q^{-2}$, $s=q^2$. In the last case, if $\bq_{\I_{2,5}}$ is in $\hlist$, then either $\bq_{\I_{2,5}}$ is finite or else $q^{-3}\in\{q,q^2,-1\}$. All in all, $q\in\Gf$.

If $\bq_{|\{1,2,3,5\}}$ is in Row 8, 12 or 13, then the diagram of $\bq$ is:
\begin{align*}
&\xymatrix@C-15pt@R-15pt{ & & & \overset{-1}{\circ} \ar @{-}[dr]^{q^2} \ar @{-}[dl]_{q^{-1}} & \\
\overset{v}{\circ} \ar @{-}[rr]^{u} & & \overset{t}{\circ} \ar @{-}[rr]^{q^{-1}}  & & \overset{-1}{\circ} \ar @{-}[rr]^{r} & & \overset{s}{\circ},}
&
&t\in\{-1, q\}, \, t^2u=q, \, v\in\{-1, u^{-1}\}.
\end{align*}
Criterium \ref{crit-A} for $(2,3;1)$ gives 
$$ \xymatrix@C-10pt@R-10pt{ & & \overset{-1}{\circ} \ar @{-}[d]^{q} & & 
\\ \overset{v}{\circ} \ar @{-}[rr]^{u} & & \overset{-tq^{-1}}{\circ} \ar @{-}[rr]^{r} & & \overset{s}{\circ}.}$$
\begin{itemize}[leftmargin=*]
\item This diagram does not belong to rows 5 neither 12 since $q\ne -1$. 
\item If this diagram is the first one of Row 13, then $t=v=s=q$, $u=r=q^{-1}$. Now $\bq_{\I_{2,5}}$ is a tadpole whose vertex of degree 3 is labeled with $-1$ and one of the vertices of the triangle labeled with $q\ne -1$. We check that $\bq_{\I_{2,5}}$ is not parametric, so it is finite, which implies that $q\in\Gf$.
\item If this diagram is the second one of Row 13, then $t=-1$, $v=s=q^{-1}$, $u=r=q$, so $\bq_{\I_{2,5}}$ is a tadpole whose vertices of the triangle are labeled with $q^{-1}$. If $\bq_{\I_{2,5}}$ belongs to Row 9 or Row 13, then $q^2=q^{-1}$. If $\bq_{\I_{2,5}}$ is in Row 14, then $q^2=-1$. Otherwise, $\bq_{\I_{2,5}}$ is finite. In any case, $q\in\Gf$.
\end{itemize}
Otherwise, the diagram above is finite, so $q\in\Gf$.

Finally, assume that $\bq_{|\{1,2,3,5\}}$ is in Row 14. Then the diagram of $\bq$ is:
\begin{align*}
&\xymatrix@C-15pt@R-15pt{ & & & \overset{-1}{\circ} \ar @{-}[dl]_{q} \ar @{-}[rd]^{-q^{-1}} & \\
	\overset{q}{\circ} \ar @{-}[rr]^{q^{-1}} & & \overset{-1}{\circ} \ar @{-}[rr]^{-1}  & & \overset{-1}{\circ} \ar @{-}[rr]^{r} & & \overset{s}{\circ},}
&
&\xymatrix@C-15pt@R-15pt{ & & & \overset{-1}{\circ} \ar @{-}[dl]_{-1} \ar @{-}[rd]^{-q^{-1}} & \\
	\overset{q}{\circ} \ar @{-}[rr]^{q^{-1}} & & \overset{-1}{\circ} \ar @{-}[rr]^{q}  & & \overset{-1}{\circ} \ar @{-}[rr]^{r} & & \overset{s}{\circ},}
\end{align*}
for some $r,s\ne 0,1$. Criterium \ref{crit-A} for $(2,3;1)$ gives, respectively,
\begin{align*}
&\xymatrix@C-10pt@R-10pt{ & & \overset{-1}{\circ} \ar @{-}[d]^{-1} & & \\
	\overset{q}{\circ} \ar @{-}[rr]^{q^{-1}} & & \overset{-1}{\circ} \ar @{-}[rr]^{r} & & \overset{s}{\circ},}
&
&\xymatrix@C-10pt@R-10pt{ & & \overset{-1}{\circ} \ar @{-}[d]^{q^{-1}} & & \\
	\overset{q}{\circ} \ar @{-}[rr]^{q^{-1}} & & \overset{q}{\circ} \ar @{-}[rr]^{r} & & \overset{s}{\circ}.}
\end{align*}
The first diagram does not belong to Rows 5, 12, 13 of \cite[Table 3]{H-full} since $q\ne 1$, so it is not parametric. 
The second is parametric if and only if $r=q^{-1}$, $s=q$ and $q^a=-1$ for $a\in\{2,3\}$.
Otherwise they must be finite. In any case, $q\in\Gf$.

\begin{enumerate}[leftmargin=*,label=\rm{(\Alph*)}]
\setcounter{enumi}{1}
\item\label{item:rank5-parametric-shapeB} The shape is 
$\xymatrix@C-18pt@R-18pt{ & {\circ} \ar @{-}[ld] \ar @{-}[rd] & &  &  &  & \\
{\circ} \ar @{-}[rr] & & {\circ} \ar @{-}[rr] & & {\circ}\ar @{-}[rr] & & {\circ}   }$.
\end{enumerate}

If $\bq_{|\{1,2,3,5\}}$ is in Row 9, then the diagram of $\bq$ is one of the following:
\begin{align*}
&\xymatrix@C-15pt@R-15pt{& \overset{-1}{\circ} \ar @{-}[dl]_{q^{-1}} \ar @{-}[rd]^{q^2} & & & \\
\overset{q}{\circ} \ar @{-}[rr]^{q^{-1}}  & & \overset{-1}{\circ} \ar @{-}[rr]^{q^{-2}} & & \overset{q^2}{\circ} \ar @{-}[rr]^{r} && \overset{s}{\circ},}
&
&\xymatrix@C-15pt@R-15pt{& \overset{-1}{\circ} \ar @{-}[dl]_{q^{-3}} \ar @{-}[rd]^{q^2} & & & \\
\overset{-1}{\circ} \ar @{-}[rr]^{q}  & & \overset{-1}{\circ} \ar @{-}[rr]^{q^{-2}} & & \overset{q^2}{\circ} \ar @{-}[rr]^{r} && \overset{s}{\circ},}
\end{align*}
where $r,s\ne 0,1$. The diagram of $\rho_2\bq$ is as in \ref{item:rank5-parametric-shapeA}, so either $\bq\in\hlist$ or $q\in\Gf$.

If $\bq_{|\{1,2,3,5\}}$ is in Row 8, 12 or 13, then the diagram of $\bq$ is:
\begin{align*}
&\xymatrix@C-15pt@R-15pt{& \overset{-1}{\circ} \ar @{-}[dl]_{q^2} \ar @{-}[rd]^{q^{-1}} & & & \\
\overset{-1}{\circ} \ar @{-}[rr]^{q^{-1}}  & & \overset{t}{\circ} \ar @{-}[rr]^{u} & & \overset{v}{\circ} \ar @{-}[rr]^{r} && \overset{s}{\circ},}
&
&t\in\{-1, q\}, \, t^2u=q, \, v\in\{-1, u^{-1}\}.
\end{align*}
We apply Criterium \ref{crit-A} for $(2,3;1)$ and we should obtain a diagram in $\hlist$. If this diagram is finite, then $q\in\Gf$. Otherwise, the diagram is a parametric tadpole, which must belong to Rows 8, 12 or 13; thus $v^2r=q$, $s\in\{-1, r^{-1}\}$, and $\bq$ belongs to \cite[Table 4 -- Rows 9, 10]{H-full}.

Assume that $\bq_{|\{1,2,3,5\}}$ is in Row 14. Then the diagram of $\bq$ is as follows:
\begin{align*}
&\xymatrix@C-15pt@R-15pt{& \overset{-1}{\circ} \ar @{-}[dl]_{-q^{-1}} \ar @{-}[rd]^{q} & & & \\
	\overset{-1}{\circ} \ar @{-}[rr]^{-1}  & & \overset{-1}{\circ} \ar @{-}[rr]^{q^{-1}} & & \overset{q}{\circ} \ar @{-}[rr]^{r} && \overset{s}{\circ},}
& & r,s\ne 0,1.
\end{align*}
Again, the diagram of $\rho_2\bq$ is as in \ref{item:rank5-parametric-shapeA}, so either $\bq\in\hlist$ or $q\in\Gf$.

\begin{enumerate}[leftmargin=*,label=\rm{(\Alph*)}]
\setcounter{enumi}{2}
\item\label{item:rank5-parametric-shapeC} The shape is 
$\xymatrix@C-18pt@R-18pt{ &  & & &{\circ}  \ar @{-}[d]   &  & \\
{\circ} \ar @{-}[rr] & & {\circ} \ar @{-}[rr] & & {\circ}\ar @{-}[rr] & & {\circ}   }$.
\end{enumerate}
Here, $\bq_{\I_{2,5}}$ belongs to Rows 5, 12 or 13. We apply Criterium \ref{crit-A} for $(2,3;1)$ and obtain a diagram in $\hlist$ whose shape is a tripod. If the diagram is finite, then $q\in \Gf$. Otherwise the diagram is parametric, so it belongs to rows 5, 12 or 13; thus $q_{11}\in\{\qti{12}^{-1}, -1\}$ and $q_{22}^2\qti{12}\qti{23}=1$ by Remark \ref{rem:rank4-parametric-tadpoles}, which implies that $\bq$ is in \cite[Table 4 -- Rows 8,10]{H-full}.

\begin{enumerate}[leftmargin=*,label=\rm{(\Alph*)}]
\setcounter{enumi}{3}
\item\label{item:rank5-parametric-shapeD} The shape is 
$\xymatrix@C-18pt@R-18pt{{\circ} \ar @{-}[rr] & & {\circ} \ar @{-}[rr] & & {\circ}\ar @{-}[rr] & & {\circ}  \ar @{-}[rr] & & {\circ} }$.
\end{enumerate}
Assume first that either $\bq_{|\I_{4}}$ or $\bq_{|\I_{2,5}}$ belongs to Rows 8, 9, 12, 13 or 14. Applying a suitable chain of reflections $\rho_i$ we obtain a matrix $\bq'$ whose shape is as in \ref{item:rank5-parametric-shapeA} or \ref{item:rank5-parametric-shapeB}; by the proofs of these cases, either the parameter $q$ belongs to $\Gf$ or else $\bq'$ is in $\hlist$, so the same holds for $\bq$.

Assume now that $\bq_{|\I_{4}}$  belongs to Rows 2, 7 or 11. The diagram of $\bq$ is
\begin{align*}
&\xymatrix@C-5pt{\overset{q_{11}}{\circ} \ar @{-}[r]^{\qti{12}} & \overset{q_{22}}{\circ} \ar @{-}[r]^{\qti{23}}  & \overset{q_{33}}{\circ} \ar @{-}[r]^{q^{-2}} & \overset{q}{\circ} \ar @{-}[r]^{r} & \overset{s}{\circ},}
& & \text{or} &
&\xymatrix@C-5pt{\overset{q}{\circ} \ar @{-}[r]^{q^{-2}} & \overset{q_{22}}{\circ} \ar @{-}[r]^{\qti{23}}  & \overset{q_{33}}{\circ} \ar @{-}[r]^{\qti{34}} & \overset{q_{44}}{\circ} \ar @{-}[r]^{r} & \overset{s}{\circ},}
\end{align*}
where either $q_{ii}=-1$ or else $q_{ii}\qti{i\,i\pm1}=1$, and $r,s\ne 0,1$. Set $t=q_{22}\qti{23}q_{33}$.
\begin{itemize}[leftmargin=*]
\item Criterium \ref{crit-A} for $(2,3;1)$ in the first diagram gives
$\xymatrix@C-5pt{\overset{q_{11}}{\circ} \ar @{-}[r]^{\qti{12}} & \overset{t}{\circ} \ar @{-}[r]^{q^{-2}} & \overset{q}{\circ} \ar @{-}[r]^{r} & \overset{s}{\circ}}$. If this diagram is finite, then $q\in\Gf$. Otherwise it is parametric, so $q_{11}=q^2$, $\qti{12}=q^{-2}$, $t=q^2$, $r=q^{-1}$ and either $s=q$ (Row 4) or $s=-1$ (Row 9): in the first case $\bq$ is Cartan of affine type, a contradiction with Theorem \ref{thm:Cartan}, while for the second $\bq_{|\I_{2,5}}$ belongs to Row 9 and was already analysed.
\item Criterium \ref{crit-A} for $(2,3;1)$ in the second diagram gives
$\xymatrix@C-5pt{\overset{q}{\circ} \ar @{-}[r]^{q^{-2}} & \overset{t}{\circ} \ar @{-}[r]^{\qti{34}} & \overset{q_{44}}{\circ} \ar @{-}[r]^{r} & \overset{s}{\circ}}$. If this diagram is finite, then $q\in\Gf$. Otherwise it is parametric, so it belongs to Rows 2, 7 or 11, with $q_{44}^2\qti{34}r=q^2$ and either $sr=1$ or $s=-1$; in any case, $\bq$ belongs to \cite[Table 4, Rows 3, 4]{H-full}.
\end{itemize}

Next we assume that $\bq_{|\I_{4}}$ belongs to Row 4. The diagram of $\bq$ is
\begin{align*}
&\xymatrix@C-5pt{\overset{q^2}{\circ} \ar @{-}[r]^{q^{-2}} & \overset{q^{2}}{\circ} \ar @{-}[r]^{q^{-2}}  & \overset{q}{\circ} \ar @{-}[r]^{q^{-1}} & \overset{q}{\circ} \ar @{-}[r]^{r} & \overset{s}{\circ},}
& & \text{or} &
&\xymatrix@C-5pt{\overset{q}{\circ} \ar @{-}[r]^{q^{-1}} & \overset{q}{\circ} \ar @{-}[r]^{q^{-2}}  & \overset{q^2}{\circ} \ar @{-}[r]^{q^{-2}} & \overset{q^2}{\circ} \ar @{-}[r]^{r} & \overset{s}{\circ}.}
\end{align*}
\begin{itemize}[leftmargin=*]\renewcommand{\labelitemi}{$\diamond$}
\item Criterium \ref{crit-A} for $(3,4;1)$ in the first diagram gives
$\xymatrix@C-5pt{\overset{q^2}{\circ} \ar @{-}[r]^{q^{-2}} & \overset{q^{2}}{\circ} \ar @{-}[r]^{q^{-2}}  & \overset{q}{\circ} \ar @{-}[r]^{r} & \overset{s}{\circ}}$. If this diagram is finite, then $q\in\Gf$. Otherwise it is parametric, so $q_{11}=q^2$, $\qti{12}=q^{-2}$, $t=q^2$, $r=q^{-1}$ and either $s=q$ (Row 4) or $s=-1$ (Row 9): we discard the first option by Theorem \ref{thm:Cartan}, while for the second $\bq_{|\I_{2,5}}$ belongs to Row 9 and was already considered.
\item Criterium \ref{crit-A} for $(3,4;1)$ in the second diagram gives
$\xymatrix@C-5pt{\overset{q}{\circ} \ar @{-}[r]^{q^{-1}} & \overset{q}{\circ} \ar @{-}[r]^{q^{-2}} & \overset{q^2}{\circ} \ar @{-}[r]^{r} & \overset{s}{\circ}}$. Suppose that this diagram is parametric; then it belongs to Row 4, i.e.~ $r=q^{-2}$, $s=q^2$, which implies that $\bq$ is Cartan of affine type, a contradiction with Theorem \ref{thm:Cartan}. Thus this diagram is finite, so $q\in\Gf$.
\end{itemize}

Next we assume that $\bq_{|\I_{4}}$ belongs to Row 3. The diagram of $\bq$ is
\begin{align*}
&\xymatrix@C-5pt{\overset{q}{\circ} \ar @{-}[r]^{q^{-1}} & \overset{q}{\circ} \ar @{-}[r]^{q^{-1}}  & \overset{q}{\circ} \ar @{-}[r]^{q^{-2}} & \overset{q^2}{\circ} \ar @{-}[r]^{r} & \overset{s}{\circ},}
& & \text{or} &
&\xymatrix@C-5pt{\overset{q^2}{\circ} \ar @{-}[r]^{q^{-2}} & \overset{q}{\circ} \ar @{-}[r]^{q^{-1}}  & \overset{q}{\circ} \ar @{-}[r]^{q^{-1}} & \overset{q}{\circ} \ar @{-}[r]^{r} & \overset{s}{\circ}.}
\end{align*}
\begin{itemize}[leftmargin=*]\renewcommand{\labelitemi}{$\heartsuit$}
\item Criterium \ref{crit-A} for $(2,3;1)$ in the first diagram gives 
$\xymatrix@C-5pt{\overset{q}{\circ} \ar @{-}[r]^{q^{-1}} & \overset{q}{\circ} \ar @{-}[r]^{q^{-2}} & \overset{q^2}{\circ} \ar @{-}[r]^{r} & \overset{s}{\circ}}$. Suppose that this diagram is parametric; then it belongs to Row 4, i.e.~ $r=q^{-2}$, $s=q^2$, which implies that $\bq$ is Cartan of affine type, a contradiction with Theorem \ref{thm:Cartan}. Thus this diagram is finite, so $q\in\Gf$.

\item Criterium \ref{crit-A} for $(2,3;1)$ in the second diagram gives
$\xymatrix@C-5pt{\overset{q^2}{\circ} \ar @{-}[r]^{q^{-2}} & \overset{q}{\circ} \ar @{-}[r]^{q^{-1}} & \overset{q}{\circ} \ar @{-}[r]^{r} & \overset{s}{\circ}}$. If this diagram is finite, then $q\in\Gf$. Otherwise it is parametric, so it belongs to Rows 3 or 8: 
$q_{44}^2\qti{34}r=q^2$ and either $sr=1$ or $s=-1$; that is, $\bq$ belongs to \cite[Table 4, Rows 7, 9]{H-full}.
\end{itemize}

The same holds if $\bq_{|\I_{2,5}}$ belongs to Rows 2, 3, 4, 7 or 11, up to reflection. Thus the remaining case is when both $\bq_{|\I_{4}}$ and $\bq_{|\I_{2,5}}$ belong to Rows either 1, 6 or 10. In this case, $\bq$ belongs to \cite[Table 4, Rows 1, 2]{H-full}.
\epf

\subsection{Proof of the conjecture in rank 5}

\begin{theorem}\label{thm:rank5-conj-true}
A rank 5 matrix $\bq$ is in $\hlist$ if and only if $\gdim \B_{\bq}<\infty$.
\end{theorem}
\pf
Assume that $\gdim \B_{\bq}<\infty$. If one of the connected subdiagrams of rank four is parametric not evaluated in $\Gf$, then $\bq$ is in $\hlist$ by Lemma
\ref{lemm:reduction-to-Gf-rank5}. Otherwise all subdiagrams are either finite or else parametric evaluated in $\Gf$ and we apply \S \ref{gap:criterian} to check that $\bq$ is in $\hlist$.
\epf

\section{Rank 6}\label{sec:rank6}

In this section we prove the conjecture in rank 6, following the same strategy as above.
For this section, we specially recall the Criterium \ref{crit-A} for $(i,j;1)$: we refer to its application saying that we {\it collapse} vertices $i$ and $j$.

\subsection{Some forbidden shapes}\label{subsec:forbidden6}
We start by discarding matrices whose diagram contains (up to deletion of edges) a graph of the following shape:
\begin{align}\label{eqn:forbidden}
\begin{aligned}
\xymatrix@C-18pt@R-18pt{& & & {\circ}  \ar @{-}[d]   &  &  & \\
& & & {\circ} \ar @{-}[ld] \ar @{-}[rd]   &  &  & \\
{\circ} \ar @{-}[rr] & &{\circ} \ar @{-}[rr] & & {\circ} \ar @{-}[rr] & & {\circ}  }
\end{aligned}
\end{align}

\begin{lemma}\label{lem:rk6-triangle-extended-3-vertices}
Assume there are three vertices $i,j,k$ such that $\qti{ij}, \qti{ik}, \qti{jk}\neq 1$ and their degrees are bigger or equal than three. Then $\gdim\B_{\bq}=\infty$.
\end{lemma}
\pf
Let $n(i)\neq j,k$, respectively, $n(j)\neq i,k$, $n(k)\neq i,j$ be a neighbour of $i$, respectively $j,k$. If two of these neighbours coincide, say $n(i)=n(j)$, then $\gdim\B_{\bq}=\infty$ by Lemma \ref{lem:4-cycles} as $i,j,k,n(i)$ determine a 4-cycle.

If $n(i), n(j), n(k)$ are pairwise different, then we restrict to the subdiagram of rank 6 with vertices $i,j,k$ and $n(i), n(j), n(k)$. If two of these neighbours are connected by an edge, then $\gdim\B_{\bq}=\infty$ by Lemma \ref{lem:4-cycles}. 

It thus remain to consider a diagram as in \eqref{eqn:forbidden}. The three subdiagrams of rank 5 obtained by removing one neighbour $n(i), n(j), n(k)$ are finite or $\gdim\B_{\bq}=\infty$ by Theorem \ref{thm:rank5-conj-true} as there are no parametric diagrams in $\hlist$ with a triangle in the middle. 
When these three subdiagrams are finite, we  use \gap\ to show that $\gdim\B_{\bq}=\infty$, see \ref{gap:star5}.
\epf

Next we show that when two vertices have degree bigger than two, they necessarily belong to a triangle in the middle of the diagram; in particular they have degree three.

\begin{lemma}\label{lem:rk6-2-vertices-of-degree-3}
Assume there are two vertices $i\ne j$ of degree bigger or equal than three. If $\gdim\B_{\bq}<\infty$, then the following conditions hold:
\begin{enumerate}[leftmargin=*,label=\rm{(\roman*)}]
\item\label{item:rk6-2-vertices-of-degree-3-i} The degree of $i$, respectively $j$, is 3 and $\qti{ij}\ne 1$.
\item\label{item:rk6-2-vertices-of-degree-3-ii} There exists $k\ne i,j$ joined to both $i$ and $j$. Moreover $k$ has no other neighbours except $i,j$.
\item\label{item:rk6-2-vertices-of-degree-3-iii} The degree of every $\ell\ne i,j$ is 1 or 2.
\end{enumerate}
\end{lemma}
\pf
We may assume that $i$ and $j$ are vertices of degree $\ge 3$ with minimal distance between them. By Proposition \ref{prop:rk5-degree-four} the degrees of $i$ and $j$ must be exactly 3. Also, $\theta \ge 5$, since otherwise $i$ and $j$ should share two neighbours, which is a contradiction with Lemma \ref{lem:4-cycles}. 

Let $\{m_1,m_2,m_3\}$, respectively $\{n_1,n_2,n_3\}$, be the set of neighbours of $i$, respectively $j$.
We argue by induction on $\theta$: the case $\theta=5$ is already considered in \S \ref{subsec:eval-param-rk5}. Fix $\theta\ge 6$.
Suppose that $m_s\ne j$ for all $s\in\I_3$, thus $\qti{ij}=1$. As the diagram is connected, there exists a path joining $i$ and $j$ of length $\ge 2$: 
we fix a path of minimal length and assume that this path starts with the edge $im_3$ followed by $m_3k$ for some $k\in\I-\{i\}$. Note that $m_1,m_2\notin \{n_1,n_2,n_3\}$: otherwise we also have $m_3\in \{n_1,n_2,n_3\}$ by the minimality of the length and the diagram has a $4$-cycle, a contradiction. Also notice that $m_3$ has degree two and $\qti{ik}=1$ by assumption on the distance between $i$ and $j$.

Collapsing $i$ and $m_3$ we obtain a matrix $\mathbf{r}=(r_{\ell m})$ of rank $\theta-1$, where the set of neighbours of the new vertex $\star=im_3$ is $\{m_1, m_2, k\}$. Indeed,
if $m_3$ is not connected with $m_1$ nor $m_2$, then 
\begin{align*}
\rti{\star k}&=\qti{ik}\qti{m_3k}=\qti{m_3k}\ne 1,  & \rti{\star m_s}&=\qti{im_s}q_{m_3m_s}=\qti{im_s}\ne 1, & s&=1,2,
\end{align*}
and $\rti{\star \ell}=1$ for any other $\ell$. Otherwise assume that $m_3$ is connected with $m_1$; hence $m_3$ is not connected with $m_2$ and 
\begin{align*}
\rti{\star k}&=\qti{m_3k}\ne 1, & \rti{\star m_1}&=\qti{im_1}q_{m_3m_1}=\qti{im_3}^{-1}\ne 1, & \rti{\star m_2}&=\qti{im_2}\ne 1
\end{align*}
and $\rti{\star \ell}=1$ for any other $\ell$.
The set of neighbours of $j$ is either $\{n_1,n_2,n_3\}$ if $k\ne j$, or $\{n_1,n_2,n_3,im_3\}-\{m_3\}$ if $k=j$. By inductive hypothesis $im_3$ and $j$ must be connected, that is $\qti{ij}\qti{m_3j}\ne 1$, which forces $\qti{m_3j}\ne 1$ and so $j=k$: that is, $m_3$ is a neighbour of $j$, say $m_3=n_3$. Also, by inductive hypothesis $im_3$ and $j$ share a neighbour: as the neighbours of $im_3$ and $j$ are $\{m_1, m_2, j\}$ and $\{n_1,n_2,im_3\}$ respectively, we get a contradiction. Thus \ref{item:rk6-2-vertices-of-degree-3-i} holds and we may fix $m_3=j$, $n_3=i$.

\smallbreak

Suppose now that $m_s\ne n_t$ for all $s,t\in\I_2$. Collapsing $i$ and $j$ we obtain a matrix $\bq'$ such that the vertex $ij$ has four neighbours: $m_1$, $m_2$, $n_1$, $n_2$. This is a contradiction with Proposition \ref{prop:rk5-degree-four}. We may assume that $m_2=n_2$. By Lemma \ref{lem:rk6-triangle-extended-3-vertices} $m_2$ is not connected to any other vertex $\ell\ne i,j$, so \ref{item:rk6-2-vertices-of-degree-3-ii} holds.

\smallbreak

As $\theta \ge 6$, there exists $\ell \ne i,j,m_1,n_1,m_2$ connected with one of these five vertices. Now $\ell$ is connected either with $m_1$ or $n_1$ as $\ell$ is not connected with $i,j,m_2$: we assume that $\ell$ is connected with $m_1$. Moreover the set of neighbours of $m_1$ is $\{i,\ell\}$ since otherwise we obtain a vertex of degree $\ge 4$ by collapsing $i$ and $m_1$. Collapsing $\ell$ and $m_1$ we get a matrix $\bq'$ where the sets of neighbours of any vertex different from $m_1$ and $\ell$ does not change and the set of neighbours of $m_1\ell$ is the union of the sets of neighbours of $m_1$ and $\ell$ for $\bq$ removing themselves: by inductive hypothesis $i$ and $j$ are the unique vertices of degree $\ge 3$, so \ref{item:rk6-2-vertices-of-degree-3-iii} also holds.
\epf

\subsection{On parametric subdiagrams}\label{subsec:eval-param-rk6}

Once again, we look into the admissible shapes of the underlying graph. Up to permutations, we may assume that $\bq'\coloneqq \bq_{|\I_5}$ is of one of the shapes determined in \ref{subsec:eval-param-rk5}. Let us fix $j\in\I_5$ be such that $\qti{j6}\neq 1$.

If $\bq'$ is of shape 
$\xymatrix@C-18pt@R-18pt{& & & {\circ} \ar @{-}[ld] \ar @{-}[rd]   &  &  & \\
{\circ} \ar @{-}[rr] & &{\circ} \ar @{-}[rr] & & {\circ} \ar @{-}[rr] & & {\circ}  }$, then $j\neq 5$ by Lemma \ref{lem:rk6-triangle-extended-3-vertices}. As well, $j\neq 2,3$ by Proposition \ref{prop:rk5-degree-four}. Thus $j$ can be only 1 or 4, and not both simultaneously to avoid cycles. This is the last figure in \eqref{eqn:shapes6}.

If in turn $\bq'$ is of shape  
$\xymatrix@C-18pt@R-18pt{ & {\circ} \ar @{-}[ld] \ar @{-}[rd]   &  &  & & & \\
{\circ} \ar @{-}[rr] & & {\circ} \ar @{-}[rr] & & {\circ} \ar @{-}[rr] & & {\circ}  }$ and $j=5$, then we recover the figure above as 6 cannot be joined with any other vertex, to avoid cycles. The same holds for $j=1$. Also, $j\neq 2$ since we cannot have degree 4. The same holds for $j=3$ by Lemma \ref{lem:rk6-2-vertices-of-degree-3}. Option $j=4$ gives the fourth figure in \eqref{eqn:shapes6} below.

Now, assume $\bq'$ is of shape  
$\xymatrix@C-18pt@R-18pt{ &  & &  &{\circ} \ar @{-}[d]    & & \\
{\circ} \ar @{-}[rr] & & {\circ} \ar @{-}[rr] & & {\circ} \ar @{-}[rr] & & {\circ}  }$. Notice that $j\neq 3$ to avoid a vertex of degree 4; also $j\neq 2$ by Lemma \ref{lem:rk6-2-vertices-of-degree-3}. If $j=4$, then $\qti{i6}= 1$ for all $i\neq 4$ to discard cycles and get the third shape in \eqref{eqn:shapes6}; case $j=5$ is symmetric. If $j=1$, then we get the second figure.

Finally, if  $\bq'$ is the line
$\xymatrix@C-18pt@R-18pt{{\circ} \ar @{-}[rr] & & {\circ} \ar @{-}[rr] & & {\circ} \ar @{-}[rr] & & {\circ} \ar @{-}[rr] & & {\circ} }$, then there are two cases. On the one hand, if there is $k\neq j$ such that $\qti{k6}\neq 1$, then $j$ and $k$ are connected: $\qti{jk}\neq 1$, to avoid $n$-cycles, $n\geq 4$. We get one of the figures already analyzed. On the other, $j$ is the single vertex with $\qti{j6}\neq 1$ and we obtain one of three shapes without triangles in \eqref{eqn:shapes6}.

We have thus obtained the following five possible shapes:
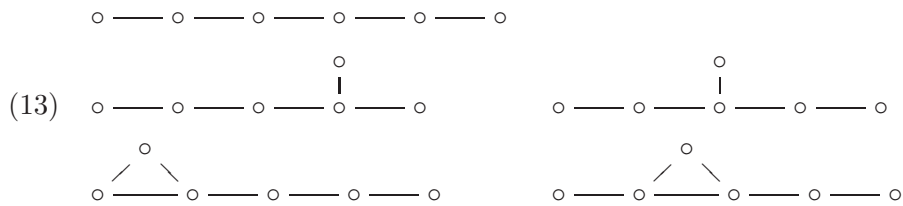
\begin{align}\label{eqn:shapes6}
\begin{aligned}
&  \xymatrix@C-18pt@R-18pt{
{\circ} \ar @{-}[rr] & & {\circ} \ar @{-}[rr] & & {\circ}\ar @{-}[rr] & & {\circ}  \ar @{-}[rr] & & {\circ} \ar @{-}[rr] & & {\circ}}  \\
& \xymatrix@C-18pt@R-18pt{  & & & & & &{\circ}  \ar @{-}[d]   & & \\
{\circ} \ar @{-}[rr] & & {\circ} \ar @{-}[rr] & & {\circ}\ar @{-}[rr] & & {\circ}  \ar @{-}[rr] & & {\circ} } &
& \xymatrix@C-18pt@R-18pt{ &  & &  & {\circ}  \ar @{-}[d] &  &  &  & \\
{\circ} \ar @{-}[rr] & &{\circ} \ar @{-}[rr] & & {\circ} \ar @{-}[rr] & & {\circ}\ar @{-}[rr] & & {\circ}   } \\
& \xymatrix@C-18pt@R-18pt{ & {\circ} \ar @{-}[ld] \ar @{-}[rd] & &  &  &  &  &  & \\
{\circ} \ar @{-}[rr] & & {\circ} \ar @{-}[rr] & & {\circ}\ar @{-}[rr] & & {\circ}  \ar @{-}[rr] & & {\circ}  } &
& \xymatrix@C-18pt@R-18pt{& & & {\circ} \ar @{-}[ld] \ar @{-}[rd]   &  &  & &  & \\
{\circ} \ar @{-}[rr] & &{\circ} \ar @{-}[rr] & & {\circ} \ar @{-}[rr] & & {\circ} \ar @{-}[rr] & & {\circ}  }
\end{aligned}
\end{align}

\begin{remark}\label{rem:collapsing-parametric}
Let $\bq$ a connected parametric matrix of rank $\theta\ge 5$ with parameter $q$ and $i$, $j$ two connected vertices of degree $\ge 2$. Collapsing $\alpha_i+\alpha_j$ we obtain a connected parametric diagram of rank $\theta-1$ with parameter $q^{\pm 1}$. The proof follows case-by-case by looking at \cite[Rows 1--4, 7--10]{H-full}.
\end{remark}

\begin{lemma}\label{lemm:reduction-to-Gf-rank6}
Let $\bq=(q_{ij})_{i,j\in\I_6}$ be such that a connected subdiagram of rank five is parametric evaluated in $q\in\Bbbk^{\times}$. 
Then either $\bq$ is in $\hlist$, $\gdim \B_{\bq}=\infty$, or else $q\in\Gf$.
\end{lemma}
\pf
We follow the same strategy as for the proof of Lemma \ref{lemm:reduction-to-Gf-rank5}. Assume that $\gdim \B_{\bq}<\infty$: By Theorem \ref{thm:rank5-conj-true}, all connected subdiagrams of rank five and those obtained by applying the criteria belong to $\hlist$. Also, the shape of the diagram is one of the five shapes in \eqref{eqn:shapes6}.

\begin{itemize}[leftmargin=*]\renewcommand{\labelitemi}{$\circ$}
\item If the shape of $\bq$ is the fifth in \eqref{eqn:shapes6}, then we collapse $3$ and $4$ to obtain a diagram of rank 5 as in Case \ref{item:rank5-parametric-shapeA} of the proof of Lemma \ref{lemm:reduction-to-Gf-rank5}. This diagram is necessarily finite since it has \emph{a triangle in the middle}, so $q\in\Gf$.
\item If the shape of $\bq$ is the fourth in \eqref{eqn:shapes6}, then we again collapse $3$ and $4$ to obtain a diagram of rank 5 as in Case \ref{item:rank5-parametric-shapeB} with a connected parametric subdiagram of rank 4 with parameter $q$, see Remark \ref{rem:collapsing-parametric}. An analysis as in Case \ref{item:rank5-parametric-shapeB} gives that either $\bq$ is in Rows 9--10 of \cite[Table 4]{H-full} or $q\in\Gf$.
\item If the shape of $\bq$ is the third in \eqref{eqn:shapes6}, then by collapsing $3$ and $4$ we obtain a diagram of rank 5 as in Case \ref{item:rank5-parametric-shapeC} with a connected parametric subdiagram of rank 4 with parameter $q$. Arguing as in Case \ref{item:rank5-parametric-shapeC}, either $\bq$ is in Row 16 of \cite[Table 4]{H-full} or $q\in\Gf$.
\item If the shape of $\bq$ is the second in \eqref{eqn:shapes6}, then we obtain a diagram of rank 5 as in Case \ref{item:rank5-parametric-shapeC} with a connected parametric subdiagram of rank 4 with parameter $q$ by collapsing $3$ and $4$. An analysis as in Case \ref{item:rank5-parametric-shapeC} gives that either $\bq$ is in Rows 8 or 10 of \cite[Table 4]{H-full}, or $q\in\Gf$.
\item If the shape of $\bq$ is the first one of \eqref{eqn:shapes6}, then by collapsing $3$ and $4$ we obtain a line of rank 5 with a connected parametric subdiagram of rank 4 with parameter $q$. An argument as in Case \ref{item:rank5-parametric-shapeD} gives that either $\bq$ is in Rows 1, 2, 3, 4, 7, 9 or 10 of \cite[Table 4]{H-full}, or $q\in\Gf$.
\end{itemize}
Thus, in any case, either $\bq$ is in $\hlist$ or else $q\in\Gf$.
\epf

\subsection{Proof of the conjecture in rank 6}

\begin{theorem}\label{thm:rank-6}
A rank 6 matrix $\bq$ is in $\hlist$ if and only if $\gdim \B_{\bq}<\infty$.
\end{theorem}
\pf
Assume that $\gdim \B_{\bq}<\infty$. If one of the connected subdiagrams of rank four is parametric not evaluated in $\Gf$, then $\bq$ is in $\hlist$ by Lemma
\ref{lemm:reduction-to-Gf-rank6}. Otherwise all subdiagrams are either finite or else parametric evaluated in $\Gf$ and we apply \S \ref{gap:criterian} to check that $\bq$ is in $\hlist$.
\epf

\section{Rank 7}\label{sec:rank7}

In this section we focus on diagrams of rank 7. 

\subsection{Some forbidden shapes}\label{subsec:forbidden7}
As in previous section we start by discarding diagrams of some particular shapes.

\begin{lemma}\label{lem:rk7-forbidden}
Assume $\bq$ has any of the following two shapes:
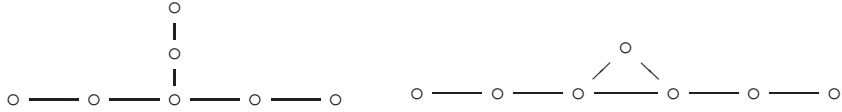
\begin{align*}
&  \xymatrix@C-18pt@R-18pt{ &  & &  & {\circ}  \ar @{-}[d]  & &   &  & \\
&  & &  & {\circ}  \ar @{-}[d]  & &   &  & \\
{\circ} \ar @{-}[rr] & &{\circ} \ar @{-}[rr] & & {\circ} \ar @{-}[rr] & & {\circ}\ar @{-}[rr] & & {\circ}   }  &
& \xymatrix@C-18pt@R-18pt{& & & &&    &  &  & &  & \\
& & & && {\circ} \ar @{-}[ld] \ar @{-}[rd]   &  &  & &  & \\
{\circ} \ar @{-}[rr] & &{\circ} \ar @{-}[rr] & &{\circ} \ar @{-}[rr] & & {\circ} \ar @{-}[rr] & & {\circ} \ar @{-}[rr] & & {\circ}  }
\end{align*}
Then $\gdim\B_{\bq}=\infty$.
\end{lemma}
\pf
We may assume, in any case, that all connected rank 6 subdiagrams of $\bq$ belong to $\hlist$, by Theorem \ref{thm:rank-6}.  
In the left hand case, these three subdiagrams are of shape $\xymatrix@C-18pt@R-18pt{ &  & &  & {\circ}  \ar @{-}[d]  & &   &  & \\
{\circ} \ar @{-}[rr] & &{\circ} \ar @{-}[rr] & & {\circ} \ar @{-}[rr] & & {\circ}\ar @{-}[rr] & & {\circ}   } $.  If at least two of them are of Cartan type $E_6$, then so is the third and $\gdim \B_{\bq}=\infty$ by Theorem \ref{thm:Cartan}. Otherwise, the diagram of $\bq$ is constructed by pasting two finite diagrams. We discard these diagrams using \gap, see \ref{gap:forbidden7}.

In the right hand, the diagram is obtained by pasting two diagrams with a triangle in the middle, which are necessarily finite. Once again, we remove them in \ref{gap:forbidden7} using \gap.
\epf

\subsection{On parametric subdiagrams}\label{subsec:eval-param-rk7}

We investigate the possible shapes of the underlying graph. We assume that $\bq'\coloneqq \bq_{|\I_6}$ is of one of the shapes in \ref{subsec:eval-param-rk6}. 
We fix $j\in\I_6$ with $\qti{j7}\neq 1$.

If the shape of $\bq'$ corresponds to the first diagram in \eqref{eqn:shapes6}, then $7$ can be connected either to two neighboring vertices, or just with $j$. In any case, we obtain one of the diagrams depicted below in \eqref{eqn:shapes7} using Lemma \ref{lem:rk7-forbidden}.

If $\bq'$ has the shape of the second diagram in \eqref{eqn:shapes6}, then $j\neq 4$ since the degree is bounded by 3. If $j=5$ (or, symmetrically, $j=6$), then $\qti{i7}=1$ for all $i\neq j$, and we get the third diagram in \eqref{eqn:shapes7}. Observe that $j\neq 2,3$ by Lemma \ref{lem:rk6-2-vertices-of-degree-3}. If $j=1$, then we get the second diagram in \eqref{eqn:shapes7}.

Now, if $\bq'$ is the third diagram in \eqref{eqn:shapes6}, then $j\neq 3$ as above. As well, $j=6$ is forbidden since either we get a long cycle if $\qti{i7}\neq 1$ for some $i\neq j$ or we get the figure in the left in Lemma \ref{lem:rk7-forbidden}. We discard cases $j=2,4$ by Lemma \ref{lem:rk6-2-vertices-of-degree-3} as before. If $j=1$, symmetrically $j=5$, then $\qti{i7}=1$ for $i\neq j$ and we arrive to the third diagram in  \eqref{eqn:shapes7}.

Assume $\bq'$ is the fourth diagram in \eqref{eqn:shapes6}. First $j\neq 2$ by the degree argument as above. If $j=1$ or symmetrically $j=6$, then the unique admissible shape is the fifth in  \eqref{eqn:shapes7}, again taking into account the restriction in the cycles. Once more, $j\neq 3,4$, by Lemma \ref{lem:rk6-2-vertices-of-degree-3}. If $j=5$, we get the fourth diagram in  \eqref{eqn:shapes7}.

If, finally, $\bq'$ is the fifth diagram in \eqref{eqn:shapes6}, then $j\neq 2,3$ because of the restrictions in the degree. As well, $j\neq 6$ since this figure is forbidden by Lemma \ref{lem:rk6-triangle-extended-3-vertices}. Similarly, $j\neq 1$, as this implies that $q_{i7}\neq 1$ for every $i\neq 1$ and we get the forbidden shape in the right in Lemma \ref{lem:rk7-forbidden}, and $j\neq 4$ by Lemma \ref{lem:rk6-2-vertices-of-degree-3}. Hence $j=5$ and we get the fifth diagram in  \eqref{eqn:shapes7}.

We have obtained five possible shapes:
\begin{align}\label{eqn:shapes7}
\begin{aligned}
& \xymatrix@C-22pt@R-18pt{{\circ} \ar @{-}[rr]  & &
{\circ} \ar @{-}[rr] & & {\circ} \ar @{-}[rr] & & {\circ}\ar @{-}[rr] & & {\circ}  \ar @{-}[rr] & & {\circ} \ar @{-}[rr] & & {\circ}}  \\
& \xymatrix@C-20pt@R-18pt{  && & & & & & &{\circ}  \ar @{-}[d]   & & \\
{\circ} \ar @{-}[rr]  & & {\circ} \ar @{-}[rr] & & {\circ} \ar @{-}[rr] & & {\circ}\ar @{-}[rr] & & {\circ}  \ar @{-}[rr] & & {\circ} } &
& \xymatrix@C-20pt@R-18pt{ & & &  & &  & {\circ}  \ar @{-}[d] &  &  &  & \\
{\circ} \ar @{-}[rr]  & & {\circ} \ar @{-}[rr] & &{\circ} \ar @{-}[rr] & & {\circ} \ar @{-}[rr] & & {\circ}\ar @{-}[rr] & & {\circ}   } \\
& \xymatrix@C-20pt@R-18pt{ & {\circ} \ar @{-}[ld] \ar @{-}[rd] & &  &  &  &  &  & && \\
{\circ} \ar @{-}[rr] & & {\circ} \ar @{-}[rr] & & {\circ}\ar @{-}[rr] & & {\circ}  \ar @{-}[rr] & & {\circ} \ar @{-}[rr] & & {\circ} } &
& \xymatrix@C-20pt@R-18pt{& & & {\circ} \ar @{-}[ld] \ar @{-}[rd]   &  &  & &  & && \\
{\circ} \ar @{-}[rr] & &{\circ} \ar @{-}[rr] & & {\circ} \ar @{-}[rr] & & {\circ} \ar @{-}[rr] & & {\circ}\ar @{-}[rr] & & {\circ}  }
\end{aligned}
\end{align}

\begin{lemma}\label{lemm:reduction-to-Gf-rank7}
Let $\bq=(q_{ij})_{i,j\in\I_7}$ be such that a connected subdiagram of rank six is parametric evaluated in $q\in\Bbbk^{\times}$.
Then either $\bq$ is in $\hlist$, $\gdim \B_{\bq}=\infty$, or else $q\in\Gf$.
\end{lemma}
\pf
Analogous to Lemma \ref{lemm:reduction-to-Gf-rank6}, as we have analogous shapes and the proof again follows by collapsing $3$ and $4$ as in loc.~cit.
\epf

\subsection{Proof of the conjecture in rank 7}

\begin{theorem}\label{thm:rank-7}
A rank 7 matrix $\bq$ is in $\hlist$ if and only if $\gdim \B_{\bq}<\infty$.
\end{theorem}
\pf
Assume that $\gdim \B_{\bq}<\infty$. If one of the connected subdiagrams of rank six is parametric not evaluated in $\Gf$, then $\bq$ is in $\hlist$ by Lemma
\ref{lemm:reduction-to-Gf-rank7}. Otherwise all subdiagrams are either finite or else parametric evaluated in $\Gf$ and we apply \S \ref{gap:criterian} to check that $\bq$ is in $\hlist$.
\epf

\section{Rank greater than 7}\label{sec:rank-ge-8}

Finally we deal with diagrams of big rank. Due to the rigidity on the admissible shapes we can now work inductively on the rank.

The implicit numeration of any generalized Dynkin diagram is from the left to the
right and from bottom to top; otherwise, the numeration appears below the vertices.
Recall that we say that we {\it collapse} vertices $i$ and $j$ when we apply Criterium \ref{crit-A} for $(i,j;1)$.

\subsection{Large diagrams with a triangle}

We start by considering the case in which the diagram contains a triangle. We will show that the only admissible shape is a (long) tadpole.

\begin{remark}\label{rem:list-ge7-mjk-ge-2}
Let $\bq$ be a matrix in $\hlist$ such that $\theta\ge 7$ and there are $j\ne k\in\I_{\theta}$ such that $m_{jk}\ge 2$. Then, by inspection, $m_{jk}=2$ and $\bq$ belongs to \cite[Table 4, Rows 3--7 \& 9--10]{H-full}. Moreover, the degree of every vertex is $1$ or $2$ (that is, the diagram is a line) and either $j$ or $k$ has no other neighbours; in other words, exactly one of the following holds:
\begin{itemize}
\item $\qti{ji}=1$ for all $i\ne j,k$ when $\bq$ is in Rows 3--6,
\item $\qti{ki}=1$ for all $i\ne j,k$ when $\bq$ is in Rows 7, 9 or 10.
\end{itemize}
\end{remark}

\begin{remark}\label{rem:list-ge7-vertex-degree-3}
Let $\bq$ be a braiding matrix in $\hlist$ such that $\theta\ge 7$ and there exist a vertex $i$ of degree 3. Then $\bq$ is in \cite[Table 4, Rows 8--10, 20--22]{H-full}.
\end{remark}

\begin{lemma}\label{lem:rank-ge-8-line&triang}
Fix $\theta\ge 8$ and assume that the diagram of $\bq$ contains a triangle. 
Then either $\gdim \B_{\bq}=\infty$ or else $\bq$ is in $\hlist$ (more precisely, the diagram belongs to \cite[Table 4, rows 9-10]{H-full}).
\end{lemma}

\pf
The proof is recursive on $\theta$. Suppose that the statement holds for $\theta\ge7$ and we want to prove it for $\theta+1\ge 8$. Notice that the case $\theta=7$ was proved in Theorem \ref{thm:rank-7}.
We assume that $\gdim \B_{\bq}<\infty$. This implies that $\bq_{|\I_{2,\theta+1}}$ is in $\hlist$.
Moreover, up to permutation of $\I_{2,\theta+1}$, we may assume that the triangle is contained in $\I_{2,\theta+1}$ and the underlying diagram is
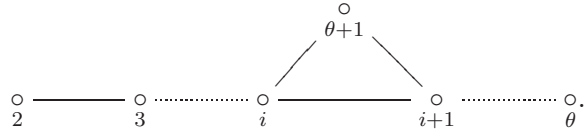
\begin{align*}
\xymatrix@C-10pt@R-10pt{ & & & & & \underset{\theta+1}{\circ} \ar @{-}[dl] \ar @{-}[rd] & & & \\
\underset{2}{\circ} \ar @{-}[rr] & & \underset{3}{\circ} \ar @{.}[rr]  & & \underset{i}{\circ} \ar @{-}[rr] & & \underset{i+1}{\circ} \ar @{.}[rr] & & \underset{\theta}{\circ}.}
\end{align*}

\begin{step}
Up to permutation, the underlying diagram of $\bq$ is
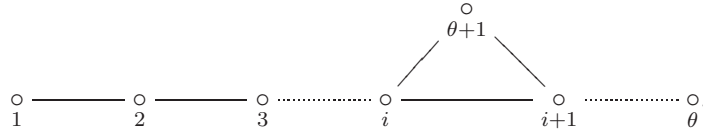
\begin{align}\label{eq:extended-triangle}
\xymatrix@C-10pt@R-10pt{ & & & & & & & \underset{\theta+1}{\circ} \ar @{-}[dl] \ar @{-}[rd] & & & \\
\underset{1}{\circ} \ar @{-}[rr] & & \underset{2}{\circ} \ar @{-}[rr] & & \underset{3}{\circ} \ar @{.}[rr]  & & \underset{i}{\circ} \ar @{-}[rr] & & \underset{i+1}{\circ} \ar @{.}[rr] & & \underset{\theta}{\circ}.}
\end{align}
\end{step}
\pf
As the diagram is connected, there exists $j\in\I_{2,\theta+1}$ such that $\qti{1j}\ne 1$.

If $i=\theta-1$, then $j\notin\I_{3,\theta-1}$ by Lemma \ref{lem:rk6-2-vertices-of-degree-3}, hence $j\in\{2,\theta,\theta+1\}$. 
This Lemma also implies that $1$ is not connected with any other vertex $\ne j$ and the diagram is as in \eqref{eq:extended-triangle}.

If $i<\theta-1$, then $j\in\{2,\theta\}$ by Lemma \ref{lem:rk6-2-vertices-of-degree-3} and $1$ is not connected with any other vertex $\ne j$. Thus, the diagram is as in \eqref{eq:extended-triangle}.
\epf

\begin{step}\label{Step2}
$i=\theta-1$, that is, the diagram is a long tadpole.
\end{step}
\pf
By symmetry we may assume that $3\le i\le \theta-1$. Suppose that $i<\theta-1$. 

Assume that $\theta+1=8$. Then either $\bq_{|\I_{2,8}}$  or $\bq_{|\I_{6}\cup\{8\}}$ is the second diagram in \cite[Table 4, row 21]{H-full} and the diagram of $\bq$ is one of the following, for some $q,r\in\Bbbk-\{1\}$ and $\zeta\in\G_3'$:
\begin{align*}
\xymatrix@C-10pt@R-10pt{ & & & & & & & & & \overset{-1}{\circ} \ar @{-}[dl]_{\zeta} \ar @{-}[rd]^{\zeta} & & & \\
\overset{q}{\circ} \ar @{-}[rr]^{r} & &\overset{\zeta}{\circ} \ar @{-}[rr]^{\overline{\zeta}} & & \overset{\zeta}{\circ} \ar @{-}[rr]^{\overline{\zeta}} & & \overset{\zeta}{\circ} \ar @{-}[rr]^{\overline{\zeta}}  & & \overset{-1}{\circ} \ar @{-}[rr]^{\zeta} & & \overset{-1}{\circ} \ar @{-}[rr]^{\overline{\zeta}} & & \overset{\zeta}{\circ},}
\\
\xymatrix@C-10pt@R-10pt{ & & & & & & & \overset{-1}{\circ} \ar @{-}[dl]_{\zeta} \ar @{-}[rd]^{\zeta} & & & & & \\
\overset{\zeta}{\circ} \ar @{-}[rr]^{\overline{\zeta}} & & \overset{\zeta}{\circ} \ar @{-}[rr]^{\overline{\zeta}} & & \overset{\zeta}{\circ} \ar @{-}[rr]^{\overline{\zeta}}  & & \overset{-1}{\circ} \ar @{-}[rr]^{\zeta} & & \overset{-1}{\circ} \ar @{-}[rr]^{\overline{\zeta}} & & \overset{\zeta}{\circ}
\ar @{-}[rr]^{r} & &\overset{q}{\circ}.}
\end{align*}
For the second, $\bq_{|\I_{2,8}}$ is not in $\hlist$ by Lemma \ref{lem:rk7-forbidden}. Thus we assume that the diagram of $\bq$ is the first. Applying Criterium \ref{crit-A} for $(2,3;1)$ we obtain a diagram with a triangle in the middle, which must be the second diagram in \cite[Table 4, row 21]{H-full}. Thus $q=\zeta$, $r=\overline{\zeta}$. We apply Criterium \ref{crit-C} and obtain the diagram
\begin{align*}
\xymatrix@C-10pt@R-10pt{ \overset{\zeta}{\underset{1}{\circ}} \ar @{-}[rr]^{\overline{\zeta}} & &\overset{\zeta}{\underset{23}{\circ}} \ar @{-}[rr]^{\overline{\zeta}} & & \overset{-1}{\underset{45}{\circ}} \ar @{-}[rr]^{\zeta} & & \overset{-1}{\underset{67}{\circ}} \ar @{-}[rr]^{\zeta}  & & \overset{-1}{\underset{8}{\circ}} \ar @{-}[rr]^{\overline{\zeta}} & & \overset{\zeta}{\underset{56}{\circ}} \ar @{-}[rr]^{\overline{\zeta}} & & \overset{\zeta}{\underset{34}{\circ}},}
\end{align*}
which is not in $\hlist$. Hence we get a contradiction, so $i=6$.

\smallbreak

Now, if $\theta+1>8$, then the subdiagram with vertices in $\I_{2,\theta+1}$ does not belong to $\hlist$ and we get a contradiction. 

In any case we get that $i=\theta-1$ and the diagram is a tadpole.
\epf

By inductive hypothesis and Step \ref{Step2}, the diagram of $\bq_{|\I_{2,\theta+1}}$ is one of the parametric tadpoles in 
\cite[Table 4, rows 9-10]{H-full}. Applying
Criterium \ref{crit-A} for $(\theta-2,\theta-1;1)$ we obtain the following diagram:
\begin{align*}
&\xymatrix@C-10pt@R-10pt{ & & & & & & & \overset{-1}{\circ} \ar @{-}[dl]_{q^{-1}} \ar @{-}[rd]^{q^2} & \\
\underset{q_{11}}{\circ} \ar @{-}[rr]^{\qti{12}} & & \underset{q_{22}}{\circ} \ar @{-}[rr]^{\qti{23}} & & \underset{q_{33}}{\circ} \ar @{.}[rr]  & & \underset{s}{\circ} \ar @{-}[rr]^{q^{-1}} & & \underset{-1}{\circ},}
\end{align*}
where $q\in\Bbbk-\{\pm 1\}$ and $s:=q_{\theta-2\,\theta-2}\qti{\theta-2\,\theta-1}q_{\theta-1\,\theta-1}\in\{-1,q\}$. As this diagram is in $\hlist$, 
then it belongs to the same rows, which forces 
\begin{align*}
\qti{12}&=q_{22}^{-2}\qti{23}^{-1}, \qquad q_{11}\in \{-1,\qti{12}^{-1}\},
\end{align*}
and therefore $\bq$ also belongs to \cite[Table 4, rows 9-10]{H-full}.
\epf

\subsection{Proof of the conjecture in rank $>7$}

Finally we answer Conjecture \ref{conjecture} for diagrams of big rank.

\begin{theorem}\label{thm:rank-ge-8}
Let $\theta\ge 8$. A matrix $\bq$ of rank $\theta$ is in $\hlist$ if and only if $\gdim \B_{\bq}<\infty$.
\end{theorem}
\pf
We prove by induction on $\theta\ge 7$ that if a matrix $\bq$ is such that $\gdim \B_{\bq}<\infty$, then $\bq$ is in $\hlist$: recall that the opposite implication always holds as explained in the Introduction. The case $\theta=7$ is Theorem \ref{thm:rank-7}.

Let $\theta\ge 8$. By Theorem \ref{thm:Cartan} we may assume that $\bq$ is not of Cartan type.
By Lemma \ref{lem:rank-ge-8-line&triang} we may assume that $\bq$ does not contain a triangle. Up to permutation of the vertices we may assume that $\bq_{|\I_{\theta-1}}$ is connected: By Theorem \ref{thm:rank-7} or inductive hypothesis, this subdiagram belongs to $\hlist$. 
Thus we study each possibility for this subdiagram.

\begin{description}[leftmargin=5pt]
\item[Rows 9-10] $\bq$ is Weyl equivalent to a braiding matrix $\bq'$ whose diagram contains a triangle, and $\gdim \B_{\bq'}=\gdim \B_{\bq}<\infty$. Now, by Lemma \ref{lem:rank-ge-8-line&triang} $\bq'$ belongs to $\hlist$, and so does $\bq$.
\end{description}

\smallbreak

For the remaining cases, namely \cite[Table 4, rows 1-8, 20, 22]{H-full}, we may further assume that the diagram does not contain a triangle and $\theta$ is joined to exactly one vertex $i\in\I_{\theta-1}$ since there are no cycles of size $\ge 4$ by Proposition \ref{prop:n-cycles}. Thus we need to determine the values of $q_{\theta\theta}$ and $\qti{i\theta}$.

\smallbreak

\begin{description}[leftmargin=5pt]
\item[Row 20] Here $\theta=8$, the diagram of $\bq_{|\I_{7}}$ is
\begin{align*}
&\xymatrix@C-10pt@R-10pt{ & & & & & & \underset{q}{\circ} \ar @{-}[d]_{q^{-1}} & &  & & \\
\underset{q}{\circ} \ar @{-}[rr]^{q^{-1}} & & \underset{q}{\circ} \ar @{-}[rr]^{q^{-1}} & & \underset{q}{\circ} \ar @{-}[rr]^{q^{-1}} & & \underset{q}{\circ} \ar @{-}[rr]^{q^{-1}} & & \underset{q}{\circ} \ar @{-}[rr]^{q^{-1}}  & & \underset{q}{\circ},}
\end{align*}
and we look for the possible values of $q_{88}$, $\qti{i8}$. If $i\ne 1$, then $\bq_{|\I_{2,8}}$ is connected so it is an extension of a diagram of Cartan type $E_6$: either $\bq_{|\I_{2,8}}$ is of Cartan type $E_7$, or else belongs to \cite[Table 4, row 21]{H-full}, where $i=6$ in both cases. The first case is of affine Cartan type so $\gdim\B_{\bq}=\infty$ by Theorem \ref{thm:Cartan}, a contradiction. For the second case $\bq$ is Weyl equivalent to a braiding matrix $\bq'$ whose diagram contains a triangle and we get a contradiction by Lemma \ref{lem:rank-ge-8-line&triang}. Thus $i=1$ and $\bq_{|\I_{8}-\{6\}}$ is an extension of a diagram of Cartan type $D_6$ whose underlying graph is of type $D_7$. By inspection, $\qti{18}=q^{-1}$ and $q_{88}\in\{q,-1\}$. If $q_{88}=q$, then $\bq$ is of Cartan type $E_8$, so it is in $\hlist$. If $q_{88}=-1$, then $\bq$ is Weyl equivalent to a braiding matrix $\bq'$ whose diagram contains a triangle and we apply Lemma \ref{lem:rank-ge-8-line&triang}.

\item[Row 22] Here $\theta=9$ and the diagram of $\bq_{|\I_{8}}$ is of Cartan type $E_8$. The analysis is completely as in Row 20; moreover, in this case 
no diagram is allowed since any extension of Cartan type is necessarily not finite.

\item[Row 8] We recall that the diagram of $\bq_{|\I_{\theta-1}}$ is
\begin{align*}
&\xymatrix@C-10pt@R-10pt{ & & & & & & \underset{q}{\circ} \ar @{-}[d]_{q^{-1}} & & \\
\underset{q}{\circ} \ar @{-}[rr]^{q^{-1}} & & \underset{q}{\circ} \ar @{-}[rr]^{q^{-1}} & & \underset{q}{\circ} \ar @{.}[rr]  & & \underset{q}{\circ} \ar @{-}[rr]^{q^{-1}} & & \underset{q}{\circ},}
\end{align*}
for some $q\ne 1$. By Proposition \ref{prop:rk5-degree-four}, $i\ne \theta-3$, and by Lemma \ref{lem:rk6-2-vertices-of-degree-3}, $i\notin\I_{2,\theta-4}$. Thus $i\in\{1,\theta-2,\theta-1\}$: the cases $\theta-2$ and $\theta-1$ are analogous.

If $i=\theta-2$, then $\bq_{|\I_{2,\theta}}$ has a connected diagram with a vertex ($\theta-3$) of degree there and two \emph{ramifications} of length $\ge 2$. By Remark \ref{rem:list-ge7-vertex-degree-3}, $\bq_{|\I_{2,\theta}}$ is of Cartan type. Then $\bq$ is so, and Theorem \ref{thm:Cartan} applies.

If $i=1$, then we collapse $1$ and $2$ to get a diagram of rank $\theta-1$ where the vertex $\theta-3$ still has degree 3. By Remark \ref{rem:list-ge7-vertex-degree-3} this diagram belongs to Rows 8 or 10, so $q_{\theta\theta}\in\{-1,q\}$ and $\qti{1\theta}=q^{-1}$. Thus $\bq$ is in Rows 8 or 10.

\item[Row 7] In this case, the diagram of $\bq_{|\I_{\theta-1}}$ is
\begin{align*}
&\xymatrix@C-10pt@R-10pt{ \underset{q}{\circ} \ar @{-}[rr]^{q^{-1}} & & \underset{q}{\circ} \ar @{-}[rr]^{q^{-1}} & & \underset{q}{\circ} \ar @{.}[rr]  & & \underset{q}{\circ} \ar @{-}[rr]^{q^{-2}} & & \underset{q^2}{\circ},}
\end{align*}
for some $q\in\Bbbk-\{0,\pm 1\}$. Suppose that $i\ge 3$. Then $\bq_{|\I_{2,\theta}}$ (which belongs to $\hlist$ by inductive hypothesis) has either a diagram with a vertex of degree 3 or else $i=\theta-1$. But this is a contradiction with Remark \ref{rem:list-ge7-mjk-ge-2} since $m_{\theta-2 \,\theta-1}=2$. Hence $i\le 2$.

If $i=2$, then we collapse $2$ and $3$ and obtain a diagram where the new vertex $23$ has degree 3 and we still have $m_{\theta-2 \,\theta-1}=2$. This is a contradiction with Remark \ref{rem:list-ge7-mjk-ge-2}. Hence $i=1$: we collapse $1$ and $2$ to get a diagram of rank $\theta-1$ where $m_{\theta-2 \,\theta-1}=2$. By Remark \ref{rem:list-ge7-mjk-ge-2} this diagram belongs to Rows 7, 9 or 10, so $q_{\theta\theta}\in\{-1,q\}$ and $\qti{1\theta}=q^{-1}$. Thus $\bq$ belongs to the same rows.

\item[Rows 5 \& 6] In this case the diagram of $\bq_{|\I_{\theta-1}}$ is
\begin{align*}
&\xymatrix@C-10pt@R-10pt{ \underset{q_{11}}{\circ} \ar @{-}[rr]^{\qti{12}} & & \underset{q_{22}}{\circ} \ar @{-}[rr]^{\qti{23}} & & \underset{q_{33}}{\circ} \ar @{.}[rr]  & & \underset{q_{\theta-2\, \theta-2}}{\circ} \ar @{-}[rr]^{-\zeta} & & \underset{\zeta}{\circ},}
\end{align*}
where $\zeta\in\G_3'$, $q_{jj}\in\{-1,\zeta^{\pm 1}\}$, $q_{jj}\qti{j \, j\pm 1}=1$ if $q_{jj}\ne -1$ and $\qti{j\,j-1}=\qti{j\, j+1}^{-1}=q^{\pm 1}$ if $q_{jj}=-1$. As in the previous case, $i<3$ by Remark \ref{rem:list-ge7-mjk-ge-2} since $m_{\theta-1 \,\theta-2}=2$, the case $i=2$ gives a contradiction by collapsing $2$ and $3$, and for $i=1$ we collapse $1$ and $2$ to get a diagram of rank $\theta-1$ where $m_{\theta-1 \,\theta-2}=2$, $q_{\theta-1\,\theta-1}=\zeta$ and $\qti{\theta-2\, \theta-1}=-\zeta$. Hence this diagram belongs to Rows 5--6, and $\bq$ belongs to the same rows.

\item[Rows 3 \& 4] Here the diagram of $\bq_{|\I_{\theta-1}}$ is
\begin{align*}
&\xymatrix@C-10pt@R-10pt{ \underset{q_{11}}{\circ} \ar @{-}[rr]^{\qti{12}} & & \underset{q_{22}}{\circ} \ar @{-}[rr]^{\qti{23}} & & \underset{q_{33}}{\circ} \ar @{.}[rr]  & & \underset{q_{\theta-2\, \theta-2}}{\circ} \ar @{-}[rr]^{q^{-2}} & & \underset{q}{\circ},}
\end{align*}
where $q\notin\{0,\pm 1\}$, $q_{jj}\in\{-1,q^{\pm 2}\}$, $q_{jj}\qti{j \, j\pm 1}=1$ if $q_{jj}\ne -1$ and $\qti{j\,j-1}=\qti{j\, j+1}^{-1}=q^{\pm 1}$ if $q_{jj}=-1$. Again, $i<3$ by Remark \ref{rem:list-ge7-mjk-ge-2} since $m_{\theta-1 \,\theta-2}=2$, the case $i=2$ gives a contradiction by collapsing $2$ and $3$, and for $i=1$ we collapse $1$ and $2$ to get a diagram of rank $\theta-1$ where $m_{\theta-1 \,\theta-2}=2$. Hence this diagram belongs to Rows 3--4, and $\bq$ belongs to the same rows.

\item[Rows 1 \& 2] That is, the diagram of $\bq_{|\I_{\theta-1}}$ is
\begin{align*}
&\xymatrix@C-10pt@R-10pt{ \underset{q_{11}}{\circ} \ar @{-}[rr]^{\qti{12}} & & \underset{q_{22}}{\circ} \ar @{-}[rr]^{\qti{23}} & & \underset{q_{33}}{\circ} \ar @{.}[rr]  & & \underset{q_{\theta-2\, \theta-2}}{\circ} \ar @{-}[rr]^{\qti{\theta-2\, \theta-1}} & & \underset{q_{\theta-1\, \theta-1}}{\circ},}
\end{align*}
where $q_{jj}\in\{-1,q^{\pm 1}\}$, $q_{jj}\qti{j \, j\pm 1}=1$ if $q_{jj}\ne -1$ and $\qti{j\, j-1}=\qti{j\, j+1}^{-1}=q^{\pm 1}$ if $q_{jj}=-1$,
for $q=q_{11}^2\qti{12}\ne 1$. Recall that $i\in\I_{\theta-1}$ is the unique vertex such that $\qti{i\theta}\ne 1$. 

If $3\le i\le \theta-2$, then $\bq_{|\I_{i,i,\theta}}$ has a ramification, leading to a case already analysed above. The same happens for $i=\theta-2$ by looking at  $\bq_{|\I_{\theta}-\{\theta-1\}}$.

If $i=\theta-1$ and either $m_{\theta-1\, \theta}\ge 2$ or $m_{\theta\, \theta-1}\ge 2$, then $\bq_{|\I_{2,\theta}}$ belongs to Rows 3--7, 9 or 10 by Remark \ref{rem:list-ge7-mjk-ge-2}. Thus $\bq$ is in $\hlist$ either by Lemma \ref{lem:rank-ge-8-line&triang} or the analysis of the corresponding row above. The same holds if $i=1$ by looking at $\bq_{|\I_{\theta}-\{\theta-1\}}$.

Finally, assume that $i\in\{1,\theta-1\}$ and $m_{i\theta}=m_{\theta i}=1$: in this case $\bq$ belongs to Rows 1 \& 2 again.
\end{description}

\noindent Hence $\bq$ is in $\hlist$ for all possible cases.
\epf

\appendix

\section{The GAP lemmas}

Here we provide a description of our calculations in \gap. This is not a line-by-line analysis but rather we point out the main ideas and processes. 
We explain what we do in each file, including part of the code. The complete \gap \, files are stored in the authors' webpages for reference.

\subsection{Basic setting} 

\subsubsection{basic.g}\label{gap:basic}  

In this file we introduce the basic definitions we shall use along the other files, such as the set $\G_f$, see \eqref{eq:defn-G-f} and our function \textsl{Zet}:
\begin{lstlisting}
gap> Zet:=function(l)
> local d,f;
> f:=[];
> for d in l do
> if d in f then
> else
> Add(f,d);
> fi;
> od;
> return f;
> end;
function( l ) ... end
\end{lstlisting}
This replaces \gap-integrated function \texttt{Set}
and serves better for our purposes. For $\G_f$, we build our field \textsl{K}, along with the roots of unity as
\begin{lstlisting}
gap> q:=Indeterminate(Rationals,"q");;
gap> K:=Field(q);;
gap> F1:=One(K);;
gap> F:=n->E(n)*F1;;
gap> C:=function(n)
> local l,k;
> l:=[];
> for k in [1..n] do
> if F(n)^k<>F1 then
> Add(l,F(n)^k);
> fi;
> od;
> return l;
> end;
function( n ) ... end
\end{lstlisting}
Thus we define
\begin{lstlisting}
gap> Gf:=Zet(Concatenation(C(10),C(12),C(18)));;
\end{lstlisting}

\subsubsection{rank3.g}\label{gap:rank3} 
We introduce a function \texttt{Hcheckr3} to determine if a given rank 3 diagram is in $\hlist$. We list all triangles and all lines. Generic diagrams are evaluated in $\G_f$.

A line is represented by a vector
\begin{align*}
&\linethree{v_1}{a_1}{v_2}{a_2}{v_3} & &\longleftrightarrow & & d=[d[1],d[2],d[3],d[4],d[5]],
\end{align*}
where $v_1\coloneqq d[1], v_2\coloneqq d[3], v_3\coloneqq d[5]$ are the label of the vertices and $a_1\coloneqq d[2], a_2\coloneqq d[4]$ are the labels of edges connecting $v_1$ with $v_2$ and $v_2$ with $v_3$, respectively. 

A triangle, in turn, is a vector 
\begin{align*}
& \tri{v_1}{a_1}{v_2}{a_2}{v_3}{a_3} & &\longleftrightarrow & & d=[d[1],d[2],d[3],d[4],d[5],d[6]],
\end{align*}
where $a_3\coloneqq d[6]$ is the label of the edge connecting $v_3$ with $v_1$. When needed, we represent a line as a vector $[d[1],d[2],d[3],d[4],d[5],d[6]]$ with $d[6]=1$:
These ideas can be represented graphically as
\begin{align*}
& \linethree{v_1}{a_1}{v_2}{a_2}{v_3} & &\longmapsto &
& \xymatrix@C-4pt@R-18pt{ & \overset{v_3}{\underset{\ }{\circ}}  \ar @{-}[rd]^{a_2} & \\
\overset{v_1}{\underset{\ }{\circ}} \ar @{.}[ru]^{1} \ar @{-}[rr]_{a_1} & & \overset{v_2}{\underset{\ }{\circ}}}
\end{align*}

\subsection{Rank 4}

A diagram of rank 4 is represented as a 10-uple 
\[
d=[d[1],d[2],\dots,d[10]].
\] We think of $d[1]$, $d[3]$, $d[5]$ and $d[7]$ as the vertices of a square and the rest of entries are all possible edges on the full graph as in the diagram:
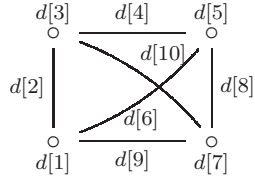
\begin{align*}
&\xymatrix@C40pt@R-5pt{ \overset{d[3]}{\circ} \ar @{-}[d]_{d[2]} \ar @{-}[r]^{d[4]} \ar@{-}@/^/[rd]^{d[10]}& 
\overset{d[5]}{\circ}\ar @{-}[d]^{d[8]} \\
\underset{d[1]}{\circ} \ar@{-}[r]_{d[9]} \ar@{-}@/_/[ru]_{d[6]\qquad\qquad} & \underset{d[7]}{\circ}}
\end{align*}
We represent the lines as the subset of vectors with $d[6]=d[9]=d[10]=1$, namely:
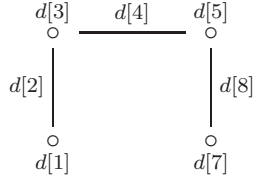
\begin{align*}
&\xymatrix@C40pt@R-5pt{ \overset{d[3]}{\circ} \ar @{-}[d]_{d[2]} \ar @{-}[r]^{d[4]} & 
\overset{d[5]}{\circ}\ar @{-}[d]^{d[8]} \\
\underset{d[1]}{\circ} & \underset{d[7]}{\circ}}
\end{align*}
Similarly, tadpoles are vectors with $d[9]=d[10]=1$ and tripods are those with $d[2]=d[9]=d[10]=1$.
Namely, tadpoles and tripods are considered, respectively as the subdiagrams
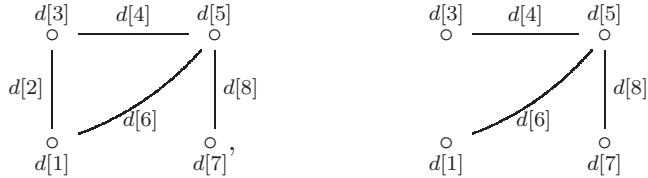
\begin{align*}
&\xymatrix@C40pt@R-5pt{ \overset{d[3]}{\circ} \ar @{-}[d]_{d[2]} \ar @{-}[r]^{d[4]} & 
\overset{d[5]}{\circ}\ar @{-}[d]^{d[8]} \\
\underset{d[1]}{\circ}  \ar@{-}@/_/[ru]_{d[6]\qquad\qquad} & \underset{d[7]}{\circ},}
&\xymatrix@C40pt@R-5pt{ \overset{d[3]}{\circ}  \ar @{-}[r]^{d[4]} & 
\overset{d[5]}{\circ}\ar @{-}[d]^{d[8]} \\
\underset{d[1]}{\circ}  \ar@{-}@/_/[ru]_{d[6]\qquad\qquad} & \underset{d[7]}{\circ}}
\end{align*}

\subsubsection{rank4t.g}\label{gap:rank4t} We list all the rank 4 diagrams in $\hlist$. 
We include permutations. 
Generic diagrams are evaluated in $\G_f$. 

All of these diagrams are vectors with 10 entries, as explained above. They are either lines (\texttt{\#li}), tripods (\texttt{\#tr}) ($D_4$-like diagrams)  or tadpoles (\texttt{\#ta}). Some examples of parametric diagrams are:
\begin{lstlisting}
gap> u:=F1;;
gap> d1q:=q->[q,q^-1,q,q^-1,q,F1,q,q^-1,F1,F1];;
gap> #li
gap> d2q:=function(q)
> if q<>-F1 then return
> [q^2,q^-2,q^2,q^-2,q^2,F1,q,q^-2,F1,F1];
> else return 0;
> fi;
> end;
function( q ) ... end
gap> #li
gap> d5q:=function(q)
> if q<>F1 then return
> [q,u,q,q^-1,q,q^-1,q,q^-1,u,u];
> else return 0;
> fi;
> end;
function( q ) ... end
gap> #tr
gap> d8dq:=function(q)
> if q^2<>F1 then return
> [-u,q^2,-u,q^-1,q,q^-1,q,q^-1,u,u];
> else return 0;
> fi;
> end;
function( q ) ... end
gap> #ta
\end{lstlisting}
We compile them all in separate lists. For tripods:
\begin{lstlisting}
gap> tr4q:=[d5q,d12fq,d13cq,d13dq];;
gap> tr4qq:=[];;
gap> for d in tr4q do
>       for n in Gf do
>               if d(n)<>0 then
>               Add(tr4qq,d(n));
>               fi;
>       od;
> od;
\end{lstlisting}
A finite tripod is 
\begin{lstlisting}
gap> z:=F(3);;
gap> d20i:=[z,u,z^2,z,-u,z^2,-u,z,u,u];;
gap> #tr
\end{lstlisting}
We consider the same diagram for each corresponding root of 1:
\begin{lstlisting}
gap> z:=F(3)^2;;
gap> d20ip:=[z,u,z^2,z,-u,z^2,-u,z,u,u];;
gap> #tr
\end{lstlisting}
Thus we obtain the complete lists of each type of diagram, with permutations.
\begin{lstlisting}
gap> tr4f:=[d18d,d18e,d18f,d20i,d20j,d18dp,d18ep,d18fp,
d20ip,d20jp];;
gap> Htr4s:=Zet(Concatenation(tr4qq,tr4f));;
gap> Htr4p:=[];;
gap> for d in Htr4s do
> Add(Htr4p,Permuted(d,(1,3)(4,6)));
> Add(Htr4p,Permuted(d,(1,3,7)(6,4,8)));
> Add(Htr4p,Permuted(d,(1,7,3)(6,8,4)));
> Add(Htr4p,Permuted(d,(3,7)(4,8)));
> Add(Htr4p,Permuted(d,(1,7)(6,8)));
> od;
\end{lstlisting}

\subsubsection{criteria4.g}\label{gap:criteria4}

In this file we deal with 10-uples representing the rank four diagrams introduced in the previous section. We settle the criteria from \S,  . As well we write down functions that reflect a diagram at each vertex, to apply the criteria to reflections.

We start loading the files with the basic data and $\hlist$ of rank three and four diagrams in $\hlist$:
\begin{lstlisting}
gap> Read("/home/.../basic.g");
gap> Read("/home/.../rank3.g");
gap> Read("/home/.../rank4t.g");
\end{lstlisting}
We shall make use of the following criteria. First we check if a diagram is of Cartan type:
\begin{lstlisting}
gap> ic;
function( d ) ... end
\end{lstlisting}
Next we build Criteria \ref{crit-A}, \ref{crit-D}, \ref{crit-D} for reflections and \ref{crit-E} (including reflections): 
\begin{lstlisting}
gap> cri2;
function( d ) ... end
gap> cri3b;
function( d ) ... end
gap> cri3;
function( d ) ... end
gap> cri4;
function( d ) ... end
\end{lstlisting}
As well, we check if a given diagram is a square, or a reflection of it. Next, if a diagram is a tadpole, or a reflection of it, or if the reflection at vertex 3 returns a tadpole or a square:
\begin{lstlisting}
gap> checksq;
function( d ) ... end
gap> checkrhosq;
function( d ) ... end
gap> checkta;
function( d ) ... end
gap> checkrhota;
function( d ) ... end
gap> checkrho3ta;
function( d ) ... end
\end{lstlisting}
The point is that once we discard squares we can discard a given diagram if a reflection of it becomes a square. Idem for tadpoles.

\medbreak
\label{page:squares}
Now we start our analysis. We build all squares pasting triangles $d$ and $e$ (six-tuples) with a matching side ($d[1]=e[1]$, $d[5]=e[5]$, $d[6]=e[6]$):
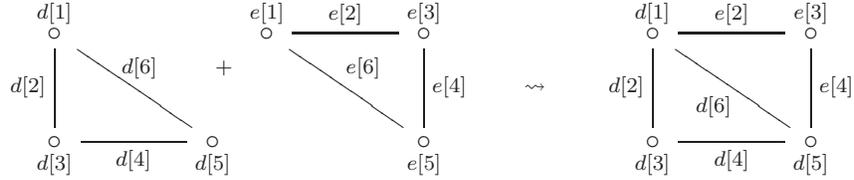
\begin{align*}
&\xymatrix@C40pt@R-5pt{ \overset{d[1]}{\circ} \ar @{-}[d]_{d[2]}  \ar@{-}[rd]^{d[6]\qquad +}& 
\\
\underset{d[3]}{\circ} \ar@{-}[r]_{d[4]}  & \underset{d[5]}{\circ}}
\xymatrix@C40pt@R-5pt{ \overset{e[1]}{\circ}  \ar @{-}[r]^{e[2]} \ar@{-}[rd]^{e[6]}& 
\overset{e[3]}{\circ}\ar @{-}[d]^{e[4]\qquad \rightsquigarrow} \\
 & \underset{e[5]}{\circ}}
&\xymatrix@C40pt@R-5pt{ \overset{d[1]}{\circ} \ar@{-}[rd]_{d[6]} \ar @{-}[d]_{d[2]} \ar @{-}[r]^{e[2]} & 
\overset{e[3]}{\circ}\ar @{-}[d]^{e[4]} \\
\underset{d[3]}{\circ} \ar@{-}[r]_{d[4]} & \underset{d[5]}{\circ}}
\end{align*}
and keep those such that the new lines generated are in $\hlist$, namely:
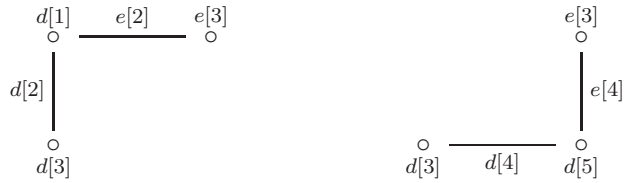
\begin{align*}
&\xymatrix@C40pt@R-5pt{ \overset{d[1]}{\circ} \ar @{-}[d]_{d[2]} \ar@{-}[r]^{e[2]} & \overset{e[3]}{\circ}
 \\
\underset{d[3]}{\circ} &}
&\xymatrix@C40pt@R-5pt{ & 
\overset{e[3]}{\circ}\ar @{-}[d]^{e[4]} \\
 \underset{d[3]}{\circ}  \ar @{-}[r]_{d[4]} & \underset{d[5]}{\circ}}
\end{align*}
These are the ``squares with a diagonal''. We also consider ``clean squares'' (no diagonal) and ``complete squares''. 
For the later, we paste, when possible, a square  with a diagonal as above (pasting of two triangles) with an extra triangle:
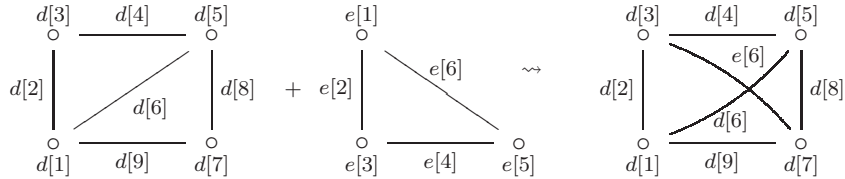
\begin{align*}
&\xymatrix@C40pt@R-5pt{ \overset{d[3]}{\circ} \ar @{-}[d]_{d[2]} \ar @{-}[r]^{d[4]} & 
\overset{d[5]}{\circ}\ar @{-}[d]^{d[8]\quad  +} \\
\underset{d[1]}{\circ} \ar@{-}[r]_{d[9]} \ar@{-}[ru]_{d[6]} & \underset{d[7]}{\circ}}
\xymatrix@C40pt@R-5pt{ \overset{e[1]}{\circ} \ar @{-}[d]_{e[2]}  \ar@{-}[rd]^{e[6]\qquad \rightsquigarrow}& 
\\
\underset{e[3]}{\circ} \ar@{-}[r]_{e[4]}  & \underset{e[5]}{\circ}}
&\xymatrix@C40pt@R-5pt{ \overset{d[3]}{\circ} \ar @{-}[d]_{d[2]} \ar @{-}[r]^{d[4]} \ar@{-}@/^/[rd]^{e[6]}& 
\overset{d[5]}{\circ}\ar @{-}[d]^{d[8]} \\
\underset{d[1]}{\circ} \ar@{-}[r]_{d[9]} \ar@{-}@/_/[ru]_{d[6]\qquad\qquad} & \underset{d[7]}{\circ}}
\end{align*}
 and check if the new triangle generated by this process is in $\hlist$:
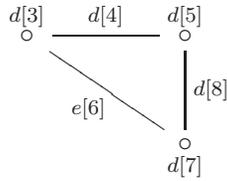
 \begin{align*}
& \xymatrix@C40pt@R-5pt{ \overset{d[3]}{\circ}  \ar @{-}[r]^{d[4]} \ar@{-}[rd]_{e[6]}& 
 \overset{d[5]}{\circ}\ar @{-}[d]^{d[8]} \\
 & \underset{d[7]}{\circ}}
 \end{align*}
The clean squares are generated by pasting two lines $d$ and $e$ with $d[1]=e[1]$ and $d[5]=e[5]$: \label{page:clean-squares}
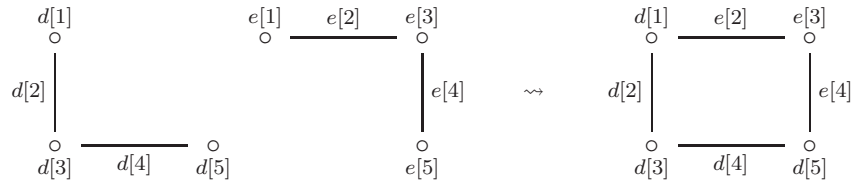
 \begin{align*}
&\xymatrix@C40pt@R-5pt{ \overset{d[1]}{\circ} \ar @{-}[d]_{d[2]}  & 
\\
\underset{d[3]}{\circ} \ar@{-}[r]_{d[4]}  & \underset{d[5]}{\circ}}
\xymatrix@C40pt@R-5pt{ \overset{e[1]}{\circ}  \ar @{-}[r]^{e[2]}& 
\overset{e[3]}{\circ}\ar @{-}[d]^{e[4]\qquad \rightsquigarrow} \\
& \underset{e[5]}{\circ}}
&\xymatrix@C40pt@R-5pt{ \overset{d[1]}{\circ} \ar @{-}[d]_{d[2]} \ar @{-}[r]^{e[2]} & 
\overset{e[3]}{\circ}\ar @{-}[d]^{e[4]} \\
\underset{d[3]}{\circ} \ar@{-}[r]_{d[4]} & \underset{d[5]}{\circ}}
\end{align*}
and checking that the new couple of lines thus generated are in $\hlist$:
 \begin{align*}
&\linethree{d[3]}{d[2]}{d[1]}{e[2]}{e[3]} & &\linethree{d[3]}{d[4]}{d[5]}{e[4]}{e[3]}
\end{align*}

We collect all squares in a  list \texttt{allsq} of 6169 diagrams which reduces to zero after applying the criteria. Namely 6002 diagrams are not of Cartan type, only 12 remain after Criterium \ref{crit-A}, this reduces to 4 after applying this criterium to reflections. No diagrams are left after applying Criterium \ref{crit-D}.

\medbreak

Next we deal with tadpoles. For the construction, we paste a triangle with a line sharing two vertices and the corresponding edge, and check that the new line generated by this process is in $\hlist$ of rank 3. This gives a list of 1964 diagrams. We check that all tadpoles from $\hlist$ are included.
After Criterium \ref{crit-A}, 738 diagrams remain, which reduce to 344 diagrams by applying this Criterium to their reflections. We now remove diagrams from $\hlist$ and we are left with 16 diagrams, all of them are discarded by applying Criterium \ref{crit-D}.

\medspace

For tripods, we paste two lines sharing two vertices and the corresponding edge, and again check that the new line is in $\hlist$ of rank 3.
We obtain 32481diagrams, which we check that include all tripods in $\hlist$. 
We remove those that reflect (at vertex 3) into tadpoles/squares, already analysed: we are left with 1183 diagrams. 
Next we keep those which are not of Cartan type, and $\hlist$ reduces to 964 diagrams. 
We remove those which after reflection $\rho_1$, $\rho_2$ or $\rho_4$ followed by $\rho_3$, turn into tadpoles/squares: we have only 44 diagrams left.
Now we apply Criterium \ref{crit-A} and 2 diagrams survive, but they fall down on Criterium \ref{crit-D}.

\medspace

Finally lines are constructed by pasting two lines: these are 15431. 
Again, we check that they include all lines in $\hlist$ and we remove them, keeping 13520 diagrams.
We exclude succesively remove: those of Cartan type (13374 diagrams survive), lines that reflect into tadpoles after $\rho_2$ or $\rho_3$ (796),
lines that fall down on Criterium \ref{crit-A} (84) and the same Criterium after a reflection (76).
The remaining lines are discarded after Criterium \ref{crit-D} (56) and Criterium \ref{crit-E}.

\subsection{Higher ranks: 5, 6 \& 7}

The strategy for this part is analogous to the one in rank four above; thus we choose not to overload the text with as many details. We construct all possible diagrams by pasting diagrams of lesser rank. We check that this includes all diagrams in $\hlist$.
We build the criteria based on \S \ref{sec:criteria}. We discard non-admissible shapes and show that the only diagrams that survive the criteria are already in $\hlist$.

\subsubsection{rank4.g, rank5.g, rank6.g, rank7.g}\label{gap:rankn} We list all the rank $n=4,5,6,7$ diagrams in $\hlist$. 
We include permutations. Generic diagrams are evaluated in $\G_f$.
These files only depend on the files \texttt{basic.g}. 

We remark that:
\begin{itemize}[leftmargin=*]
\item Lines are vectors $[d[1],\dots,d[2n-1]]$ of length $2n-1$ representing:
\begin{align*}
&\xymatrix@C-10pt@R-10pt{ \underset{d[1]}{\circ} \ar @{-}[rr]^{d[2]} & & \underset{d[3]}{\circ} \ar @{-}[rr]^{d[4]} & & \underset{d[5]}{\circ} \ar @{.}[rr]  & & \underset{d[2n-3]}{\circ} \ar @{-}[rr]^{d[2n-2]} & & \underset{d[2n-1]}{\circ}}
\end{align*}
\item $D_n$-like diagrams are vectors $[d[1],\dots, d[2n-1]]$ of length $2n-1$:
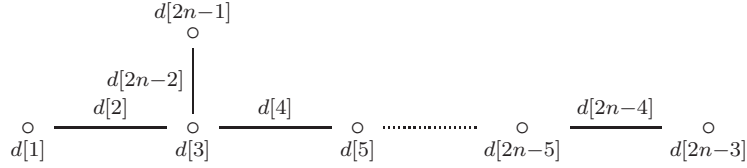
\begin{align*}
&\xymatrix@C-10pt@R-10pt{ & & \overset{d[2n-1]}{\circ} \ar @{-}[d]_{d[2n-2]}  & & & & & & \\
\underset{d[1]}{\circ} \ar @{-}[rr]^{d[2]} & & \underset{d[3]}{\circ} \ar @{-}[rr]^{d[4]} & & \underset{d[5]}{\circ} \ar @{.}[rr]  & & \underset{d[2n-5]}{\circ} \ar @{-}[rr]^{d[2n-4]} & & \underset{d[2n-3]}{\circ}}
\end{align*}
\item Tadpoles are vectors of length $2n$: 
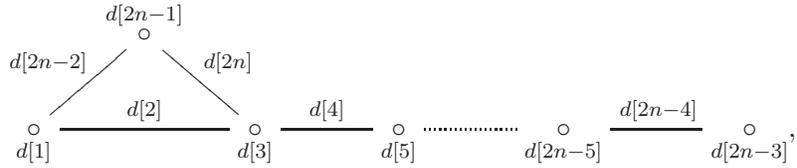
\begin{align*}
&\xymatrix@C-10pt@R-10pt{ &  \overset{d[2n-1]}{\circ} \ar @{-}[dl]_{d[2n-2]} \ar @{-}[rd]^{d[2n]} & & & & & &  & \\
\underset{d[1]}{\circ} \ar @{-}[rr]^{d[2]} & & \underset{d[3]}{\circ} \ar @{-}[rr]^{d[4]} & & \underset{d[5]}{\circ} \ar @{.}[rr]  & & \underset{d[2n-5]}{\circ} \ar @{-}[rr]^{d[2n-4]} & & \underset{d[2n-3]}{\circ},}
\end{align*}
\item Diagrams with a triangle in the middle ($n\geq 5$) are vectors of length $2n$:
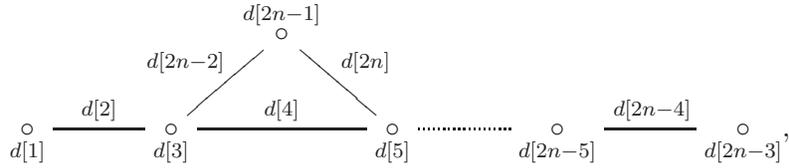
\begin{align*}
&\xymatrix@C-10pt@R-10pt{ & & & \overset{d[2n-1]}{\circ} \ar @{-}[dl]_{d[2n-2]} \ar @{-}[rd]^{d[2n]} & & & &  & \\
\underset{d[1]}{\circ} \ar @{-}[rr]^{d[2]} & & \underset{d[3]}{\circ} \ar @{-}[rr]^{d[4]} & & \underset{d[5]}{\circ} \ar @{.}[rr]  & & \underset{d[2n-5]}{\circ} \ar @{-}[rr]^{d[2n-4]} & & \underset{d[2n-3]}{\circ},}
\end{align*}
\item Diagrams for $E_n$, $n=6,7$, are vectors of length $2n-1$:
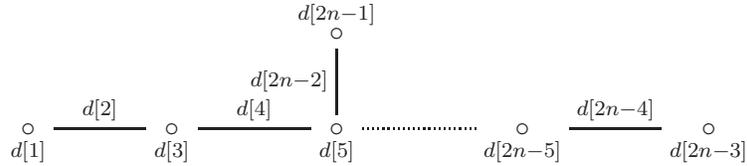
\begin{align*}
&\xymatrix@C-10pt@R-10pt{ & & & & \overset{d[2n-1]}{\circ} \ar @{-}[d]_{d[2n-2]}  & & & & \\
\underset{d[1]}{\circ} \ar @{-}[rr]^{d[2]} & & \underset{d[3]}{\circ} \ar @{-}[rr]^{d[4]} & & \underset{d[5]}{\circ} \ar @{.}[rr]  & & \underset{d[2n-5]}{\circ} \ar @{-}[rr]^{d[2n-4]} & & \underset{d[2n-3]}{\circ}}
\end{align*}

\end{itemize}
The rank 4 diagrams are imported from \texttt{rank4t.g}, and adjusted to these new shapes. For instance, for lines:
\begin{lstlisting}
gap> Read("/home/.../rank4t.g");
gap> Hli4o:=Hli4;;
gap> Hli4:=[];;
gap> for d in Hli4o do
> Add(Hli4,[d[1],d[2],d[3],d[4],d[5],d[8],d[7]]);
> od;
\end{lstlisting}
For $n\geq 5$, parametric diagrams are written as generating vectors, then later evaluated in $\G_f$.
Let us fix $n=5$ for example, cases 6 and 7 are completely analogous. We write down parametric lines generators
\begin{lstlisting}
gap> li5ag;
function( e1, e2, e3, e4, e5, h ) ... end
gap> li5bg;
function( e2, e3, e4, e5, h ) ... end
gap> li5cg;
function( e1, e2, e3, e4, e5 ) ... end
gap> li5dg;
function( e2, e3, e4, e5, h ) ... end
\end{lstlisting}
corresponding to the parametric lines in 
\begin{itemize}
\item[\texttt{a}.] Rows 1 and 2,
\item[\texttt{b}.] Rows 3 and 4,
\item[\texttt{c}.] Rows 5 and 6,
\item[\texttt{d}.] Rows 7, 9 and 10,
\end{itemize}
respectively. We evaluate these functions on $e_i\in \{0,1\}$ and $h\in\G_f$, creating lists: 
\begin{lstlisting}
gap>li5qa;;
gap>li5qb;;
gap>li5qc;;
gap>li5qd;;
\end{lstlisting}
We collect all resulting evaluated parametric lines in a list
\begin{lstlisting}
gap>li5q:=Zet(Concatenation(li5qa,li5qb,li5qc,li5qd));;
\end{lstlisting}
Then, we write down finite lines on a list
\begin{lstlisting}
gap>li5f;;
\end{lstlisting}
and finally we consider the union of both sets, adding the permutations.
\begin{lstlisting}
gap> Hli5s:=Zet(Concatenation(li5q,li5f));;
gap> Hli5p:=[];;
gap> for d in Hli5s do
> Add(Hli5p,Permuted(d,(1,9)(2,8)(3,7)(4,6)));
> od;
gap> Hli5:=Zet(Concatenation(Hli5s,Hli5p));;
\end{lstlisting}

A similar approach is followed for $D_n$-like parametric diagrams in $\hlist$. For example, we write down parametric $D_5$-type generators
\begin{lstlisting}
gap> d5g;                      
function( e2, e3, e4, h ) ... end
\end{lstlisting}
corresponding to the diagrams of this shape in Rows 8 and 10
and collect all parametric and finite diagrams in lists
\begin{lstlisting}
gap> d5listq;;
gap> d5listf;;
\end{lstlisting}
Finally, we get all diagrams in a single list (including permutations):
\begin{lstlisting}
gap> Hd5s:=Zet(Concatenation(d5listf,d5listq));;
gap> Hd5p:=[];;
gap> for d in Hd5s do
> Add(Hd5p,Permuted(d,(1,9)(2,8)));
> od;
gap> Hd5:=Zet(Concatenation(Hd5s,Hd5p));;
\end{lstlisting}
Generic tadpoles are loaded in an analogous way:
\begin{lstlisting}
gap> Hta5s:=Zet(Concatenation(Hta5q,Hta5f));;
gap> Hta5p:=[];;
gap> for d in Hta5s do
> Add(Hta5p,Permuted(d,(1,9)(2,10)));
> od;
gap> Hta5:=Zet(Concatenation(Hta5s,Hta5p));;
\end{lstlisting}
Same applies for diagrams of type $E_6$ and $E_7$.

\subsubsection{criteria5.g, criteria6.g, criteria7.g}\label{gap:criterian}

In these files we generate $\hlist$ of all possible diagrams of rank $n=5,6,7$. For each admissible shape as in \eqref{eq:rank5-possible-shapes}, \eqref{eqn:shapes6},
\eqref{eqn:shapes7} respectively, we build the diagram by pasting appropriate diagrams in $\hlist$ of rank $n-1=4,5,6$. 
We check that all diagrams from $\hlist$, which we have collected in \texttt{rank}$n$\texttt{.g}, are included in this new list, and remove them.
Then we apply different criteria to the remaining \emph{bad diagrams} until we are able to eliminate them all.

\medbreak

Diagrams of rank 5 with a triangle in the middle are obtained by pasting two tadpoles of rank 4 sharing the triangle, and keeping those for which the new line belongs to $\hlist$ of rank 4:
\begin{lstlisting}
tm5:=[];;
for d in Hta4 do
 for e in Hta4 do
  if d[1]=e[3] and d[2]=e[2] and d[3]=e[1] and d[6]=e[8]
   and d[7]=e[7] and d[8]=e[6]
   and [e[5],e[4],d[1],d[2],d[3],d[4],d[5]] in Hli4
  then Add(tm5,
   [e[5],e[4],d[1],d[2],d[3],d[4],d[5],d[6],d[7],d[8]]);
  fi;
 od;
od;     
\end{lstlisting}
These are 194 diagrams: we check that $\hlist$ contains all rank 5 diagrams of this shape in $\hlist$ and remove them (168 \emph{bad} diagrams left).

In turn rank 5 tadpoles are constructed by pasting a tadpole of rank 4 and a line, in such a way that they share two vertices and the corresponding edge, and the new rank four line thus obtained is in $\hlist$. These are 728 diagrams that include all tadpoles in  $\hlist$: we remove them and remain with 464.

For $D_5$ type diagrams, we paste two lines and check that the $D_4$ tripod is in $\hlist$. We keep 710 bad diagrams from the 1035 diagrams obtained, which include all rank five tripods from  $\hlist$.

Lines of rank 5 are obtained by pasting two lines of rank 4. These are 6115. We check that they contain rank 5 lines from $\hlist$ and remove them: 3136 lines remain.

For ranks 6 and 7, diagrams with a triangle in the middle, tadpoles, $D_n$ type diagrams and lines are constructed in the same way. In turn, diagrams of the shape $E_6$ are constructed by pasting two $D_5$ type diagrams and checking that the new line is in $\hlist$. For $E_7$, we paste an $E_6$ diagram with a line and check if the $D_6$ diagram generated is in $\hlist$.

We get, for ranks 5, 6 and 7 respectively (where ``triangle'' stands for triangle in the middle):
\begin{table}[H]
{\begin{tabular}{|l|c|c||c|c||c|c|}
\hline
diagram & total & bad& total & bad& total & bad\\
\hline
tadpoles  & 728 & 464 & 1048 & 560 & 1948 & 988\\
\hline
tripods  & 1035 & 710 & 1237 & 716 & 2067 & 1098 \\
\hline
lines  & 6115 & 3136 & 8577 & 2776 & 16617 & 5088\\
\hline
triangles & 194 &  168 & 42 & 32 &  20  & 18 \\
\hline
$E_\ast$ & - & - & 263 & 196 & 121 & 82 \\
\hline
\end{tabular}}
\end{table}

Now we apply the criteria to discard all remaining diagrams.

\medbreak

For rank 5 tadpoles, we apply Criterium \ref{crit-A} and keep 96, which are actually 48 up to permutations. We remove those that reflect via $\rho_2$ to a triangle in the middle that does not pass the Criterium \ref{crit-A} for this shape, and get 10 diagrams left. Next, we remove those that reflect via $\rho_2$ to a $D_5$ diagram whose reflection $\rho_1$ does not pass the Criterium \ref{crit-A}: 8 diagrams are left. Now we sequentially check if reflections $\rho_3$ and $\rho_4$ of these diagrams had been already removed. This is indeed the case, as no diagrams are left. 

Next we look at diagrams with a triangle in the middle. We apply Criterium \ref{crit-A} and remove permutations. Two diagrams remain, which are discarded by using $\rho_2$: we obtain tadpoles (just analyzed) not in $\hlist$. 

For $D_5$, 48 diagrams remain after Criterium \ref{crit-A}. Next, we check if $\rho_2$ creates a triangle (tadpole or in the middle) and use previous criteria: 10 diagrams are left after removing permutations. Next we apply the criteria to reflections $\rho_1$ and $\rho_5$ when the vertex is -1. Four diagrams are left: two of them reflect via $\rho_2$ to diagrams with a triangle in the middle that are not in $\hlist$ and the other two are discarded using $\rho_5$.

As for lines, 1482 remain after removing permutations and lines of Cartan type. Furthermore, 603 still remain if we remove those that reflect into a diagram with a triangle. And 106 are left after Criterium \ref{crit-A}. We apply now this criteria to the reflections $\rho_1,\dots,\rho_4$ at vertices labeled with $-1$ and to $\rho_5$ (any label): we get 36 lines, that become 2 after double reflections and zero after triple reflections. 

\medbreak

On rank 6, we are left with 10 diagrams with a triangle in the middle after Criterium \ref{crit-A}. Then we apply Criterium \ref{crit-C}, with reflections, which gives a diagram of type $D_5$ in all cases. No diagrams survive. 

We have 16 tadpoles after Criterium \ref{crit-A} and 8 after applying this criterium to reflections. Via $\rho_2$, we obtain diagrams with a triangle in the middle which are not in $\hlist$.

For $E_6$, we combine Criterium \ref{crit-A} with one that discard diagrams that reflect to one containing a triangle. This leads to 28 diagrams. We apply this combination of criteria to reflections and obtain two diagrams. For both, $\rho_3\rho_6$ gives a star, which is a forbidden shape.

The same combination of criteria leads to 53 remaining diagrams of type $D_6$. They fall down on applying this criteria to reflections. 

For lines, we remove those that reflect into a diagram with a triangle, that gives 1768 lines left, of which 1622 remain after discarding those of Cartan type. After Criterium \ref{crit-A} only 48 remain. Now we combine these criteria and apply it to reflections: we get 8 lines, that fall down on the criteria applied to double reflections.

\medbreak

For rank 7, we have 8 diagrams with a triangle in the middle after Criterium \ref{crit-A}. No diagrams are left after applying Criterium \ref{crit-F} for $\alpha_1$,  $\alpha_2+\alpha_3$, $\alpha_3+\alpha_4$, $\alpha_4+\alpha_5$, $\alpha_6$, $\alpha_7$, including suitable reflections: they turn into $E_6$ diagrams that are not in $\hlist$.

Criterium \ref{crit-A} for tadpoles leads to 8 diagrams. Criterium \ref{crit-F} for $\alpha_1+\alpha_2$, $\alpha_2+\alpha_3$, $\alpha_3+\alpha_4$, $\alpha_5$, $\alpha_6$, $\alpha_7$ give rank 6 lines which are not in $\hlist$.

For $E_7$, we are left 12 diagrams after Criterium \ref{crit-A} and discarding those that reflect into a diagram with a triangle. After multiple reflections, no diagrams survive. This combination of criteria leave 10  $D_7$ diagrams standing and again fall down on the same criteria applied to reflections. 

We have 3332 lines by removing those that reflect into diagrams with a triangle and 3186 which are not of Cartan type. After Criterium \ref{crit-A} only 8 remain.  Applying all of these criteria to reflection removes all lines.

\subsection{Forbidden shapes}

We include three extra files in which we deal with some forbidden shapes: vertices of valence four (or greater), star-shaped diagrams, and rank 7 diagrams, see \S \ref{subsec:highdegree}, \S \ref{subsec:forbidden6} and \S \ref{subsec:forbidden7} respectively.

\subsubsection{valence4.g}\label{gap:valence4}

Recall from \S \ref{subsec:highdegree} that we are necessarily are bounded to the following shapes, to avoid cycles: \emph{bowtie} diagrams $\bowtie$, \emph{semidirect product} diagrams $\rtimes$  and \emph{cross product} diagrams $\times$.

First we build all 108 diagrams $\bowtie$ by pasting two tadpoles of rank four and checking that the new tadpoles are in $\hlist$. Criterium \ref{crit-A} for $\alpha_4+\alpha_5$ and $\alpha_1+\alpha_2$ leaves no diagrams left. 

Next we turn to 92 diagrams $\rtimes$ by pasting a tadpole with a tripod and checking that the new tadpole and tripod are in $\hlist$. We apply Criterium \ref{crit-A} for $\alpha_1+\alpha_2$ (which gives a rank four tripod), and  both $\alpha_3+\alpha_4$ and $\alpha_3+\alpha_5$ (which gives rank four tadpoles). No diagrams are left. 

Finally, we build 294  diagrams $\times$ by pasting two tripods and checking that the new two tripods are in $\hlist$ (and are not both of Cartan type). Now Criterium \ref{crit-A} for $\alpha_3+\alpha_5$, $\alpha_3+\alpha_4$, $\alpha_2+\alpha_3$ and $\alpha_1+\alpha_3$. We discard those that reflect at vertex three, when this is $-1$, to a diagram that contains a triangle. We have now 24 diagrams.  We check if the reflection at outer vertices with $-1$ give again a diagram in $\hlist$ and and no diagrams survive.

\subsubsection{star5.g}\label{gap:star5}

We show that there cannot be a star-shaped diagram as in:
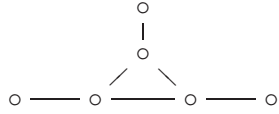
\begin{align*}
\begin{aligned}
\xymatrix@C-18pt@R-18pt{& & & {\circ}  \ar @{-}[d]   &  &  & \\
& & & {\circ} \ar @{-}[ld] \ar @{-}[rd]   &  &  & \\
{\circ} \ar @{-}[rr] & &{\circ} \ar @{-}[rr] & & {\circ} \ar @{-}[rr] & & {\circ}  }
\end{aligned}
\end{align*}
We construct all possible such diagrams by appropriately pasting two diagrams of rank 5 with a triangle in the middle and checking that the new diagram with that shape is $\hlist$. These are 8 diagrams.  We check that all of them the vertices in the triangle are all -1. We see that the first and fourth diagrams coincide up to a choice of a root of 1. It falls down on Criterium \ref{crit-F} for $\alpha_1$, $\alpha_2$, $\alpha_3+\alpha_4$, $\alpha_3+\alpha_5$, $\alpha_5+\alpha_6$ which gives a vertex (corresponding to $\alpha_2$) of valence 4 (already discarded).
The remaining 6 diagrams are actually a single one up to permutations and choices of roots. The second diagram, for instance, falls down on reflecting as a vertex labeled with $-1$, as it gives a new star whose triangle contains a vertex not labeled with $-1$. 

\subsubsection{forbidden7.g}\label{gap:forbidden7}

We build all extended $E_6$ diagrams by pasting two finite diagrams of type $E_6$ and checking that the new $E_6$ diagram is in $\hlist$, as explained in the proof of Lemma \ref{lem:rk7-forbidden}. We have the four diagrams on the left in \eqref{eqn:forbidden-diagr}, with $\zeta\in\G_3'\cup\G_4'$. Now, the reflection $\rho_5$ creates a diagram not in $\hlist$ (vertices 4 and 5 with label -1):
\begin{align}\label{eqn:forbidden-diagr}
&  \xymatrix@C-18pt@R-18pt{ &  & &  & \overset{\zeta}{\circ}  \ar @{-}[d]^{\overline{\zeta}}  & &   & & & &  & &  & &  & \overset{\zeta}{\circ}  \ar @{-}[d]^{\overline{\zeta}}  & &   &  & \\
&  & &  & \overset{\zeta}{\circ}  \ar @{-}[d]^{\overline{\zeta}}  & &   & & \ar@{<~>}[rrr]^{\rho_5} & & & &  & &  & \overset{\zeta}{\circ}  \ar @{-}[d]^{\overline{\zeta}}  & &   &  & \\
\overset{\zeta}{\circ} \ar @{-}[rr]^{\overline{\zeta}} & &\overset{\zeta}{\circ} \ar @{-}[rr]^{\overline{\zeta}} & & \overset{\zeta}{\circ} \ar @{-}[rr]^{\overline{\zeta}} & & \overset{\zeta}{\circ}\ar @{-}[rr]^{\overline{\zeta}} & & \overset{-1}{\circ} & &  & \overset{\zeta}{\circ} \ar @{-}[rr]^{\overline{\zeta}} & &\overset{\zeta}{\circ} \ar @{-}[rr]^{\overline{\zeta}} & & \overset{\zeta}{\circ} \ar @{-}[rr]^{\overline{\zeta}} & & \overset{-1}{\circ}\ar @{-}[rr]^{\zeta} & & \overset{-1}{\circ} }  
\end{align}

In turn, we build all diagrams of rank 7 with a triangle in the middle by pasting two diagrams of rank 6, both with a triangle in the middle, and checking that the line generated in this way is in $\hlist$. After removing permutations, these are 6 diagrams. They fall down on  Criterium \ref{crit-A} for (3,4;1): it leads to $E_6$ diagrams not in $\hlist$.

\end{document}